\documentclass[numbers=enddot,12pt,final,onecolumn,notitlepage]{scrartcl}
\pdfoutput=1
\usepackage[hyperfootnotes=false]{hyperref}
\usepackage[all]{hypcap}
\usepackage{amsmath}
\usepackage{amsthm}
\usepackage{amsfonts}
\usepackage[table]{xcolor}
\usepackage{ytableau}
\usepackage{mathdots}
\usepackage{verbatim}
\usepackage{tabu}
\usepackage{algorithm,algorithmicx,algpseudocode}
\usepackage{algpseudocode}
\usepackage[margin=1in]{geometry}
\usepackage{enumitem}
\usepackage{arydshln}
\usepackage{setspace}
\usepackage{multirow}
\usepackage{tikz}
\spacing{1}
\theoremstyle{definition}
\newtheorem{defin}{Definition}[section]
\newtheorem{example}[defin]{Example}
\newtheorem{alg}[defin]{Algorithm}
\theoremstyle{plain}
\newtheorem{theorem}[defin]{Theorem}
\newtheorem{prop}[defin]{Proposition}
\newtheorem{lemma}[defin]{Lemma}
\newtheorem{cor}[defin]{Corollary}
\newtheorem{remark}[defin]{Remark}
\newcommand{\ZZ}{\mathbb{Z}}
\newcommand{\NN}{\mathbb{N}}
\newcommand{\ScrC}{\mathcal{C}}
\newcommand{\ScrB}{\mathcal{B}}
\newcommand{\Cylpar}{\operatorname{Cylpar}}
\newcommand{\SSCT}{\operatorname{SSCT}}
\newcommand{\wt}{\operatorname{wt}}
\newcommand{\Insert}{\operatorname{Insert}}
\newcommand{\OneStepMulti}{\operatorname{OneStepMulti}}
\newcommand{\FullMulti}{\operatorname{FullMulti}}
\newcommand{\FM}{\operatorname{FM}}

\newcommand{\ReverseOneStepMulti}{\operatorname{ReverseOneStepMulti}}
\newcommand{\ReverseFullMulti}{\operatorname{ReverseFullMulti}}
\newcommand{\RFM}{\operatorname{RFM}}
\newcommand{\Neg}{\operatorname{Neg}}
\newcommand{\Flip}{\operatorname{Flip}}
\newcommand{\F}{\operatorname{F}}
\newcommand{\Id}{\operatorname{I}}
\newcommand{\CRSK}{\operatorname{CRSK}}
\newcommand{\CylindricRSK}{\operatorname{CylindricRSK}}
\newcommand{\InverseCylindricRSK}{\operatorname{InverseCylindricRSK}}
\newcommand{\Arr}{\operatorname{Arr}}
\newcommand{\Cyl}{\operatorname{Cyl}}
\newcommand{\hg}{\operatorname{hg}}
\newcommand{\QuA}{\operatorname{QuA}}
\newenvironment{subtab}[1][]
{\renewcommand{\arraystretch}{.9}\begin{tabular}{c}}
{\end{tabular}\renewcommand{\arraystretch}{1}}
\setlist[itemize]{itemsep=0mm}
\setlist[enumerate]{itemsep=0mm}

\begin{document}
\vspace*{\fill}
\begin{center}
{\Large \textbf{Cylindric Young Tableaux and their Properties}}

{\large
\vspace{.7cm}
Eric Neyman

Princeton University\footnote{Member of the Princeton undergraduate class of 2019. Most of this paper was written while the author was a student at Montgomery Blair High School in Silver Spring, Maryland.}

Princeton, New Jersey, USA

\vspace{.7cm}
\textit{Mentor: Darij Grinberg}

\textit{Massachusetts Institute of Technology}

\textit{Cambridge, Massachusetts, USA}}
\end{center}
\vspace*{\fill}

\newpage
\tableofcontents

\newpage
\begin{center}
{\large \textbf{Cylindric Young Tableaux and their Properties}}
\end{center}

\begin{abstract}
\textsc{Abstract.} Cylindric Young tableaux are combinatorial objects that first appeared in the 1990s. A natural extension of the classical notion of a Young tableau, they have since been used several times, most notably by Gessel and Krattenthaler \cite{GesKra} and by Alexander Postnikov \cite[\S3]{Post}. Despite this, relatively little is known about cylindric Young tableaux. This paper is an investigation of the properties of this object. In this paper, we extend the Robinson-Schensted-Knuth correspondence, a well-known and very useful bijection concerning regular Young tableaux, to be a correspondence between pairs of cylindric tableaux. We use this correspondence to reach further results about cylindric tableaux. We then establish an interpretation of cylindric tableaux in terms of a game involving marble-passing. Next, we demonstrate a generic method to use results concerning cylindric tableaux in order to prove results about skew Young tableaux. We finish with a note on Knuth equivalence and its analog for cylindric tableaux.
\end{abstract}


\newpage\section{Introduction}
In 1997, Gessel and Krattenthaler introduced cylindric semistandard Young tableaux as a modification of the concept of semistandard Young tableaux \cite{GesKra}. They are essentially semistandard Young tableaux on a cylinder (the precise structure will be discussed in the following section). Hence, cylindric Young tableaux do not have a top row or a bottom row and, in general, are more ``symmetric" than regular Young tableaux.

Being a natural extension of the notion of Young tableaux, cylindric tableaux have been studied by mathematicians during the last two decades. Among these are Alexander Postnikov \cite[\S3]{Post} and Peter McNamara \cite{McN}, as well as Jennifer Morse and Anne Schilling \cite[\S3.1]{MorSch}. The purpose of this paper is to expand this field through new results and applications.

In this paper, we begin by defining multi-insertion and reverse multi-insertion for cylindric tableaux, processes analogous to row-insertion and row-deletion for regular tableaux. We then prove a cylindric analog of the row-bumping lemma, a useful fact in tableau theory. We proceed to describe and prove an analog of the Robinson-Schensted-Knuth (RSK) correspondence for cylindric tableaux. This is the principal result of the paper, as the RSK correspondence for regular tableaux, which is a bijection between pairs of tableaux and matrices, has a variety of combinatorial applications \cite[\S4.3]{Ful}. We adapt one of these applications --- the Cauchy identity --- to cylindric Schur polynomials (defined in the following section). We proceed to prove a surprising symmetry property of our cylindric RSK correspondence.

We then demonstrate an interpretation of cylindric tableaux in terms of people passing marbles in circles. We note a possible application of this interpretation to information theory. Next, we show how results concerning cylindric tableaux may be used to prove analogous results concerning regular tableaux (including skew Young tableaux).

Finally, we define a natural variation on Knuth equivalence for words that represent cylindric tableaux, which we call cyclic Knuth equivalence. We demonstrate that, despite the naturality of this extension, cyclic Knuth equivalence is not a useful construct for cylindric tableaux, as all tableaux of the same weight are cyclic Knuth equivalent under our definition. We believe that finding a variation on Knuth equivalence that does not place all cylindric tableaux of the same content in the same equivalence class would be an important development in cylindric tableau theory.

\section{Preliminary Definitions}
\begin{defin}
Fix positive integers $k$ and $n$, with $n > k$, for the rest of this paper, excluding the examples given in the paper. The \textit{cylinder} $\ScrC_{k,n}$ is the quotient $\ZZ^2/(-k,n-k)\ZZ$. In other words, $\ScrC_{k,n}$ is the set of equivalence classes of points modulo a shift by the $\ZZ$ vector $(-k,n - k)$ \cite[\S3]{Post}.\footnote{The fact that the shift is described as $-k$ and $n - k$, instead of $-k$ and $m$ for some $m$, is a standard in cylindric tableau theory \cite[\S3]{Post}; the reason for this is not in the scope of this paper.}
\end{defin}
\begin{defin}
A \textit{cylindric partition} $\lambda$ on $\ScrC_{k,n}$ is a weakly decreasing sequence of integers $ \dots, \lambda_{-1}, \lambda_{0}, \lambda_{1}, \dots$, infinite in both directions, such that for any integer $m$, $\lambda_m = \lambda_{m+k} + n - k$.
\end{defin}
In this paper, we will think of cylindric partitions in terms of their pictures in the plane --- analogues of Young diagrams. Drawn below is a cylindric partition; here, $k = 7$ and $n = 12$, and the part of the sequence of the partition that is shown below is $6$, $4$, $2$, $2$, $2$, $2$, $1$, $1$, $-1$, $-3$, $-3$, $-3$, $-3$, $-4$.
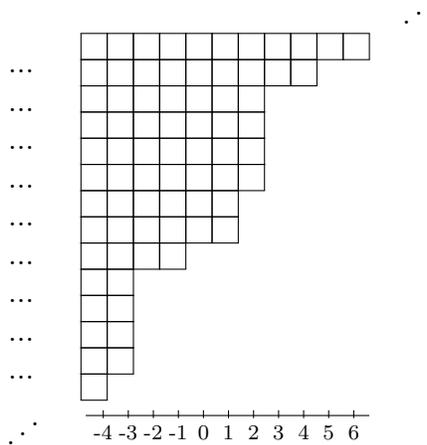
\begin{center}
\def\arraystretch{0.5}
\begin{tabular}{rcl}
&&$\iddots$\\
\multirow{3}{1em}{... ... ... ... ... ... ... ... ...}&
\ytableausetup{smalltableaux}
\ydiagram{11,9,7,7,7,7,6,6,4,2,2,2,2,1}&\\
$\iddots$&\begin{tikzpicture}
\hspace{.1cm}
\begin{scriptsize}
\draw (0.1,0) -- (3.87,0);
\foreach \x  in {1,...,11}
\draw[xshift=.333 *\x cm] (0pt,2pt) -- (0pt,-1pt) node[below,fill=white] {\the\numexpr\x-5\relax};
\end{scriptsize}
\end{tikzpicture}&
\end{tabular}
\end{center}
\begin{defin}
The set of cylindric partitions  will be denoted $\Cylpar$.
\end{defin}
In this paper, the term ``partition" will be used to refer to cylindric partitions, unless stated otherwise. The term ``regular" will be used as the negation of ``cylindric" (e.g. non-cylindric partitions will be called regular partitions).
\begin{defin}
A \textit{point} is a pair of integers $(x,y)$. In the diagrams in this paper, a point $(x,y)$ is represented on the Cartesian plane by a square of side length $1$. Note that the positive $x$-axis points downward and the positive $y$-axis points to the right, as usual in the theory of Young diagrams. We say that a point $P = (x,y)$ is in a partition $\lambda$ (denoted $P \in \lambda$) if $y \le \lambda_x$. Visually, $P$ is in $\lambda$ if it lies inside (i.e. to the left of the right boundary of) $\lambda$ when $\lambda$ is drawn on the plane as above.
\end{defin}
\begin{defin}
A \textit{plane row} is a set of all points with the same $x$-coordinate. We say that the point $(x,y)$ is in plane row $x$. A \textit{row} is the projection of a plane row onto $\ScrC_{k,n}$. Thus, row $x$ (i.e. the projection of plane row $x$ onto the cylinder) is the same as row $x + mk$ for any integer $m$.
\end{defin}
\begin{defin}
A \textit{plane column} is a set of all points with the same $y$-coordinate. We say that the point $(x,y)$ is in plane column $y$. A \textit{column} is the projection of a plane column onto $\ScrC_{k,n}$. Thus, column $y$ (i.e. the projection of plane column $y$ onto the cylinder) is the same as column $y + m(n - k)$ for any integer $m$.
\end{defin}
\begin{defin}
A \textit{box} is the projection of a point $(x,y)$ onto $\ScrC_{k,n}$. Thus, the box with coordinates $(x,y)$ is the same as the box with coordinates $(x - mk,y + m(n - k))$, for any $m \in \ZZ$.
\end{defin}
\begin{defin}
For any point $P$, $\pi(P)$ is the projection of $P$ onto $\ScrC_{k,n}$. For any box $B$ in $\ScrC_{k,n}$, $\pi^{-1}(B)$ is the set of all points $P$ such that $\pi(P) = B$. For any plane row $r$, $\pi(r)$ is the projection of $r$ onto $\ScrC_{k,n}$. For any row $s$ in $\ScrC_{k,n}$, $\pi^{-1}(s)$ is the set of all plane rows $r$ such that $\pi(r) = s$.
\end{defin}
\begin{defin} \label{pbrdef}
For any box $B$ and plane row $r$ such that $B \in \pi(r)$, we define $\pi^{-1}_r(B)$ to be the point in $\pi^{-1}(B)$ that is in plane row $r$.
\end{defin}
\begin{defin}
A box $B$ is in a partition $\lambda$ (denoted $B \in \lambda$) if $P \in \lambda$ for some $P \in \pi^{-1}(B)$. Note that, because of the periodicity of cylindric partitions, if $P \in \lambda$ for some $P \in \pi^{-1}(B)$, then $P \in \lambda$ for all $P \in \pi^{-1}(B)$.
\end{defin}
\begin{defin}
Let $\lambda$ and $\mu$ be two cylindric partitions. We say that $\mu \subseteq \lambda$ if, for all integers $m$, we have $\mu_m \le \lambda_m$.
\end{defin}
\begin{defin}
For any cylindric partitions $\lambda$ and $\mu$ such that $\mu \subseteq \lambda$, a box $B$ is in $\lambda / \mu$ (denoted $B \in \lambda / \mu$) if $B \in \lambda$, but $B \not \in \mu$. (Note that the set of boxes in $\lambda / \mu$ is finite, since there are $k$ rows and in each row the $i$'th row has
$\lambda_i - \mu_i$ boxes in $\lambda / \mu$.) Similarly, point $P$ is in $\lambda / \mu$ (also denoted $P \in \lambda / \mu$) if $P \in \lambda$, but $P \not \in \mu$.
\end{defin}
\begin{defin}
Given two cylindric partitions $\lambda$ and $\mu$ such that $\mu \subseteq \lambda$, a \textit{semistandard cylindric tableau} with \textit{outer shape} $\lambda$ and \textit{inner shape} $\mu$ is a map $R$ from the set of all boxes in $\lambda / \mu$ to a totally ordered set $A$ such that
\begin{enumerate}[label=(\alph*)]
\item $R(\pi((x,y_1))) \le R(\pi((x,y_2)))$ for any $x$, $y_1$, and $y_2$ such that $(x,y_1)$ and $(x,y_2)$ are in $\lambda / \mu$, and $y_1 < y_2$; and
\item $R(\pi((x_1,y))) < R(\pi((x_2,y)))$ for any $x_1$, $x_2$, and $y$ such that $(x_1,y)$ and $(x_2,y)$ are in $\lambda / \mu$, and $x_1 < x_2$.
\end{enumerate}
The inner and outer shapes of a semistandard cylindric tableau are regarded as part of the tableau's data; that is, the tableau ``remembers" its $\lambda$ and $\mu$. A semistandard cylindric tableau can be drawn on the plane, with each square $(x,y)$ holding the entry that the corresponding box $\pi((x,y))$ maps to under the tableau (see the diagram below). Visually, a semistandard cylindric tableau's entries increase weakly from left to right along its rows and increase strictly from top to bottom along its columns. We say that $R$ is \textit{bounded} by $\lambda$ and $\mu$, and that the \textit{shape} of $R$ is $\lambda / \mu$.\footnote{While it is often convenient to think of a tableau's shape as the set of boxes between two partitions (hence the notation), $\lambda / \mu$ here is a notation for the pair $(\lambda, \mu)$; the pair of partitions may carry more information than simply the set of boxes between them.} We call $A$ the \textit{alphabet} of $R$. Frequently, alphabets of semistandard cylindric tableaux are $\ZZ^+$ or $\ZZ$. For the rest of this paper, the alphabets of all tableaux will be implicit. The elements of an alphabet are referred to as \textit{letters}. For any particular semistandard cylindric tableau, the image of a box under the tableau is referred to as the \textit{entry} in the box; the images of the boxes of $\lambda / \mu$ under the tableau are collectively called the \textit{entries} of the tableau.
\end{defin}
Note that, if $\mu \not \subseteq \lambda$, then there are no semistandard cylindric tableaux of shape $\lambda / \mu$.

For the rest of this paper, unless stated otherwise, the word ``tableau" will be used to refer to semistandard cylindric tableaux.

Below is an example of a cylindric tableau. In this example, $k = 3$ (i.e. the vertical period of the tableau is $3$) and $n = 6$ (which means that the horizontal period of the tableau is $6 - 3 = 3$). The entry in a given box of a cylindric tableau in the diagram is the image of the box under the tableau map. In such diagrams, the top row that is drawn is row $0$ of the tableau.
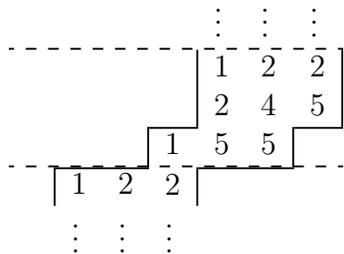
\begin{center}
\begin{tabular}{m{.2cm}m{.2cm}m{.2cm}m{.2cm}m{.2cm}m{.2cm}m{.2cm}}
&&&&$\vdots$&$\vdots$&$\vdots$\\
\hdashline&&&\multicolumn{1}{c|}{}&1&2&\multicolumn{1}{c|}{2}\\
&&&\multicolumn{1}{c|}{}&2&4&\multicolumn{1}{c|}{5}\\
\cline{4-4}\cline{7-7}&&\multicolumn{1}{c|}{}&1&5&\multicolumn{1}{c|}{5}\\
\hdashline\cline{2-3}\cline{5-6}\multicolumn{1}{c|}{}&1&2&\multicolumn{1}{c|}{2}\\
&$\vdots$&$\vdots$&$\vdots$
\end{tabular}
\end{center}
Notice that the top row that we draw is repeated at the bottom of the diagram. In this paper, we will draw tableaux as they are drawn above: the top row will be repeated at the bottom, separated by a dashed line. This convention allows us to distinguish, for example, between the following two tableaux, the first of which has $k = 1$ and the second of which has $k = 2$.
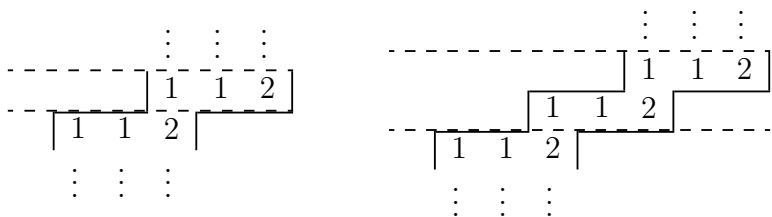
\begin{center}
\begin{tabular}{m{.2cm}m{.2cm}m{.2cm}m{.2cm}m{.2cm}m{.2cm}}
&&&$\vdots$&$\vdots$&$\vdots$\\
\hdashline&&\multicolumn{1}{c|}{}&1&1&\multicolumn{1}{c|}{2}\\
\hdashline\cline{2-3}\cline{5-6}\multicolumn{1}{c|}{}&1&1&\multicolumn{1}{c|}{2}\\
&$\vdots$&$\vdots$&$\vdots$
\end{tabular}
\hspace{1cm}
\begin{tabular}{m{.2cm}m{.2cm}m{.2cm}m{.2cm}m{.2cm}m{.2cm}m{.2cm}m{.2cm}}
&&&&&$\vdots$&$\vdots$&$\vdots$\\
\hdashline&&&&\multicolumn{1}{c|}{}&1&1&\multicolumn{1}{c|}{2}\\
\cline{4-5}\cline{7-8}&&\multicolumn{1}{c|}{}&1&1&\multicolumn{1}{c|}{2}\\
\hdashline\cline{2-3}\cline{5-6}\multicolumn{1}{c|}{}&1&1&\multicolumn{1}{c|}{2}\\
&$\vdots$&$\vdots$&$\vdots$
\end{tabular}
\end{center}
\begin{defin}
The set of semistandard cylindric tableaux of shape $\lambda / \mu$, for $\lambda, \mu \in \Cylpar$, will be denoted $\SSCT(\lambda / \mu)$.
\end{defin}
\begin{defin}
Given $R \in \SSCT(\lambda / \mu)$, a box $B$ is in $R$ (denoted $B \in R$) if $B \in \lambda / \mu$.
\end{defin}
\begin{defin} \label{wtdef}
Let $R$ be a cylindric tableau with alphabet $A$. The \textit{weight} of $R$ (also known as the \textit{content} of $R$) is the map $\wt(R):A \rightarrow \NN$, where $\wt(R)(i)$ is the size of the preimage of $i$ under $R$ (that is, it is the number of boxes that, in $R$, contain $i$). Here, we use $\NN$ to represent the set of natural numbers, including $0$.

When $A = \{1,2,3, \dots\}$, the weight of $R$ is written as the sequence $(\wt(R)(1)$, $\wt(R)(2)$, $ \dots)$, and can be truncated at a place where all following terms are zero. For example, the cylindric tableau shown below has weight $(1,3,1,0,1)$.
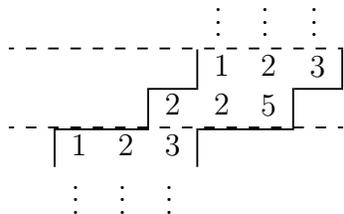
\begin{center}
\begin{tabular}{m{.2cm}m{.2cm}m{.2cm}m{.2cm}m{.2cm}m{.2cm}m{.2cm}}
&&&&$\vdots$&$\vdots$&$\vdots$\\
\hdashline&&&\multicolumn{1}{c|}{}&1&2&\multicolumn{1}{c|}{3}\\
\cline{4-4}\cline{7-7}&&\multicolumn{1}{c|}{}&2&2&\multicolumn{1}{c|}{5}\\
\hdashline\cline{2-3}\cline{5-6}\multicolumn{1}{c|}{}&1&2&\multicolumn{1}{c|}{3}\\
&$\vdots$&$\vdots$&$\vdots$
\end{tabular}
\end{center}
\end{defin}
\begin{defin}
Given an alphabet $A$, for every $a \in A$, let $x_a$ be a variable. For any two $a$, $b \in A$, $x_a$ and $x_b$ are distinct variables. We will denote the family $(x_a)_{a \in A}$ by $\mathbf{x}$; it will be called a \textit{variable set}.
\end{defin}
\begin{defin}
Given a tableau $R$, we will define the \textit{weight monomial} of $R$ with variable set $\mathbf{x}$, denoted $\mathbf{x}^{\wt(R)}$, as follows:
\begin{equation*}
\mathbf{x}^{\wt(R)} = \prod \limits_{a \in A} {x_a}^{\wt(R)(a)} = \prod_{B \text{, box in } R} x_{R(B)}.
\end{equation*}
We use $R(B)$ above to denote the entry of $R$ in $B$.

For example, if $R$ is the tableau that is shown as an example in Definition \ref{wtdef}, then $\mathbf{x}^{\wt(R)} = x_1x_2^3x_3x_5$.
\end{defin}
\begin{defin}
Given two cylindric partitions $\mu$ and $\lambda$ such that $\mu \subseteq \lambda$, the \textit{Schur polynomial} of $\lambda / \mu$ with variable set $\mathbf{x}$, denoted $s_{\lambda / \mu}(\mathbf{x})$, is a power series of bounded degree, defined as follows:
\begin{equation*}
s_{\lambda / \mu}(\mathbf{x}) = \sum \limits_{R \in \SSCT(\lambda / \mu)} \mathbf{x}^{\wt(R)}.
\end{equation*}
\end{defin}
\begin{defin}
A \textit{horizontal strip} is a pair $\lambda / \mu$ of partitions $\mu$ and $\lambda$, with $\mu \subseteq \lambda$, such that none of the boxes in $\lambda / \mu$ are in the same column. We say that a set $S$ of boxes forms a horizontal strip if there exists a horizontal strip $h$ such that $S$ is the set of boxes in $h$.\footnote{It is worth noting that $\lambda/\mu$ is a horizontal strip if and only if $\lambda_i \ge \mu_i \ge \lambda_{i+1}$.}
\end{defin}

\section{Forward Internal Insertion and Multi-Insertion}
\subsection{Forward Insertion Algorithms and Examples}
\begin{defin}
Given $R \in \SSCT(\lambda / \mu)$, a box $B = \pi((i,j))$ is an \textit{inside cocorner} of $R$ if $B \not \in \mu$, but $\pi((i-1,j)) \in \mu$ and $\pi((i,j-1)) \in \mu$. Visually, $B$ is an inside cocorner of $R$ if it is to the right of (outside) $\mu$, but the boxes to the left of and above $B$ are inside $\mu$. Note that $B$ is not necessarily a box of $R$, as it is possible that $B$ is outside $\lambda$ as well.
\end{defin}

Given a cylindric tableau $R$ and an inside cocorner $B$ of $R$, we can \textit{internally insert} $B$ into $R$, in a process that is similar to regular row-insertion, using the algorithm below, which takes a cylindric tableau and an inside cocorner of the tableau as input and outputs a new tableau (we will later prove that the output is indeed a valid semistandard tableau). This algorithm is inspired by the internal row insertion for skew tableaux described by Bruce Sagan and Richard Stanley \cite[\S2]{SagStan}.

\begin{alg}[Internal Row-Insertion] \label{intinsert}
\begin{center}\end{center}
\vspace{.5cm}
\begin{algorithmic}[1]
Function Insert(tableau $R$, box $B$) \Comment{$B$ must be an inside cocorner of $R$.}
\State $\mu :=$ inner shape of $R$.
\State $\lambda :=$ outer shape of $R$.
\If{$B \in R$}:
\State $x :=$ entry of $R$ that is in $B$.
\EndIf.
\State Expand $\mu$ to include $B$ and remove $B$ from $R$. \label{expmu1}
\If{$B \not \in \lambda$}: \Comment{This happens only when $B$ was not in $R$ to begin with.}
\State Expand $\lambda$ to include $B$. \label{explambda}
\Else:
\While{$x \neq$ null}: \label{loopstart}
\State $r:=$ row of $B$.
\If{$x$ is greater than or equal to every entry in $R$ in row $r + 1$}:
\State $B:=$ leftmost box of row $r + 1$ that is not in $\lambda$.
\State Put $x$ in $B$ and expand $\lambda$ to include $B$. \Comment{We say that $x$ \textit{lands} in $B$.}
\State $x :=$ null.
\Else:
\State $B:=$ box of the leftmost entry of $R$ in row $r + 1$ that is greater than $x$.
\State $x':=$ entry of $R$ in $B$.
\State In $R$, replace $x'$ in $B$ with $x$. \label{replacestep}
\State $x := x'$.
\EndIf.
\EndWhile.
\EndIf. 
\State \Return $R$.
\end{algorithmic}
\end{alg}

Note that, although we use the phrase, ``internally insert $B$ into $R$," the result of the process is that $B$, which was previously in $R$, is no longer in $R$.

\begin{example} \label{insertex}
Suppose we want to apply Algorithm \ref{intinsert} to $(R,B)$, where $R$ is the tableau shown below and $B$ is the box that contains $1$ in $R$. We say that we \textit{bump} the $1$ from $B$ (or that the $1$ is \textit{bumped} from $B$).
\begin{center}
\begin{tabular}{m{.2cm}m{.2cm}m{.2cm}m{.2cm}m{.2cm}}
&&&$\vdots$&$\vdots$\\
\hdashline&&\multicolumn{1}{c|}{}&1&\multicolumn{1}{c|}{4}\\
\cline{3-3}&\multicolumn{1}{c|}{}&2&5&\multicolumn{1}{c|}{6}\\
&\multicolumn{1}{c|}{}&3&7&\multicolumn{1}{c|}{7}\\
\hdashline\cline{2-2}\cline{4-5}\multicolumn{1}{c|}{}&1&\multicolumn{1}{c|}{4}\\
&$\vdots$&$\vdots$
\end{tabular}
\end{center}
We expand $\mu$ to include the box containing the $1$ and exclude the $1$ from the tableau.
\begin{center}
\begin{tabular}{m{.2cm}m{.2cm}m{.2cm}m{.2cm}}
&&&$\vdots$\\
\hdashline&&\multicolumn{1}{c|}{}&\multicolumn{1}{c|}{4}\\
\cline{2-3}\multicolumn{1}{c|}{}&2&5&\multicolumn{1}{c|}{6}\\
\multicolumn{1}{c|}{}&3&7&\multicolumn{1}{c|}{7}\\
\hdashline\cline{3-4}\multicolumn{1}{c|}{}&\multicolumn{1}{c|}{4}\\
&$\vdots$
\end{tabular}
\end{center}
We then replace the leftmost entry greater than $1$ in the following row (which is $2$) with $1$. We say that the $2$ is \textit{bumped} by the $1$.
\begin{center}
\begin{tabular}{m{.2cm}m{.2cm}m{.2cm}m{.2cm}}
&&&$\vdots$\\
\hdashline&&\multicolumn{1}{c|}{}&\multicolumn{1}{c|}{4}\\
\cline{2-3}\multicolumn{1}{c|}{}&1&5&\multicolumn{1}{c|}{6}\\
\multicolumn{1}{c|}{}&3&7&\multicolumn{1}{c|}{7}\\
\hdashline\cline{3-4}\multicolumn{1}{c|}{}&\multicolumn{1}{c|}{4}\\
&$\vdots$
\end{tabular}
\end{center}
This process continues as follows:
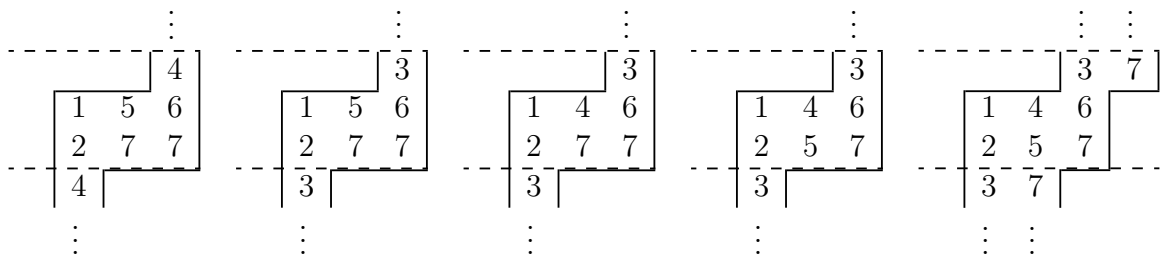
\begin{center}
\begin{tabular}{m{.2cm}m{.2cm}m{.2cm}m{.2cm}}
&&&$\vdots$\\
\hdashline&&\multicolumn{1}{c|}{}&\multicolumn{1}{c|}{4}\\
\cline{2-3}\multicolumn{1}{c|}{}&1&5&\multicolumn{1}{c|}{6}\\
\multicolumn{1}{c|}{}&2&7&\multicolumn{1}{c|}{7}\\
\hdashline\cline{3-4}\multicolumn{1}{c|}{}&\multicolumn{1}{c|}{4}\\
&$\vdots$
\end{tabular}
\hspace{.2cm}
\begin{tabular}{m{.2cm}m{.2cm}m{.2cm}m{.2cm}}
&&&$\vdots$\\
\hdashline&&\multicolumn{1}{c|}{}&\multicolumn{1}{c|}{3}\\
\cline{2-3}\multicolumn{1}{c|}{}&1&5&\multicolumn{1}{c|}{6}\\
\multicolumn{1}{c|}{}&2&7&\multicolumn{1}{c|}{7}\\
\hdashline\cline{3-4}\multicolumn{1}{c|}{}&\multicolumn{1}{c|}{3}\\
&$\vdots$
\end{tabular}
\hspace{.2cm}
\begin{tabular}{m{.2cm}m{.2cm}m{.2cm}m{.2cm}}
&&&$\vdots$\\
\hdashline&&\multicolumn{1}{c|}{}&\multicolumn{1}{c|}{3}\\
\cline{2-3}\multicolumn{1}{c|}{}&1&4&\multicolumn{1}{c|}{6}\\
\multicolumn{1}{c|}{}&2&7&\multicolumn{1}{c|}{7}\\
\hdashline\cline{3-4}\multicolumn{1}{c|}{}&\multicolumn{1}{c|}{3}\\
&$\vdots$
\end{tabular}
\hspace{.2cm}
\begin{tabular}{m{.2cm}m{.2cm}m{.2cm}m{.2cm}}
&&&$\vdots$\\
\hdashline&&\multicolumn{1}{c|}{}&\multicolumn{1}{c|}{3}\\
\cline{2-3}\multicolumn{1}{c|}{}&1&4&\multicolumn{1}{c|}{6}\\
\multicolumn{1}{c|}{}&2&5&\multicolumn{1}{c|}{7}\\
\hdashline\cline{3-4}\multicolumn{1}{c|}{}&\multicolumn{1}{c|}{3}\\
&$\vdots$
\end{tabular}
\hspace{.2cm}
\begin{tabular}{m{.2cm}m{.2cm}m{.2cm}m{.2cm}m{.2cm}}
&&&$\vdots$&$\vdots$\\
\hdashline&&\multicolumn{1}{c|}{}&3&\multicolumn{1}{c|}{7}\\
\cline{2-3}\cline{5-5}\multicolumn{1}{c|}{}&1&4&\multicolumn{1}{c|}{6}\\
\multicolumn{1}{c|}{}&2&5&\multicolumn{1}{c|}{7}\\
\hdashline\cline{4-4}\multicolumn{1}{c|}{}&3&\multicolumn{1}{c|}{7}\\
&$\vdots$&$\vdots$
\end{tabular}
\end{center}
In the final step, there is no entry greater than $7$ in the following row, and the $7$ lands at the end of the following row, to the right of the $3$. The insertion process is now complete.
\end{example}

\begin{remark} \label{boxdecrease}
Consider any box $C$ of $R$. Since on line \ref{replacestep} of Algorithm \ref{intinsert} an entry only replaces an entry greater than it, it follows that the entry in box $C$ can only decrease (or disappear) throughout the row-insertion process, and does decrease if and when an entry is bumped out of $C$ and a new entry takes its place.
\end{remark}

We now draw the final tableau above, this time with more repeating rows in the diagram.

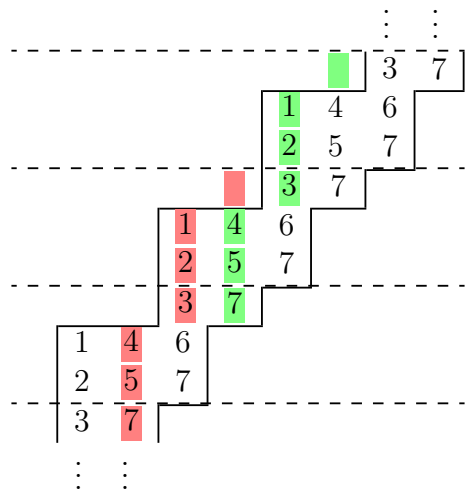
\begin{center}
\begin{tabular}{m{.2cm}m{.2cm}m{.2cm}m{.2cm}m{.2cm}m{.2cm}m{.2cm}m{.2cm}m{.2cm}}
&&&&&&&$\vdots$&$\vdots$\\
\hdashline&&&&&&\multicolumn{1}{c|}{\begin{subtab}\cellcolor{green!50}\phantom{1}\end{subtab}}&3&\multicolumn{1}{c|}{7}\\
\cline{6-7}\cline{9-9}&&&&\multicolumn{1}{c|}{}&\begin{subtab}\cellcolor{green!50}1\end{subtab}&4&\multicolumn{1}{c|}{6}\\
&&&&\multicolumn{1}{c|}{}&\begin{subtab}\cellcolor{green!50}2\end{subtab}&5&\multicolumn{1}{c|}{7}\\
\hdashline\cline{8-8}&&&&\multicolumn{1}{c|}{\begin{subtab}\cellcolor{red!50}\phantom{1}\end{subtab}}&\begin{subtab}\cellcolor{green!50}3\end{subtab}&\multicolumn{1}{c|}{7}\\
\cline{4-5}\cline{7-7}&&\multicolumn{1}{c|}{}&\begin{subtab}\cellcolor{red!50}1\end{subtab}&\begin{subtab}\cellcolor{green!50}4\end{subtab}&\multicolumn{1}{c|}{6}\\
&&\multicolumn{1}{c|}{}&\begin{subtab}\cellcolor{red!50}2\end{subtab}&\begin{subtab}\cellcolor{green!50}5\end{subtab}&\multicolumn{1}{c|}{7}\\
\hdashline\cline{6-6}&&\multicolumn{1}{c|}{}&\begin{subtab}\cellcolor{red!50}3\end{subtab}&\multicolumn{1}{c|}{\begin{subtab}\cellcolor{green!50}7\end{subtab}}\\
\cline{2-3}\cline{5-5}\multicolumn{1}{c|}{}&1&\begin{subtab}\cellcolor{red!50}4\end{subtab}&\multicolumn{1}{c|}{6}\\
\multicolumn{1}{c|}{}&2&\begin{subtab}\cellcolor{red!50}5\end{subtab}&\multicolumn{1}{c|}{7}\\
\hdashline\cline{4-4}\multicolumn{1}{c|}{}&3&\multicolumn{1}{c|}{\begin{subtab}\cellcolor{red!50}7\end{subtab}}\\
&$\vdots$&$\vdots$
\end{tabular}
\end{center}

Above in green and red are two bumping routes, which we will define shortly.

\begin{remark}
Algorithm \ref{intinsert} always ends after a finite number of steps.
\end{remark}

\begin{proof}
Let $M$ be the largest entry of $R$. Every iteration of the loop beginning on line \ref{loopstart}, $x$ increases; however, $x$ cannot be greater than $M$. Thus, the loop must terminate and the algorithm must end.
\end{proof}

\begin{defin}
Let $R$ be a tableau and $P$ be a point such that $\Insert(R, \pi(P))$ is well-defined. We will define the \textit{bumping route} of $P$ as a list of points (not boxes --- and hence the two bumping routes above are distinct) that is constructed as follows:

Suppose that we perform Algorithm \ref{intinsert} with $R$ and $\pi(P)$ as input. We will modify the algorithm for the purposes of this definition. At the very beginning, we let $s$ be the plane row of $P$, and we initialize $H$ as a list with the single entry $P = \pi^{-1}_s(\pi(B))$ (see Definition \ref{pbrdef}). At the very end of the loop starting on line \ref{loopstart}, add a line that increments $s$ by $1$.\footnote{After this augmentation, $B \in \pi(s)$.} Immediately after the new line, append $\pi^{-1}_s(B)$ to the end of $H$. The bumping route of $P$ is defined to be the list $H$ as it stands after the algorithm terminates.
\end{defin}

Intuitively, then, the bumping route is a list of points in consecutively increasing rows, such that each point corresponds to a value of $B$ in Algorithm \ref{intinsert} during its execution.

Later, we will prove two important results: that Algorithm \ref{intinsert} necessarily outputs a valid semistandard tableau and that each box in a bumping route (except the first one) is weakly left of the box before it (we say that bumping routes \textit{trend weakly left}).

It is also possible to internally row-insert multiple entries at the same time. The algorithm below is the algorithm for one-step multi-insertion. It takes a tableau and a regular insertion queue (defined below).

\begin{defin}
A \textit{queue} is a data structure that operates on a ``first-in-first-out" basis: when elements are added to a queue, they are added to the end of the queue, but when elements are removed from a queue, they are removed from the beginning of the queue. An \textit{insertion queue} is a queue of pairs, such that the first element of each pair is a letter (element of the (implicit) alphabet) and the second element of each pair is a row.
\end{defin}

\begin{defin}
An insertion queue is \textit{regular} if, for any two elements of the queue $(x_1,r)$ and $(x_2,r)$, where $x_1 < x_2$, $(x_1,r)$ comes before $(x_2,r)$ in the queue.\footnote{The fact that the second element of a pair in an insertion queue is a row should be kept in mind. For example, if $k = 3$, $(3,5)$ cannot come before $(2,2)$ in a regular insertion queue, since $2$ and $5$ would refer to the same row.}
\end{defin}

The algorithm below (the subroutine to our multi-insertion algorithm) outputs a pair, whose first element is a map from a subset of $\ScrC_{k,n}$ to the alphabet of the tableau taken as input (though the map itself is not necessarily a valid semistandard tableau), and whose second element is an insertion queue corresponding to the entries bumped from the tableau via insertion of the entries in the insertion queue that is taken as a parameter.

\begin{alg}[One-Step Multi-Insertion] \label{osmulti}
\begin{center}\end{center}
\vspace{.5cm}
\begin{algorithmic}[1]
Function OneStepMulti(tableau $R$, insertion queue $\mathfrak{q}$) \Comment{$\mathfrak{q}$ must be regular.}
\State $\lambda :=$ outer shape of $R$.
\State $\mathfrak{q}' :=$ empty insertion queue.
\While{$\mathfrak{q}$ is not empty}: \label{onesteploop}
\State Remove the first element from $\mathfrak{q}$. Let it be $(x,r)$.
\If{$x$ is greater than or equal to every entry in $R$ in row $r$}:
\State $B:=$ leftmost box of row $r$ that is not in $\lambda$.
\State Put $x$ in $B$ and expand $\lambda$ to include $B$. \Comment{$\lambda$ is not necessarily a valid partition anymore.}
\Else:
\State $B:=$ box of the leftmost entry of $R$ in row $r$ that is greater than $x$. \label{newbox}
\State $x':=$ entry of $R$ in $B$.
\State In $R$, replace $x'$ in $B$ with $x$.
\State Add $(x',r + 1)$ to $\mathfrak{q}'$. \label{putinq}
\EndIf.
\EndWhile.
\State \Return $(R,\mathfrak{q}')$. \Comment{$R$ is not necessarily a tableau.}
\end{algorithmic}
\end{alg}

\vspace{.2cm}
\begin{remark} \label{regularity}
The insertion queue returned by Algorithm \ref{osmulti} is regular.
\end{remark}

\begin{proof}
Let $(R',\mathfrak{q}') = \OneStepMulti(R,\mathfrak{q})$ for a tableau $R$ and a regular insertion queue $\mathfrak{q}$. Suppose, for contradiction, that $\mathfrak{q}'$ is not regular. Then there exist two elements of $\mathfrak{q}'$, $(y_1,r)$ and $(y_2,r)$, such that $y_1 < y_2$, but $(y_2,r)$ comes before $(y_1,r)$ in $\mathfrak{q}'$. Let $x_1$ and $x_2$ be the entries that bumped out $y_1$ and $y_2$, respectively. Then $(x_2,r-1)$ and $(x_1,r-1)$ were in $\mathfrak{q}$, with $(x_2,r-1)$ coming first. Since $\mathfrak{q}$ is regular, it follows that $x_2 < x_1$. We also know that $x_1 < y_1$, since $x_1$ bumps out $y_1$. It follows that $x_2 < x_1 < y_1 < y_2$. Since $y_1 < y_2$, $y_1$ is to the left of $y_2$ in $R$, and $y_1 > x_2$. It follows that $y_2$ is not the leftmost entry of $R$ in its row that is greater than $x_2$ at the time that $x_2$ is to be inserted, which is a contradiction, because then $x_2$ does not bump out $y_2$. Thus, $\mathfrak{q}'$ is a regular insertion queue.
\end{proof}

\begin{remark} \label{qequiv}
Let $R$ be a cylindric tableau and let $\mathfrak{q}_1$ and $\mathfrak{q}_2$ be two regular insertion queues that are permutations of each other. Let $(R_1, \mathfrak{q}'_1) = \OneStepMulti(R, \mathfrak{q}_1)$ and $(R_2, \mathfrak{q}'_2) = \OneStepMulti(R, \mathfrak{q}_2)$. Then $R_1 = R_2$ and $\mathfrak{q}'_1$ is a permutation of $\mathfrak{q}'_2$.
\end{remark}

\begin{proof}
Consider a particular row $r$ and all of the pairs of $\mathfrak{q}_1$ and $\mathfrak{q}_2$ to be placed in $r$ (in other words, all pairs whose second element is $r$). These are the same pairs, because $\mathfrak{q}_1$ is a permutation of $\mathfrak{q}_2$. Furthermore, regularity defines an ordering (least to greatest) among these pairs based on their first elements, so these pairs are also in the same order. Since insertion into row $r$ is only affected by the entry being inserted into row $r$ and the entries already in row $r$, it follows that row $r$ is the same in $R_1$ and $R_2$ after the two insertion queues have been processed. Since this is true for all $r$, it follows that $R_1 = R_2$.

Since $R_1 = R_2$, the same entries were bumped out of $R$ to produce $R_1$ and $R_2$ (though not necessarily in the same order). These bumped-out entries, along with the row numbers of the rows immediately below the ones from which they were bumped out, constitute the pairs in $\mathfrak{q}'_1$ and $\mathfrak{q}'_2$, respectively. These pairs being the same, it follows that $\mathfrak{q}'_1$ is a permutation of $\mathfrak{q}'_2$.
\end{proof}

Finally, using Algorithm \ref{osmulti}, we can internally row-insert multiple boxes into a tableau simultaneously. The following algorithm does this, taking a tableau and a set of boxes that forms a horizontal strip as input (with the precondition that, when the inner shape of the tableau is expanded to include these boxes, it remains a valid partition) and outputs a tableau (we will prove later that the output is indeed a valid semistandard tableau).

\begin{alg}[Full Multi-Insertion] \label{fullmulti}
\begin{center}\end{center}
\vspace{.5cm}
\begin{algorithmic}[1]
Function FullMulti(tableau $R$, set $S$ of boxes) \Comment{The inner shape of $R$ plus the boxes in $S$ must be a valid partition; no box in $S$ is in the inner shape of $R$. Also, $S$ must form a horizontal strip.}
\State $\mu :=$ inner shape of $R$.
\State $\lambda :=$ outer shape of $R$.
\State $\mathfrak{q}_0 :=$ empty insertion queue.
\State Choose any integer $r$. \label{forchooser} \Comment{We prove later that the choice of $r$ is immaterial.}
\State $h := r$.
\While{$h \neq r + k$}: \Comment{$k$ is the vertical period of $R$.}
\State $L :=$ list of boxes in row $r$ in $S$, from left to right.
\While{$L$ is not empty}:
\State $B :=$ first element of $L$.
\State Remove $B$ from $L$.
\If{$B \in R$}:
\State $x :=$ entry of $R$ in $B$.
\State Put $(x,h + 1)$ into $\mathfrak{q}_0$. \label{putinq0}
\State Expand $\mu$ to include $B$ and remove $B$ from $R$.
\Else: \label{fordegen}
\State Expand $\mu$ and $\lambda$ to include $B$ and remove $B$ from $R$. \label{expmulambda}
\EndIf.
\EndWhile.
\State $h := h + 1$.
\EndWhile.
\State $R_0 := R$. \label{r0}
\State $i := 0$.
\While{$\mathfrak{q}_i$ is not empty}: \label{callosmloop}
\State $(R_{i+1},\mathfrak{q}_{i+1}) := \OneStepMulti{(R_i,\mathfrak{q}_i)}$. \label{callos} \Comment{This line inserts all elements of $\mathfrak{q}_i$ into $R_i$ and denotes the resulting tableau and queue as $R_{i + 1}$ and $\mathfrak{q}_{i + 1}$, respectively.}
\State $i := i + 1$.
\EndWhile.
\State $R := R_i$.
\State \Return $R$.
\end{algorithmic}
\end{alg}

\begin{example}
Let $R$, drawn below, be the input tableau into Algorithm \ref{fullmulti}, and let the red boxes in the diagram below constitute $S$.
\begin{center}
$R=$
\begin{tabular}{m{.2cm}m{.2cm}m{.2cm}m{.2cm}m{.2cm}m{.2cm}m{.2cm}}
&&&&$\vdots$&$\vdots$&$\vdots$\\
\hdashline&&&\multicolumn{1}{c|}{}&2&3&\multicolumn{1}{c|}{5}\\
\cline{4-4}\cline{6-7}&&\multicolumn{1}{c|}{}&\begin{subtab}\cellcolor{red!50}2\end{subtab}&\multicolumn{1}{c|}{6}\\
\cline{2-3}\cline{5-5}\multicolumn{1}{c|}{}&\begin{subtab}\cellcolor{red!50}1\end{subtab}&\begin{subtab}\cellcolor{red!50}2\end{subtab}&\multicolumn{1}{c|}{4}\\
\hdashline\multicolumn{1}{c|}{}&2&3&\multicolumn{1}{c|}{5}\\
&$\vdots$&$\vdots$&$\vdots$
\end{tabular}
\end{center}
Suppose that we pick row $0$ (the top row) to be our starting row $r$. We expand $\mu$ to include the red boxes, deleting those boxes from $R$ and putting corresponding entries into $\mathfrak{q}_0$. The $2$ is deleted from row $1$, and it is to be inserted into the following row, so we add $(2,2)$ to $\mathfrak{q}_0$. We then add $(1,0)$ (the same as $(1,3)$, since the vertical period of $R$ is $3$) and $(2,0)$ into $\mathfrak{q}_0$. We let $R_0$ be our current tableau.
\begin{center}
$\mathfrak{q}_0=(2,2),(1,0),(2,0)$
\hspace{1cm}
$R_0=$
\begin{tabular}{m{.2cm}m{.2cm}m{.2cm}m{.2cm}m{.2cm}m{.2cm}m{.2cm}}
&&&&$\vdots$&$\vdots$&$\vdots$\\
\hdashline&&&\multicolumn{1}{c|}{}&2&3&\multicolumn{1}{c|}{5}\\
\cline{6-7}&&&\multicolumn{1}{c|}{}&\multicolumn{1}{c|}{6}\\
\cline{4-4}\cline{5-5}&&\multicolumn{1}{c|}{}&\multicolumn{1}{c|}{4}\\
\hdashline\cline{2-3}\multicolumn{1}{c|}{}&2&3&\multicolumn{1}{c|}{5}\\
&$\vdots$&$\vdots$&$\vdots$
\end{tabular}
\end{center}
We now enter the One-Step Multi-Insertion subroutine. For this subroutine, set $\mathfrak{q} = \mathfrak{q}_0$, and let $\mathfrak{q}'$ be the empty queue. We remove $(2,2)$ from $\mathfrak{q}$ and insert $2$ into row $2$, bumping out the $4$ and adding $(4,0)$ to $\mathfrak{q}'$. We remove $(1,0)$ from $\mathfrak{q}$ and insert $1$ into row $0$, bumping out the $2$ and adding $(2,1)$ to $\mathfrak{q}'$. We remove $(2,0)$ from $\mathfrak{q}$ and insert $2$ into row $0$, bumping out the $3$ and adding $(3,1)$ to $\mathfrak{q}'$. Finally, with $\mathfrak{q}$ empty, we exit the subroutine. We let $\mathfrak{q}_1$ and $R_1$ be the returned insertion queue and the returned tableau, respectively; they are shown below:
\begin{center}
$\mathfrak{q}_1=(4,0),(2,1),(3,1)$
\hspace{1cm}
$R_1=$
\begin{tabular}{m{.2cm}m{.2cm}m{.2cm}m{.2cm}m{.2cm}m{.2cm}m{.2cm}}
&&&&$\vdots$&$\vdots$&$\vdots$\\
\hdashline&&&\multicolumn{1}{c|}{}&1&2&\multicolumn{1}{c|}{5}\\
\cline{6-7}&&&\multicolumn{1}{c|}{}&\multicolumn{1}{c|}{6}\\
\cline{4-4}\cline{5-5}&&\multicolumn{1}{c|}{}&\multicolumn{1}{c|}{2}\\
\hdashline\cline{2-3}\multicolumn{1}{c|}{}&1&2&\multicolumn{1}{c|}{5}\\
&$\vdots$&$\vdots$&$\vdots$
\end{tabular}
\end{center}
After another run through the subroutine, we have:
\begin{center}
$\mathfrak{q}_2=(5,1),(6,2)$
\hspace{1cm}
$R_2=$
\begin{tabular}{m{.2cm}m{.2cm}m{.2cm}m{.2cm}m{.2cm}m{.2cm}m{.2cm}}
&&&&$\vdots$&$\vdots$&$\vdots$\\
\hdashline&&&\multicolumn{1}{c|}{}&1&2&\multicolumn{1}{c|}{4}\\
\cline{7-7}&&&\multicolumn{1}{c|}{}&2&\multicolumn{1}{c|}{3}\\
\cline{4-4}\cline{5-6}&&\multicolumn{1}{c|}{}&\multicolumn{1}{c|}{2}\\
\hdashline\cline{2-3}\multicolumn{1}{c|}{}&1&2&\multicolumn{1}{c|}{4}\\
&$\vdots$&$\vdots$&$\vdots$
\end{tabular}
\end{center}
Notice that a box got added to $R_2$. This happened because, when $(3,1)$ was removed from $\mathfrak{q}$ in the subroutine, the largest entry of row $1$ was $2$, so $3$ was added at the end of row $1$ and nothing was bumped out. For the same reason, $\mathfrak{q}_2$ now has only two elements.

After a third iteration of the subroutine, we have:
\begin{center}
$\mathfrak{q}_3$ is empty
\hspace{1cm}
$R_3=$
\begin{tabular}{m{.2cm}m{.2cm}m{.2cm}m{.2cm}m{.2cm}m{.2cm}m{.2cm}}
&&&&$\vdots$&$\vdots$&$\vdots$\\
\hdashline&&&\multicolumn{1}{c|}{}&1&2&\multicolumn{1}{c|}{4}\\
&&&\multicolumn{1}{c|}{}&2&3&\multicolumn{1}{c|}{5}\\
\cline{4-4}\cline{6-7}&&\multicolumn{1}{c|}{}&2&\multicolumn{1}{c|}{6}\\
\hdashline\cline{2-3}\cline{5-5}\multicolumn{1}{c|}{}&1&2&\multicolumn{1}{c|}{4}\\
&$\vdots$&$\vdots$&$\vdots$
\end{tabular}
\end{center}
Now that $\mathfrak{q}_3$ is empty, we exit out of the loop, let $R$ equal $R_3$, and return $R$.
\end{example}

If $S$ consists of only one box, then Algorithm \ref{fullmulti} functions as Algorithm \ref{intinsert}. This means that when we prove certain properties of full multi-insertion (such as the fact that it returns a valid semistandard tableau), we will be showing the analogous properties to be true of single insertion as well.

\begin{remark} \label{fmtermin}
Algorithm \ref{fullmulti} always ends after a finite number of steps.
\end{remark}

\begin{proof}
Consider any entry that is originally bumped out of $R$. This entry will bump out a larger entry, and that entry will bump out a larger entry, and so on. However, the bumped-out entry can only be as large as the largest entry that was originally in $R$. Therefore, any chain of bumps will eventually end with a landing; when this happens to every chain, the queue becomes empty and Algorithm \ref{fullmulti} terminates.
\end{proof}

The concept of a bumping route can be extended to multi-insertion. The intuitive idea is that we keep track of an entry's element in the queue and, when we insert the entry, we put the box's corresponding point (in the following row) into the bumping route. Thus, when we perform multi-insertion on a tableau and a set of boxes, we can construct the bumping route of any point corresponding to any box in the set. Our formal extension of the definition of a bumping route to multi-insertion is below.

\begin{defin}
Let $R$ be a tableau and $S$ be a set of boxes such that $\FullMulti(R, S)$ is well-defined. Let $P$ any point such that $\pi(P) \in S$. We will define the \textit{bumping route} of $P$ with input box set $S$ as a list of points that is constructed as follows:\footnote{The input box set $S$ will be implicit in all future discussion of bumping routes.}

Suppose that we perform Algorithm \ref{fullmulti} with $R$ and $S$ as input. We will modify Algorithms \ref{osmulti} and \ref{fullmulti} for the purposes of this definition. At the very beginning, we let $s$ be the plane row of $P$, and we initialize $H$ as a list with the single entry $P$. On line \ref{putinq0} of Algorithm \ref{fullmulti}, if $B = \pi^{-1}(P)$, instead of putting $(x, h + 1)$ into $\mathfrak{q}_0$, put $(x, h + 1, s + 1)$ into the queue.\footnote{For the purposes of this definition only, we extend the definition of an insertion queue to include triples of the form (letter, row, plane row).}

In Algorithm \ref{osmulti}, if in an iteration of the loop beginning on line \ref{onesteploop} a three-element member of $\mathfrak{q}$ is removed (let it be $(x, r, s)$), on line \ref{putinq} add $(x', r + 1, s + 1)$ into $\mathfrak{q'}$ instead of $(x', r + 1)$. Furthermore, if a three-element member $(x, r, s)$ is removed from the queue, at the very end of that iteration of the loop, append $\pi^{-1}_s(B)$ to the end of $H$.

The bumping route of $P$ is defined to be the list $H$ as it stands after Algorithm \ref{fullmulti} terminates.
\end{defin}

\begin{prop} \label{validsst}
For any $R$ and $S$ satisfying the preconditions of Algorithm \ref{fullmulti}, after every bump that occurs during the execution of Algorithm \ref{fullmulti} with inputs $R$ and $S$, $R$ is a valid semistandard tableau. Also, for every bumping route created through the insertions caused by Algorithm \ref{fullmulti}, each point of the bumping route (except the first) is weakly left of the point before it --- that is, all thus generated bumping routes trend weakly left.
\end{prop}

\begin{proof}
Let $A$ be the alphabet of $R$. Define $\QuA = A \cup \{-\infty,\infty\}$. We compare $-\infty$ and $\infty$ to all elements of $\QuA$ as such: $-\infty < a < \infty$ for all $a \in \QuA$ (including $a = -\infty$ and $a = \infty$). That is, we have a quasiorder on $\QuA$: in addition to the ordering of $A$, we have $-\infty < a < \infty$ for all $a \in A$, as well as $-\infty < -\infty$, $-\infty < \infty$, and $\infty < \infty$.

In this proof, we will think of a tableau as a map from $\ScrC_{k,n}$ (not just a part of $\ScrC_{k,n}$) to $\QuA$, where the images of all boxes in the inner shape of the tableau are $-\infty$ and the images of all boxes not in the outer shape of the tableau are $\infty$. A similar technique was used by Donald Knuth in his studies of regular tableaux \cite[\S3]{Knu}. If $R$ is a tableau by our regular definition, then we will denote $\widehat{R}$ to be the corresponding tableau by our new definition. Note that, by looking at $\widehat{R}$, we know $R$ and the two partitions that bound $R$.

For this proof, we will modify line \ref{newbox} of Algorithm \ref{osmulti} as such: ``$B:=$ box of the leftmost entry of $\widehat{R}$ in row $r$ that is \emph{not less-than-or-equal-to} $x$" (note that ``greater than" and ``less than or equal to" are not opposites in our quasiorder). This will allow us to think of an element of $S$ being bumped out as $-\infty$ being inserted into the corresponding row, and it will allow us to think of a box landing in a square as bumping out $\infty$ and the $\infty$ being inserted ``infinitely far to the right" in the following row.

\begin{defin}
Given a tableau $\widehat{R}$ by our new definition, define the \textit{inner shape} of $\widehat{R}$ to be the partition containing exactly the $-\infty$'s in $\widehat{R}$. Define the \textit{outer shape} of $\widehat{R}$ to be the partition that does not contain exactly the $\infty$'s in $\widehat{R}$.
\end{defin}

Our definition above is consistent with the partitions of $R$ during the execution of Algorithm \ref{fullmulti}; by looking at the placement of $-\infty$ and $\infty$ in $\widehat{R}$, we can determine $R$'s inner and outer shape. This consistency includes line \ref{expmulambda} of Algorithm \ref{fullmulti}: $-\infty$ is placed into $B$, and thus $B$ becomes part of both the inner and outer shapes of both $\widehat{R}$ and $R$.

Suppose that we know that, at some point during the execution of Algorithm \ref{fullmulti}, the entries of $\widehat{R}$ increase weakly from left to right and strictly from top to bottom (that is, $\widehat{R}$ is semistandard). Then it follows that $R$ is bounded by two partitions (otherwise, the columns would not be strictly increasing, or the rows would fail to have the form (infinitely many $-\infty$'s, some elements of $A$, infinitely many $\infty$'s)). Thus, for the first part of our proposition, it suffices to show that, at every point during the execution of Algorithm \ref{fullmulti}, the entries of $\widehat{R}$ increase weakly from left to right and strictly from top to bottom.

We will prove our proposition by induction. We know that $R$ (and therefore $\widehat{R}$) is a valid semistandard tableau to start with and that all bumping routes (all with $0$ boxes thus far) trend weakly left. Suppose that $R$ (and $\widehat{R}$) is a valid semistandard tableau after every bump up to but not including the insertion of an entry $x$ into a box $B$ (in the case that $B \in S$, we have $x = -\infty$ for the purposes of insertion into $\widehat{R}$). Suppose further that all bumping routes created up to this point in the execution of Algorithm \ref{fullmulti} have trended weakly left. We will show that, after the insertion of $x$ into $B$, $R$ remains a valid semistandard tableau, and the bumping route that $B$ becomes a part of after this insertion still trends weakly left.

Just before the insertion of $x$ into $B$, $R$ is semistandard. Let $r$ be the row containing box $B$. Since $x$ is put into the leftmost box containing an entry not less-than-or-equal-to $x$, all boxes directly to the left of $B$ (that is, all boxes in the row of $B$ to the left of $B$) have entries that are less than or equal to $x$; furthermore, since $x$ is less than the entry previously in $B$, we have that $x$ is less than or equal to all entries directly to the right of $B$. Similarly, we have that $x$ is less than all entries directly below $B$ (in the column of $B$ below $B$). To show that $R$ continues to be a valid semistandard tableau after the insertion of $x$ into $B$, it thus remains only to show that $x$ is greater than all entries directly above $B$ when it is inserted into $B$.

Suppose that $x$ is inserted into $B$ during the formation of $\widehat{R_j}$, for some $j \ge 0$ (that is, on the $j$'th call of Algorithm \ref{osmulti} or, if $j = 0$, before line \ref{r0} of Algorithm \ref{fullmulti}).

\textbf{Case 1: $j = 0$.} In this case, the bumping route containing $B$ only has $B$, so clearly it trends weakly left. Also, by the horizontal strip restriction on $S$, we have that all entries above $B$ are $-\infty$, as desired.

\textbf{Case 2: $j > 0$.} In this case, let $B'$ be the box that is immediately above $B$ (in row $r - 1$), let $C'$ be the box in row $r - 1$ from which the $x$ that is inserted into $B$ was bumped (during the formation of $\widehat{R_{j - 1}}$), and let $C$ be the box immediately below $C'$ (in row $r$). First we show that the bumping route that includes $B$ after the insertion of $x$ into $B$ trends weakly left after $B$ is added to it --- that is, that $B$ is weakly to the left of $C$. To show this, it suffices to show that the entry in $C$ in $\widehat{R_{j - 1}}$ is not less-than-or-equal-to $x$ (because of how $B$ is chosen). (Note that the entry in $C$ in $\widehat{R_{j - 1}}$ is the entry in $C$ just before $x$ is inserted into $B$. This is because during the formation of $\widehat{R_j}$, pairs of the insertion queue with second element $r$ are processed based on the first elements' boxes in row $r - 1$ from left to right; since $x$ was in $C'$, all entries thus far processed came from boxes in $r - 1$ to the left of $C'$, and thus landed in boxes in $r$ to the left of $C$, by our inductive hypothesis about bumping routes.)

Suppose, for contradiction, that the entry in $C$ in $\widehat{R_{j - 1}}$ --- call it $y$ --- is less than or equal to $x$. Then $\widehat{R_{j - 2}}$ (or $\widehat{R}$, in the case of $j = 1$), which contains $x$ in $C'$, cannot contain $y$ in $C$. Since $y$ is in $C$ in $R_{j - 1}$, but is not in $C$ in $R_{j - 2}$, we have that $y$ is inserted into $C$ during the formation of $\widehat{R_{j - 1}}$. We know that $C \not \in S$, because $C$ is below $C'$ (and $S$ is a horizontal strip). This means that $y$ is bumped out of row $r - 1$. In fact, since $y$ is inserted into row $r$ during the formation of $R_{j - 1}$, it follows that $y$ is bumped out of row $r - 1$ during the formation of $\widehat{R_{j - 2}}$ --- a nonsensical concept and, therefore, a contradiction if $j = 1$. If $j > 1$, then by our inductive hypothesis, we know that $y$ is bumped from a box that is weakly to the right of $C'$ during the formation of $\widehat{R_{j - 2}}$. Hence, a box weakly to the right of $C'$ contains an entry that is not greater-than-or-equal-to $y$ in $\widehat{R_{j - 2}}$, while $C'$ contains $x$ in $\widehat{R_{j - 2}}$ --- a contradiction, since $y \le x$. Therefore, our inductive step holds for the part of our proposition that concerns bumping routes.

We now know that $B$ is weakly to the left of $C$ and that $B'$ is weakly to the left of $C'$. Thus, the entry in $B'$ right before the insertion of $x$ into $B$ is less than or equal to the entry in $C'$ at this time, which is in turn less than $x$ (since $x$ was previously bumped out of $C'$). Therefore, $x$ is greater than the entry in $B'$ (and all the entries directly above $B'$) when it is inserted into $B$, as desired.

Having completed our induction, we have proven our proposition.
\end{proof}

Suppose that $S$ contains only one element (call it $B$). In this case, the precondition of the inner shape of $R$ plus $B$ being a valid partition is equivalent to $B$ being an inside cocorner. Thus, with only one element in $S$, Algorithm \ref{fullmulti} is reduced to Algorithm \ref{intinsert}. Therefore, Algorithm \ref{intinsert} necessarily produces a valid semistandard tableau, and the bumping route of any $P \in \pi^{-1}(B)$ trends weakly left.

\begin{remark} \label{forimmaterial}
The output of Algorithm \ref{fullmulti} does not depend on the choice of $r$ in line \ref{forchooser} of Algorithm \ref{fullmulti}.
\end{remark}

\begin{proof}
Take two values of $r$ --- call them $r_1$ and $r_2$ --- and perform Algorithm \ref{fullmulti} on a tableau $R$ and a set $S$ of boxes using $r_1$ and $r_2$. For any $i$, let $R_i(r_1)$ be $R_i$ produced when the algorithm is run using $r = r_1$ and let $R_i(r_2)$ be $R_i$ produced when the algorithm is run using $r = r_2$. Similarly, let $\mathfrak{q}_i(r_1)$ be $\mathfrak{q}_i$ produced when the algorithm is run using $r = r_1$ and let $\mathfrak{q}_i(r_2)$ be $\mathfrak{q}_i$ produced when the algorithm is run using $r = r_2$.

We will proceed by induction on $i$. We know that $R_0(r_1) = R_0(r_2)$ since the same boxes are taken out of $R$ either way. We also know that $\mathfrak{q}_0(r_1)$ contains the same elements as $\mathfrak{q}_0(r_2)$, and that both insertion queues are regular, as entries were taken out from $R$ within any particular row from left to right (and thus from least to greatest). Suppose that we have $R_g(r_1) = R_g(r_2)$ and the queues $\mathfrak{q}_g(r_1)$ and $\mathfrak{q}_g(r_2)$ are regular and are permutations of each other (contain the same elements) for all $g < i$ for some positive $i$. Then, by Remark \ref{qequiv}, we have $R_i(r_1) = R_i(r_2)$, and the queues $\mathfrak{q}_i(r_1)$ and $\mathfrak{q}_i(r_2)$ are permutations of each other. Furthermore, both $\mathfrak{q}_i(r_1)$ and $\mathfrak{q}_i(r_2)$ are regular, by Remark \ref{regularity}. Having completed our induction, we have shown that $R_i(r_1) = R_i(r_2)$ for all $i$, including the final value that $i$ takes on. Thus, the output of Algorithm \ref{fullmulti} does not depend on the choice of $r$.
\end{proof}

\subsection{The Cylindric Row-Bumping Lemma and its Various Corollaries}
\begin{lemma} \label{bumpincrease}
Let $R$ be a tableau on which multi-insertion is performed. Then for all rows $r$ of $R$, the list of entries bumped out of $r$ in the order in which they were bumped out is weakly increasing.
\end{lemma}

\begin{proof}
It suffices to show that for any two entries bumped consecutively from a given row $s$ (``consecutively'' meaning that nothing is bumped out from $s$ in the meantime), the first of the two entries that is bumped is less than or equal to the second entry that is bumped. We will proceed by induction on the moment of time the second entry is bumped. Suppose that, for every row $s$, any two entries $e_1$ and $e_2$ that are bumped consecutively from $s$ (in that order) satisfy $e_1 \le e_2$, up to but not necessarily including the bumping of entries $a_1$ and $a_2$ (in that order) from a row $r$. We will show that $a_1 \le a_2$.

Suppose that $a_1 \in S$ (meaning that the box containing $a_1$ lies in $S$) and $a_2 \in S$. Since $a_1$ is bumped first, $a_1$ is to the left of $a_2$, and so $a_1 \le a_2$. Suppose that $a_1 \in S$ and $a_2 \not \in S$. Then $a_1$ and $a_2$ are both originally in row $r$ in $R$ ($a_2$ is originally in row $r$ because it is bumped out consecutively after $a_1$ is), with $a_1$ to the left of $a_2$. Again, it follows that $a_1 \le a_2$. Clearly, it cannot be the case that $a_1 \not \in S$ and $a_2 \in S$. (Note that this is why we do not need a base case for our induction: we do not use our inductive hypothesis for the entries originally bumped from $R$.)

Now, suppose that $a_1 \not \in S$ and $a_2 \not \in S$. Let $b_1$ and $b_2$ be the entries that bump out $a_1$ and $a_2$, respectively. We have that $b_1$ is inserted into $r$ before $b_2$ is. It follows that $b_1$ is bumped from row $r - 1$ before $b_2$ is (because queues work on a first-in-first-out basis). It follows by our inductive hypothesis that $b_1 \le b_2$. Since $b_1$ and $b_2$ are inserted consecutively into row $r$, it follows that $b_2$ is placed strictly to the right of $b_1$. Therefore, $a_1 \le a_2$, as desired.
\end{proof}

\begin{defin}
Let $H$ be a bumping route (created by an insertion algorithm, such as Algorithm \ref{intinsert} and Algorithm \ref{fullmulti}). For all plane rows $r$ such that $H$ has a point in $r$, $H(r)$ will denote that point. $L(H)$ will denote the number of points in $H$, and for $1 \le i \le L(H)$, $H_i$ will denote the $i$'th point of $H$.
\end{defin}

Define a total ordering on points as follows: $P = (x_1,y_1) \le Q = (x_2,y_2)$ if and only if (a) $x_1 < x_2$ ($P$ is above $Q$), or (b) $x_1 = x_2$ and $y_1 \ge y_2$ ($P$ and $Q$ are in the same row, and $P$ is weakly right of $Q$). (The $\ge$ sign is intended.)

\begin{theorem}[Cylindric Row-Bumping Lemma] \label{crblemma}
Let $R$ be a tableau and $S$ be a set of boxes such that $R$ and $S$ satisfy the input preconditions for Algorithm \ref{fullmulti}. Let $G$ and $H$ be two bumping routes created when Algorithm \ref{fullmulti} is performed with $R$ and $S$ as input, such that $G_1 < H_1$. Then for all plane rows $r$ such that $G(r)$ and $H(r)$ are defined, $H(r)$ is strictly to the left of $G(r)$ and the bump that extended $G$ to include $G(r)$ came after the bump that extended $H$ to include $H(r)$.
\end{theorem}

\begin{proof}
Let $s$ be the plane row that $H_1$ is in. Since $G_1 < H_1 = H(s)$, $G_1$ is in $s$ or is above $s$. If $G_{L(G)}$ is in a plane row above $s$, then all elements of $G$ are above all elements of $H$, meaning that $G(r)$ and $H(r)$ are not both defined for any $r$; in this case, we are done. Otherwise, we know that $G(s)$ is to the right of $H(s)$: if $G_1 = G(s)$, then $G_1$ must be to the right of $H_1$ in order for $G_1 < H_1$ to hold; if $G(s) = G_m$ for some $m > 1$, then $G(s)$ must be to the right of $H(s)$ because $\pi(H(s))$ is a box that was originally removed from the tableau, and thus $H(s)$ is permanently to the left of the tableau. Similarly, we know that the bump that extended $G$ to include $G(s)$ came after the bump that extended $H$ to include $H(s)$: if $G_1 = G(s)$, then $G(s)$ is to the right of $H(s)$, which means that it is bumped out later (since the original bumps proceed left to right across rows); if $G(s) = G_m$ for some $m > 1$, then $G$ reaches plane row $s$ on a later call of Algorithm \ref{osmulti}.

Suppose that $H(r)$ is to the left of $G(r)$ and that the bump that extended $G$ to include $G(r)$ came after the bump that extended $H$ to include $H(r)$ for all $r$ such that $s \le r < t$. We are to show that $H(t)$ is to the left of $G(t)$ and that the bump that extended $G$ to include $G(t)$ came after the bump that extended $H$ to include $H(t)$.

Let $u = \pi(t)$. Let $g$ and $h$ be the entries bumped out of $\pi(G(t - 1))$ and $\pi(H(t - 1))$ when $G$ and $H$ were extended to include $G(t - 1)$ and $H(t - 1)$, respectively. If $G(s) = G_m$ for some $m > 1$, then $G$ was extended to include $G(t)$ on a later call of Algorithm \ref{osmulti} than the call on which $H$ was extended to include $H(t)$; thus, in this case, the second part of our claim is obviously true. If $G(s) = G_1$, then $(h,u)$ and $(g,u)$ are in the same queue, but $(h,u)$ comes first by our inductive hypothesis; it follows that the insertion of $h$ into row $u$ came before the insertion of $g$ into row $u$. Thus, the bump that extended $G$ to include $G(t)$ did indeed come after the bump that extended $H$ to include $H(t)$.

By our inductive hypothesis, $h$ was bumped out of row $u - 1$ before $g$ was. It follows by Lemma \ref{bumpincrease} that $h \le g$. Let $h'$ be the entry in $\pi(H(t))$ at the time that $g$ is inserted into row $u$. We know that $h' \le h$ (Remark \ref{boxdecrease}). Thus, $h' \le g$. It follows that, at this time, $g$ is greater than or equal to the entry in $\pi(H(t))$ and all entries to the left of $\pi(H(t))$ in row $u$. Thus, $g$ is inserted in row $u$ strictly to the right of $H(t)$. Therefore, $H(t)$ is to the left of $G(t)$, as desired.

Having completed our induction, we have shown that $H(r)$ is to the left of $G(r)$, and that the bump that extended $G$ to include $G(r)$ came after the bump that extended $H$ to include $H(r)$, for all plane rows $r$ such that $G(r)$ and $H(r)$ are defined.
\end{proof}

\begin{cor} \label{oner}
Let $G$ and $H$ be two bumping routes such that, for some plane row $r$, $H(r)$ is strictly to the left of $G(r)$. Then $G_1 < H_1$, and for all plane rows $s$ such that $G(s)$ and $H(s)$ are defined, $H(s)$ is strictly to the left of $G(s)$.
\end{cor}

\begin{proof}
Suppose, for contradiction, that $G_1 \ge H_1$. Clearly, $G_1 \neq H_1$, because otherwise $G$ and $H$ would be the same bumping route and it would not be the case that $H(r)$ is strictly to the left of $G(r)$ for some $r$. Thus, $G_1 > H_1$, meaning that $H_1 < G_1$. By Theorem \ref{crblemma}, this implies that $G(r)$ is to the left of $H(r)$, a contradiction. Thus, $G_1 < H_1$. The rest of the corollary follows directly from Theorem \ref{crblemma}.
\end{proof}

\begin{cor} \label{onerlater}
Let $G$ and $H$ be two bumping routes such that, for some plane row $r$, the bump that extends $G$ to include $G(r)$ came after the bump that extended $H$ to include $H(r)$. Then $G_1 < H_1$, and for all plane rows $s$ such that $G(s)$ and $H(s)$ are defined, the bump that extends $G$ to include $G(s)$ came after the bump that extended $H$ to include $H(s)$.
\end{cor}

\begin{proof}
Suppose, for contradiction, that $G_1 \ge H_1$. Clearly, $G_1 \neq H_1$, because otherwise $G$ and $H$ would be the same bumping route and it would not be the case that the bump that extends $G$ to include $G(r)$ came after the bump that extended $H$ to include $H(r)$ for some $r$. Thus, $G_1 > H_1$, meaning that $H_1 < G_1$. By Theorem \ref{crblemma}, this implies that the bump that extends $H$ to include $H(r)$ came after the bump that extended $G$ to include $G(r)$, a contradiction. Thus, $G_1 < H_1$. The rest of the corollary follows directly from Theorem \ref{crblemma}.
\end{proof}

\begin{cor} \label{pointonce}
No point can be part of two different bumping routes that are created by the same application of multi-insertion. (Bumping routes of two different points are considered different, even if the points correspond to the same box.)
\end{cor}

\begin{proof}
Suppose, for contradiction, that this is false. Then there exists a point $P$ that is part of two different bumping routes $G$ and $H$. Without loss of generality, assume that $G_1 < H_1$ (we know that $G_1 \neq H_1$ because $G$ and $H$ are different bumping routes). Let $r(P)$ be the plane row that $P$ is in. Then $G(r(P))$ and $H(r(P))$ are both defined. It follows by Theorem \ref{crblemma} that $H(r(P))$ is to the left of $G(r(P))$. However, by our assumption, $H(r(P)) = G(r(P)) = P$. This is a contradiction. Thus, no point can be part of two different bumping routes.
\end{proof}

\begin{defin}
Given any bumping route $H$, the \textit{cylindric bumping route} $\pi(H)$ is the list of boxes $\pi(H_1)$, $\pi(H_2)$, $ \dots$, $\pi(H_{L(H)})$. Given any cylindric bumping route $K$, the set of bumping routes $H$ such that $\pi(H) = K$ is denoted $\pi^{-1}(K)$.
\end{defin}

\begin{cor} \label{boxonce}
No box is part of two different cylindric bumping routes. No box appears twice or more among the elements of any cylindric bumping route.
\end{cor}

\begin{proof}
Suppose for contradiction that a box $B$ is part of two different cylindric bumping routes $J$ and $K$. Let $P$ be any element of $\pi^{-1}(B)$. Then there exists a bumping route $G \in \pi^{-1}(J)$ such that $P \in G$. Similarly, there exists a bumping route $H \in \pi^{-1}(K)$ such that $P \in H$. Clearly, $G$ and $H$ are distinct bumping routes, since $\pi(G) \neq \pi(H)$. This is a contradiction, by Corollary \ref{pointonce}.

Suppose for contradiction that a box $B$ appears twice or more among the elements of a cylindric bumping route $K$. Let $x$ be the entry bumped out of $B$ as a result of the bump that extends $K$ to include $B$ for the first time. Let $y$ be the entry bumped out of $B$ as a result of the bump that extends $K$ to include $B$ for the second time (we know that this is a bump, not a landing, because $B$ is already in the tableau); let $z$ be the entry that bumps $y$ out during this bump. We know that $z \ge x$, because bumped-out entries along a bumping route increase; we know that $y \le x$ by Remark \ref{boxdecrease}. It follows that $z \ge y$, a contradiction, as $z$ bumps out $y$. Thus, no box appears twice or more among the elements of any cylindric bumping route, as desired.
\end{proof}

It follows from Corollary \ref{boxonce} that the entry that is in a particular box can only change once during an application of multi-insertion. This means that when an entry is inserted into a box, that entry remains in the box throughout the rest of the multi-insertion process. (In turn, this means that once an entry lands in a box, it cannot later be bumped from the box.) An entry cannot be bumped out more than once during an application of multi-insertion.

\begin{cor} \label{insertincrease}
Let $R$ be a tableau on which multi-insertion is performed. Then for all rows $r$ of $R$, the list of entries inserted into $r$ in the order in which they were inserted is weakly increasing.
\end{cor}

\begin{proof}
Suppose, for contradiction, that this is not the case. Then there exists a row $r$ such that an entry is inserted into $r$, and then a smaller entry is inserted into $r$. Choose any $s \in \pi^{-1}(r)$. Let $P$ and $Q$ be points in $s$ such that an entry $a$ is inserted into $P$, and then $b < a$ is inserted into $Q$. Let $G$ be the bumping route that contains $P$ and $H$ be the bumping route that contains $Q$ (we know that $G$ and $H$ are unique from Corollary \ref{pointonce}). From Corollary \ref{onerlater}, we know that $H_1 < G_1$. From Theorem \ref{crblemma}, we conclude that $P$ is to the left of $Q$.

At the moment when $Q$ becomes part of $H$, $a$ has already been inserted into $P$, and from Corollary \ref{boxonce} we know that $a$ stays in $P$ for the rest of the insertion process. This means that $b$ is inserted into $Q$ while a larger entry $a$ is in $P$, to its left. This is a contradiction. Therefore, for all rows $r$ of $R$, the list of entries inserted into $r$ in the order in which they were inserted is weakly increasing.
\end{proof}

\begin{cor} \label{startend}
Let $G$ and $H$ be two bumping routes. Then $G_1 < H_1$ if and only if $G_{L(G)} < H_{L(H)}$.
\end{cor}

\begin{proof}
We will first prove that if $G_1 < H_1$, then $G_{L(G)} < H_{L(H)}$. Suppose, for contradiction, that this is not the case --- that there exist two bumping routes $G$ and $H$ such that $G_1 < H_1$, but $G_{L(G)} \ge H_{L(H)}$. It follows that the final point of $G$ is in the same plane row as the final point of $H$, or it is in a lower plane row. If it is in the same plane row, it follows by Theorem \ref{crblemma} that $G_{L(G)}$ is to the right of $H_{L(H)}$, which means that $G_{L(G)} < H_{L(H)}$ --- a contradiction.

If the final point of $G$ is in a lower plane row than the final point of $H$, let $r$ be the plane row such that $H_{L(H)} = H(r)$. Then $G(r)$ is defined, as otherwise it would not be the case that $G_1 < H_1$ (as then $G$ would start lower than $H$). It follows by Theorem \ref{crblemma} that $H(r)$ is to the left of $G(r)$. However, since $H(r) = H_{L(H)}$, $\pi(H(r))$ is a box where an entry lands. This means that an entry lands in $\pi(G(r))$ as well. We know from Corollary \ref{boxonce} that an entry could not have landed in $\pi(G(r))$ and then have been bumped out. It follows that the insertion that extends $G$ to include $G(r)$ is a landing. This is a contradiction, as then the final point of $G$ cannot be in a lower plane row than the final point of $H$.

Now we will prove that if $G_{L(G)} < H_{L(H)}$, then $G_1 < H_1$. Suppose, for contradiction, that $G_1 \ge H_1$. We know that $G_1 \neq H_1$, because then $G$ and $H$ are the same bumping route, and $G_{L(G)} = H_{L(H)}$. If $G_1 > H_1$, then, by the first part of this corollary (which we have proved), $H_{L(H)} < G_{L(G)}$. This is a contradiction. Therefore, if $G_{L(G)} < H_{L(H)}$, then $G_1 < H_1$.
\end{proof}

\begin{defin} \label{newsetdef}
Let $\lambda$ be the outer shape of a tableau $R$. Let $S$ be a set of boxes such that $R$ and $S$ satisfy the input preconditions of Algorithm \ref{fullmulti}. Let $R' = \FullMulti(R,S)$ and $\lambda'$ be the outer shape of $R'$. Then the \textit{new set} associated with $R$ and $S$, denoted $N(R,S)$, is the set of boxes in $\lambda' / \lambda$.
\end{defin}

\begin{cor} \label{horstrip}
For any tableau $R$ and set of boxes $S$ such that $R$ and $S$ are valid inputs into Algorithm \ref{fullmulti}, $N(R,S)$ forms a horizontal strip.
\end{cor}

\begin{proof}
Suppose, for contradiction, that $N(R,S)$ does not form a horizontal strip. Then there are two boxes $B, C \in N(R,S)$, with $C$ immediately below $B$. Then there exist two points $P$ and $Q$ such that $P \in \pi^{-1}(B)$, $Q \in \pi^{-1}(C)$, and $P$ is immediately above $Q$. Let $G$ and $H$ be the bumping routes that $P$ and $Q$ are a part of, respectively; thus, $P = G_{L(G)}$ and $Q = H_{L(H)}$. If $L(H) = 1$ (meaning that $C$ was not in the tableau to begin with), $P$ cannot be above $Q$, because then $B$ would have had to be in $\lambda$ to begin with, contradicting our supposition that $B \in N(R,S)$. If $L(H) > 1$, consider $Q' = H_{L(H)-1}$. We know that $\pi(P)$ is not in the outer shape of $R$ to begin with. This means that $\pi(P)$ was added to $R$. By Corollary \ref{boxonce}, once an entry has landed in $\pi(P)$, it is never bumped out. Thus, if $H_{L(H)-1}$ were $P$ or any point directly to the right of $P$, $Q$ would end at that landing point. This means that $H_{L(H)-1}$ is in the row of $P$, but to the left of $P$ (and thus to the left of $Q$). This is a contradiction, as bumping routes trend weakly left. Thus, $N(R,S)$ necessarily forms a horizontal strip.
\end{proof}

\section{Reverse Insertion and Reverse Multi-Insertion}
\subsection{Reverse Insertion Algorithms and Examples}
Analogous to the process of row-deletion for non-cylindric tableaux, Algorithm \ref{intinsert} can be reversed in a process that we will call ``reverse row-insertion." The algorithm for reverse row-insertion is shown below. We will prove later that reverse row-insertion is indeed a process that reverses internal row-insertion, and that the process always results in a valid semistandard tableau. The algorithm takes a tableau $R$ and an outside corner $B$ of $R$, defined as follows:

\begin{defin}
Given $R \in \SSCT(\lambda / \mu)$, a box $B = \pi((i,j))$ is an \textit{outside corner} of $R$ if $B \in \lambda$, but $\pi((i+1,j)) \not \in \lambda$ and $\pi((i,j+1)) \not \in \lambda$. Visually, $B$ is an outside cocorner of $R$ if it is to the left of (inside) $\lambda$, but the boxes to the right of and below $B$ are outside $\lambda$. Note that $B$ is not necessarily a box of $R$, as it is possible that $B$ is inside $\mu$ as well.
\end{defin}

\begin{alg}[Reverse Row-Insertion] \label{revinsert}
\begin{center}\end{center}
\vspace{.5cm}
\begin{algorithmic}[1]
Function ReverseInsert(tableau $R$, box $B$) \Comment{$B$ must be an outside corner of $R$.}
\State $\mu :=$ inner shape of $R$.
\State $\lambda :=$ outer shape of $R$.
\If{$B \in R$}:
\State $x :=$ entry of $R$ that is in $B$.
\EndIf.
\State Shrink $\lambda$ to exclude $B$ and remove $B$ from $R$.
\If{$B \in \mu$}: \Comment{This happens only when $B$ was not in $R$ to begin with.}
\State Shrink $\mu$ to exclude $B$.
\Else:
\While{$x \neq$ null}: \label{revloop}
\State $r:=$ row of $B$.
\If{$x$ is less than or equal to every entry in $R$ in row $r - 1$}:
\State $B:=$ rightmost box of row $r - 1$ that is in $\mu$.
\State Put $x$ in $B$ and shrink $\mu$ to exclude $B$. \Comment{We say that $x$ \textit{lands} in $B$.}
\State $x :=$ null.
\Else:
\State $B:=$ box of the rightmost entry of $R$ in row $r - 1$ that is less than $x$.
\State $x':=$ entry of $R$ in $B$.
\State In $R$, replace $x'$ in $B$ with $x$.
\State $x := x'$.
\EndIf.
\EndWhile.
\EndIf.
\State \Return $R$.
\end{algorithmic}
\end{alg}

For an example of reverse row-insertion, the reader may look at Example \ref{insertex}, starting with the last tableau drawn as $R$ and the box with the $7$ in the top row drawn as $B$. The tableaux obtained after each insertion are the same as the tableaux in Example \ref{insertex}, except in reverse order. The algorithm ends with the original tableau in the example.

Just as a particular box's entry can only decrease during insertion, the a particular box's entry can only increase during reverse insertion, since entries always bump out smaller entries.

\begin{defin}
Bumping routes for reverse insertion are defined exactly as bumping routes are defined for insertion. Consider a point $P$ such that $\pi(P)$ was reverse-inserted into a tableau $R$. The \textit{reverse bumping route} of $P$ is a list of points, constructed as follows:
\begin{itemize}
\item Add $P$ to the reverse bumping route when $\pi(P)$ is reverse-inserted into $R$.
\item Say that a point $Q = (r,y)$ is added to the reverse bumping route when an entry $x$ is bumped out of $\pi(Q)$. When $x$ is reverse-inserted into a box $B \in \pi(r - 1)$, add the element of $\pi^{-1}(B)$ that is in plane row $r - 1$ to the bumping route.
\end{itemize}
\end{defin}

Since $x$ decreases after every iteration of the loop beginning on line \ref{revloop} and $x$ cannot be smaller than the smallest entry in the original tableau, the loop must terminate and the algorithm must end. Thus, reverse insertion, just like insertion, always ends after a finite number of steps.

\begin{defin}
An insertion queue is \textit{reverse-regular} if, for any two elements of the queue $(x_1,r)$ and $(x_2,r)$, where $x_1 < x_2$, $(x_1,r)$ comes \emph{after} $(x_2,r)$ in the queue.
\end{defin}

Just as with insertion, one can reverse-insert multiple entries at the same time. The algorithm below is analogous to Algorithm \ref{osmulti}; it takes a tableau and a reverse-regular insertion queue as input and outputs a pair consisting of a map from a subset of $\ScrC_{k,n}$ to the alphabet of the tableau taken as input (though the map itself is not necessarily a valid semistandard tableau) and an insertion queue.

\begin{alg}[Reverse One-Step Multi-Insertion] \label{revosmulti}
\begin{center}\end{center}
\vspace{.5cm}
\begin{algorithmic}[1]
Function ReverseOneStepMulti(tableau $R$, insertion queue $\mathfrak{q}$) \Comment{$\mathfrak{q}$ must be reverse-regular.}
\State $\mu :=$ inner shape of $R$.
\State $\mathfrak{q}' :=$ empty insertion queue.
\While{$\mathfrak{q}$ is not empty}:
\State Remove the first element from $\mathfrak{q}$. Let it be $(x,r)$.
\If{$x$ is less than or equal to every entry in $R$ in row $r$}:
\State $B:=$ rightmost box of row $r$ that is in $\mu$.
\State Put $x$ in $B$ and shrink $\mu$ to exclude $B$. \Comment{$\mu$ is not necessarily a valid partition anymore.}
\Else:
\State $B:=$ box of the rightmost entry of $R$ in row $r$ that is less than $x$.
\State $x':=$ entry of $R$ in $B$.
\State In $R$, replace $x'$ in $B$ with $x$.
\State Add $(x',r - 1)$ to $\mathfrak{q}'$.
\EndIf.
\EndWhile.
\State \Return $(R,\mathfrak{q}')$. \Comment{$R$ is not necessarily a tableau.}
\end{algorithmic}
\end{alg}

\begin{remark}
The insertion queue returned by Algorithm \ref{revosmulti} is reverse-regular.
\end{remark}

\begin{proof}
Let $(R',\mathfrak{q}') = \ReverseOneStepMulti(R,\mathfrak{q})$ for a tableau $R$ and a reverse-regular insertion queue $\mathfrak{q}$. Suppose, for contradiction, that $\mathfrak{q}'$ is not reverse-regular. Then there exist two elements of $\mathfrak{q}'$, $(y_1,r)$ and $(y_2,r)$, such that $y_1 > y_2$, but $(y_2,r)$ comes before $(y_1,r)$ in $\mathfrak{q}'$. Let $x_1$ and $x_2$ be the entries that bumped out $y_1$ and $y_2$, respectively. Then $(x_2,r+1)$ and $(x_1,r+1)$ were in $\mathfrak{q}$, with $(x_2,r+1)$ coming first. Since $\mathfrak{q}$ is reverse-regular, it follows that $x_2 > x_1$. We also know that $x_1 > y_1$, since $x_1$ bumps out $y_1$. It follows that $x_2 > x_1 > y_1 > y_2$. Since $y_1 < y_2$, $y_1$ is to the right of $y_2$ in $R$, and $y_1 < x_2$. It follows that $y_2$ is not the rightmost entry of $R$ in its row that is greater than $x_2$ at the time that $x_2$ is to be inserted, which is a contradiction, because then $x_2$ does not bump out $y_2$. Thus, $\mathfrak{q}'$ is a reverse-regular insertion queue.
\end{proof}

We can fully reverse-insert multiple entries from a tableau simultaneously using Algorithm \ref{revosmulti} as a subroutine, in a way analogous to Algorithm \ref{fullmulti}. This algorithm takes a tableau and a set of boxes that forms a horizontal strip as input (with the precondition that, when the outer shape of the tableau is shrunk to exclude these boxes, it remains a valid partition) and outputs a tableau (we will prove later that the output is indeed a valid semistandard tableau).

\begin{alg}[Reverse Full Multi-Insertion] \label{revfullmulti}
\begin{center}\end{center}
\vspace{.5cm}
\begin{algorithmic}[1]
Function ReverseFullMulti(tableau $R$, set $S$ of boxes) \Comment{The outer shape of $R$ minus the boxes in $S$ must be a valid partition; every box in $S$ is in the outer shape of $R$. Also, $S$ must form a horizontal strip.}
\State $\mu :=$ inner shape of $R$.
\State $\lambda :=$ outer shape of $R$.
\State $\mathfrak{q}_0 :=$ empty insertion queue.
\State Choose any integer $r$. \label{revchooser} \Comment{We prove later that the choice of $r$ is immaterial.}
\State $h := r$.
\While{$h \neq r - k$}: \Comment{$k$ is the vertical period of $R$.}
\State $L :=$ list of boxes in row $r$ in $S$, from right to left.
\While{$L$ is not empty}:
\State $B :=$ first element of $L$.
\State Remove $B$ from $L$.
\If{$B \in R$}:
\State $x :=$ entry of $R$ in $B$.
\State Put $(x,h - 1)$ into $\mathfrak{q}_0$.
\State Shrink $\lambda$ to exclude $B$ and remove $B$ from $R$.
\Else:
\State Shrink $\mu$ and $\lambda$ to exclude $B$ and remove $B$ from $R$.
\EndIf.
\EndWhile.
\State $h := h - 1$.
\EndWhile.
\State $R_0 := R$.
\State $i := 0$.
\While{$\mathfrak{q}_i$ is not empty}:
\State $(R_{i+1},\mathfrak{q}_{i+1}) := \ReverseOneStepMulti{(R_i,\mathfrak{q}_i)}$. \Comment{This line inserts all elements of $\mathfrak{q}_i$ into $R_i$ and calls the result $R_{i+1}$. $\mathfrak{q}_{i+1}$ is the insertion queue of bumped-out cells.}
\State $i := i + 1$.
\EndWhile.
\State $R := R_i$.
\State \Return $R$.
\end{algorithmic}
\end{alg}

\begin{example}
Let $R$, drawn below, be the input tableau into Algorithm \ref{revfullmulti}, and let the red boxes in the diagram below constitute $S$.
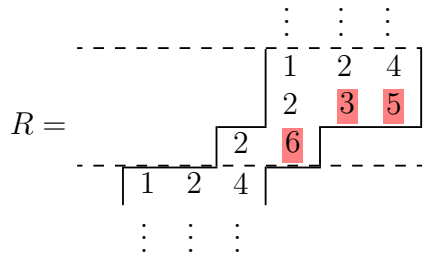
\begin{center}
$R =$
\begin{tabular}{m{.2cm}m{.2cm}m{.2cm}m{.2cm}m{.2cm}m{.2cm}m{.2cm}}
&&&&$\vdots$&$\vdots$&$\vdots$\\
\hdashline&&&\multicolumn{1}{c|}{}&1&2&\multicolumn{1}{c|}{4}\\
&&&\multicolumn{1}{c|}{}&2&\begin{subtab}\cellcolor{red!50}3\end{subtab}&\multicolumn{1}{c|}{\begin{subtab}\cellcolor{red!50}5\end{subtab}}\\
\cline{4-4}\cline{6-7}&&\multicolumn{1}{c|}{}&2&\multicolumn{1}{c|}{\begin{subtab}\cellcolor{red!50}6\end{subtab}}\\
\hdashline\cline{2-3}\cline{5-5}\multicolumn{1}{c|}{}&1&2&\multicolumn{1}{c|}{4}\\
&$\vdots$&$\vdots$&$\vdots$
\end{tabular}
\end{center}
Suppose that we pick row $1$ to be our starting row $r$. We shrink $\lambda$ to exclude the red boxes, deleting those boxes from $R$ and putting corresponding entries into $\mathfrak{q}_0$. The $5$ is deleted from row $1$, and it is to be reverse-inserted into the previous row, so we add $(5,0)$ to $\mathfrak{q}_0$. We then add $(3,0)$ to $\mathfrak{q}_0$ ($(5,0)$ is added first because entries are processed from right to left, not from left to right as done in Algorithm \ref{fullmulti}). Finally, $(6,1)$ is added to $\mathfrak{q}_0$. We let $R_0$ be our current tableau.
\begin{center}
$\mathfrak{q}_0 = (5,0), (3,0), (6,1)$
\hspace{1cm}
$R_0 =$
\begin{tabular}{m{.2cm}m{.2cm}m{.2cm}m{.2cm}m{.2cm}m{.2cm}m{.2cm}}
&&&&$\vdots$&$\vdots$&$\vdots$\\
\hdashline&&&\multicolumn{1}{c|}{}&1&2&\multicolumn{1}{c|}{4}\\
\cline{6-7}&&&\multicolumn{1}{c|}{}&\multicolumn{1}{c|}{2}\\
\cline{4-4}\cline{5-5}&&\multicolumn{1}{c|}{}&\multicolumn{1}{c|}{2}\\
\hdashline\cline{2-3}\multicolumn{1}{c|}{}&1&2&\multicolumn{1}{c|}{4}\\
&$\vdots$&$\vdots$&$\vdots$
\end{tabular}
\end{center}
We now enter the One-Step Reverse Multi-Insertion subroutine. We remove $(5,0)$ from $\mathfrak{q}$ and insert $5$ into row $0$, bumping out the $4$ (the rightmost entry strictly less than $5$) and adding $(4,2)$ to $\mathfrak{q}'$. We remove $(3,0)$ from $\mathfrak{q}$ and insert $3$ into row $0$, bumping out the $2$ and adding $(2,2)$ to $\mathfrak{q}'$. We remove $(6,1)$ from $\mathfrak{q}$ and insert $6$ into row $1$, bumping out the $2$ and adding $(2,0)$ to $\mathfrak{q}'$. Finally, with $\mathfrak{q}$ empty, we exit the subroutine. We let $\mathfrak{q}_1$ and $R_1$ be the returned insertion queue and the returned tableau, respectively; they are shown below:
\begin{center}
$\mathfrak{q}_1 = (4,2), (2,2), (2,0)$
\hspace{1cm}
$R_1 =$
\begin{tabular}{m{.2cm}m{.2cm}m{.2cm}m{.2cm}m{.2cm}m{.2cm}m{.2cm}}
&&&&$\vdots$&$\vdots$&$\vdots$\\
\hdashline&&&\multicolumn{1}{c|}{}&1&3&\multicolumn{1}{c|}{5}\\
\cline{6-7}&&&\multicolumn{1}{c|}{}&\multicolumn{1}{c|}{6}\\
\cline{4-4}\cline{5-5}&&\multicolumn{1}{c|}{}&\multicolumn{1}{c|}{2}\\
\hdashline\cline{2-3}\multicolumn{1}{c|}{}&1&3&\multicolumn{1}{c|}{5}\\
&$\vdots$&$\vdots$&$\vdots$
\end{tabular}
\end{center}
The process continues as follows:
\begin{center}
$\mathfrak{q}_2 = (2,1), (1,2)$
\hspace{1cm}
$R_2 =$
\begin{tabular}{m{.2cm}m{.2cm}m{.2cm}m{.2cm}m{.2cm}m{.2cm}m{.2cm}}
&&&&$\vdots$&$\vdots$&$\vdots$\\
\hdashline&&&\multicolumn{1}{c|}{}&2&3&\multicolumn{1}{c|}{5}\\
\cline{6-7}&&&\multicolumn{1}{c|}{}&\multicolumn{1}{c|}{6}\\
\cline{3-4}\cline{5-5}&\multicolumn{1}{c|}{}&2&\multicolumn{1}{c|}{4}\\
\hdashline\cline{2-2}\multicolumn{1}{c|}{}&2&3&\multicolumn{1}{c|}{5}\\
&$\vdots$&$\vdots$&$\vdots$
\end{tabular}
\end{center}
\begin{center}
$\mathfrak{q}_3$ is empty
\hspace{1cm}
$R_3 =$
\begin{tabular}{m{.2cm}m{.2cm}m{.2cm}m{.2cm}m{.2cm}m{.2cm}m{.2cm}}
&&&&$\vdots$&$\vdots$&$\vdots$\\
\hdashline&&&\multicolumn{1}{c|}{}&2&3&\multicolumn{1}{c|}{5}\\
\cline{4-4}\cline{6-7}&&\multicolumn{1}{c|}{}&2&\multicolumn{1}{c|}{6}\\
\cline{2-3}\cline{5-5}\multicolumn{1}{c|}{}&1&2&\multicolumn{1}{c|}{4}\\
\hdashline\multicolumn{1}{c|}{}&2&3&\multicolumn{1}{c|}{5}\\
&$\vdots$&$\vdots$&$\vdots$
\end{tabular}
\end{center}
\end{example}
Now that $\mathfrak{q}_3$ is empty, we exit out of the loop, let $R$ equal $R_3$, and return $R$.

Since entries only bump out smaller entries and no bumped-out entry can be smaller than the smallest entry originally in the tableau, any chain of bumps eventually ends. It follows that Algorithm \ref{revfullmulti} always terminates.

In the same way that the concept of a bumping route applies to multi-insertion as well as single-box insertion, the concept of a reverse bumping route applies to reverse multi-insertion as well as reverse single-box insertion.

\subsection{Relating Reverse Insertion to Forward Insertion}
\begin{theorem} \label{retrace}
Let $R$ be a tableau and $S$ be a set of boxes such that $R$ and $S$ satisfy the input preconditions of Algorithm \ref{fullmulti}. Let $R' = \FullMulti(R,S)$ and $S' = N(R,S)$ (see Definition \ref{newsetdef}). Let $B$ be a box in $S$. Choose $P \in \pi^{-1}(B)$ and let $H$ be the bumping route of $P$. Let $P' = H_{L(H)}$ and let $H'$ be the reverse bumping route of $P'$ when Algorithm \ref{revfullmulti} is called on $R'$ and $S'$. Then, for all $1 \le i \le L(H)$, $H'_i = H_{L(H) + 1 - i}$. That is, $H'$ consists of the same points as $H$, but in reverse order.
\end{theorem}

\begin{proof}
We first note that the preconditions for Algorithm \ref{revfullmulti} are satisfied: $S'$ forms a horizontal strip (Corollary \ref{horstrip}) that can be deleted from the outer shape of $R'$.

Take any $R$ and $S$ as described in the theorem statement. We know by definition that, for all $P \in \bigcup\limits_{B \in S}\pi^{-1}(B)$ (with $P'$, $H$, and $H'$ defined correspondingly as in the theorem statement), $H'_1 = P' = H_{L(H)}$. Suppose that, up to but not including the reverse-insertion of an entry $x$, all reverse bumping routes created by the reverse multi-insertion process have retraced the bumping routes created by the multi-insertion process. Let $A$ be the box containing $x$ prior to this reverse insertion. Select any $Q \in \pi^{-1}(A)$. Let $H$ be the bumping route containing $Q$ that was made by the multi-insertion algorithm and $H'$ be the reverse bumping route containing $Q$ that was made by the reverse multi-insertion algorithm. Let $j$ be such that $H'_j = Q$. (Thus, $H'_j = H_{L(H) + 1 - j}$.) We will show that $H'_{j+1} = H_{L(H) - j}$.

Let $r$ be the plane row immediately above the plane row containing $Q$. Let $C$ be the box into which $x$ is to be reverse-inserted if $H'_{j+1} = H_{L(H) - j}$. Suppose, for contradiction, that $x$ is not reverse-inserted into $C$.

We know that $x$ was bumped out of $C$ in the forward multi-insertion process. Let $y$ be the entry that replaced $x$ in $C$. By Corollary \ref{boxonce}, we know that $x$ and $y$ were the only two entries ever to be present in $C$ during the forward multi-insertion process. Since all bumping routes thus far have been retraced, it follows that no entry has thus far been reverse-inserted into $C$. Thus, $y$ is in $C$ when $x$ is to be reverse-inserted into $\pi(r)$. Since $y < x$ and $x$ is not reverse-inserted into $C$, it follows that $y$ is not the rightmost entry in $\pi(r)$ (or $r$) that is less than $x$. Let $z$ be the rightmost entry in $r$ that is less than $x$, and let $z$ be in box $D$.

Let $V$ be the box in plane row $r$ onto which $C$ projects. Let $W$ be the box in plane row $r$ onto which $D$ projects. We know that $W$ contains an entry less than $x$; however, this was not the case when $y$ was inserted into $C$ during the forward multi-insertion process (because then $z$, which is smaller than $x$, would have been to the right of $x$). Since the entry in $D$ changes, $W$ must be part of a bumping route. Call this bumping route $G$, and let $G'$ be the corresponding reverse bumping route to $G$ created during the reverse multi-insertion process. Since $z$ is still in $W$ by the time that $x$ is to be inserted into plane row $r$, the first point of $G'$ must either be below the first point of $H'$ (in which case $G'$ reaches plane row $r$ on a later call of the subroutine), or be to the left of the first point of $H'$ (in which case it comes later in the queue and continues to be later in the queue, because of the way that the elements of the queue are processed). This means that $H_{L(H)} = H'_1 < G'_1 = G_{L(G)}$.

Since $V$ is to the left of $W$, $H_r$ is to the left of $G_r$. By Corollary \ref{oner}, it follows that $G_1 < H_1$. From this, by Corollary \ref{startend}, it follows that $G_{L(G)} < H_{L(H)}$. This is a contradiction. Therefore, $x$ is reverse-inserted into $C$.

Having completed our induction, we have proven that all reverse bumping routes retrace their corresponding bumping routes during reverse multi-insertion.
\end{proof}

\begin{remark} \label{inv1}
Since all bumping routes are retraced, multi-insertion followed by reverse multi-insertion results in the original tableau. That is, for any $R$ and $S$ satisfying the preconditions of multi-insertion,
\begin{equation*}
\ReverseFullMulti(\FullMulti(R,S),N(R,S)) = R.
\end{equation*}
\end{remark}

By letting $S$ consist of one box, we see that Algorithm \ref{revinsert} retraces the bumping route of Algorithm \ref{intinsert} and returns the original tableau. This means that everything henceforth that we prove about reverse multi-insertion also applies to reverse insertion of a single box.

\begin{defin}
The \textit{negation function} $\Neg$ takes a tableau $R$ and, for every entry $x$ in $R$, changes $x$ to $\overline{x}$. For any two entries $x$ and $y$ such that $x < y$, we have $\overline{y} < \overline{x}$. If the entries of our tableau are integers, we can think of the $\Neg$ operator as taking the additive inverse of every entry in the tableau and returning the result.
\end{defin}

\begin{defin}
Define the \textit{flip function} on a cylindric partition $\lambda$ as follows: $\Flip(\lambda)$ is the partition such that a point $(x,y)$ is in $\Flip(\lambda)$ if and only if it is not in $\lambda$. Equivalently, $\Flip(\lambda)$ is the cylindric partition with sequence $(\dots, -1 - \lambda_1, -1 - \lambda_0, -1 - \lambda_{-1}, \dots)$.

Consider the map $M: \ZZ^2 \rightarrow \ZZ^2$ that takes $(x,y)$ to $(-x,-y)$. Since $M$ takes $(-k,n-k)\ZZ$ to itself, it follows that $M$ induces an endomorphism $E$ of $\ScrC_{k,n}$. If we think of a cylindric tableau as a partial map from the cylinder to a totally ordered set, then, for a cylindric tableau $R$, define $\Flip(R)$ to be the map created by applying $E$ on the cylinder, followed by $R$, followed by $\Neg$. If $R$ has inner and outer shapes $\mu$ and $\lambda$, respectively, then $\Flip(R)$ has inner and outer shapes $\Flip(\lambda)$ and $\Flip(\mu)$, respectively.

We also define the $\Flip$ function for boxes and sets of boxes. If $B$ is a box, then $\Flip(S)$ is the endomorphism $E$ described above applied to $B$. If $S$ is a set of boxes, $\Flip(S)$ is $E$ applied to $S$.

Simply put, let $R$ be a tableau and $B$ be a box of $R$ that contains $t$. Choose any $P \in \pi^{-1}(B)$ and let $P$ have coordinates $(a,b)$. Let $C = \pi((-a,-b))$. Then $C = \Flip(B)$ and, in the tableau $\Flip(R)$, $C$ contains $\overline{t}$.

Visually, $\Flip$ takes a tableau, rotates it by $180^\circ$, and negates all of the tableau's entries. Similarly, $\Flip$ takes a set of boxes and rotates it by $180^\circ$.
\end{defin}

\begin{remark} \label{flipvalidsst}
Let $R$ be a tableau. Then $\Flip(R)$ is also a tableau.
\end{remark}

\begin{proof}
Clearly, $\Flip(R)$ is bounded by two valid partitions --- a partition remains valid after a $180^\circ$ rotation; equivalently, $\Flip$ negates all elements of a partition's sequence, but also reverses their order, keeping the sequence weakly decreasing.

Let $B$ and $C$ be two boxes in $\Flip(R)$ such that $C$ is directly to the right of $B$. Then $\Flip(B)$ is directly to the right of $\Flip(C)$, so, in $R$, $\Flip(B)$ contains an entry that is weakly greater than the entry in $\Flip(C)$. This means that, in $\Flip(R)$, $B$ contains an entry that is weakly less than the entry in $C$, as desired.

Similarly, let $D$ and $E$ be two boxes in $\Flip(R)$ such that $D$ is directly below $E$. Then $\Flip(D)$ is directly above $\Flip(E)$, so, in $R$, $\Flip(E)$ contains an entry that is strictly greater than the entry in $\Flip(D)$. This means that, in $\Flip(R)$, $E$ contains an entry that is strictly less than the entry in $D$, as desired. Therefore, $\Flip(R)$ is a valid semistandard tableau.
\end{proof}

\begin{remark} \label{involflip}
$\Flip$ is an involution --- that is, $\Flip \circ \Flip = \Id$, where $\Id$ is the identity function.
\end{remark}

\begin{remark} \label{fliphorstrip}
Let $S$ be a set of boxes that forms a horizontal strip. Then $\Flip(S)$ forms a horizontal strip.
\end{remark}

\begin{theorem} \label{flipfmflip}
Let $R$ be a tableau and $S$ be a set of boxes such that $R$ and $S$ satisfy the input preconditions of Algorithm \ref{revfullmulti}. Call Algorithm \ref{revfullmulti} on $R$ and $S$; simultaneously, call Algorithm \ref{fullmulti} on $R' = \Flip(R)$ and $S' = \Flip(S)$. Then $R'_0 = \Flip(R_0)$. Furthermore, for every value of $t$, the tableau in Algorithm \ref{fullmulti} after the $t$'th insertion is the flip of the tableau in Algorithm \ref{fullmulti} after the $t$'th insertion; the same number of insertions occurs in the two algorithms.
\end{theorem}

\begin{proof}
We first note that the preconditions for Algorithm \ref{fullmulti} are satisfied: $\Flip(R)$ is a valid semistandard tableau (Remark \ref{flipvalidsst}) and $\Flip(S)$ forms a horizontal strip (Remark \ref{fliphorstrip}) that can be added to the inner shape of $\Flip(R)$ and keep the partition valid. We will show our theorem to be true by simultaneously performing steps of the reverse multi-insertion algorithm on $R$ and $S$ and the forward multi-insertion algorithm on $R'$ and $S'$.

Clearly, $R'_0 = \Flip(R_0)$, as $R_0$ is $R$ without the boxes in $S$, and $R'_0$ is $\Flip(R)$ without the boxes in $\Flip(S)$.

Let $r_1$ be the value of $r$ chosen in line \ref{revchooser} of Algorithm \ref{revfullmulti}. Then choose $-r_1$ for the value of $r$ in line \ref{forchooser} of Algorithm \ref{fullmulti}.

\begin{defin}
The \textit{corresponding element} of an element $(x,r)$ of an insertion queue is $(\overline{x},-r)$. Two insertion queues $\mathfrak{q}_1$ and $\mathfrak{q}_2$ are \textit{corresponding queues} is they contain the same number of elements and, for all $i$, the $i$'th elements of $\mathfrak{q}_1$ and $\mathfrak{q}_2$ are corresponding elements.
\end{defin}

Let $\mathfrak{q}_0$ be the original queue of elements taken out of $R$ during the reverse multi-insertion process. Let $\mathfrak{q}'_0$ be the original queue of elements taken out of $R'$ during the forward multi-insertion process. We first note that $\mathfrak{q}_0$ and $\mathfrak{q}'_0$ are corresponding queues. Since $\mathfrak{q}_0$ contains pairs containing all entries in the boxes of $S$ in $R$ and $\mathfrak{q}'_0$ contains pairs containing all entries in the boxes of $S'$ in $R'$, and these entries are in rows whose numbers are negatives of one another, we know that the elements of $\mathfrak{q}_0$ and $\mathfrak{q}'_0$ can be paired such that each element is paired with its corresponding element. In order to show that these entries are in the same order, consider two boxes $B, C \in S$ and their corresponding boxes $B' = \Flip(B)$ and $C' = \Flip(C)$. If $B$ and $C$ are in the same row, suppose, without loss of generality, that $B$ is to the left of $C$. Then $C$ is deleted from $R$ (and its entry is put into $\mathfrak{q}_0$) before $B$ is deleted from $R$. $B'$ is to the right of $C'$, so the entry in $C'$ is put in a pair that is put into $\mathfrak{q}'_0$ before the entry in $B'$ is put into a pair that is put into $\mathfrak{q}'_0$. Thus, $C$ comes before $B$ and $C'$ comes before $B'$ as desired. If $B$ and $C$ are in different rows, suppose, without loss of generality, that $B$ is deleted from $R$ before $C$ is. Since rows are processed from bottom to top beginning with row $r_1$ in the reverse multi-insertion process and they are processed from top to bottom beginning with row $-r_1$ in the forward multi-insertion process (and the rows of $B'$ and $C'$ are opposite to those of $B$ and $C$, respectively), we conclude that $B'$ is deleted from $R'$ before $C'$ is. Thus, for any two entries in $\mathfrak{q}_0$, their corresponding entries in $\mathfrak{q}'_0$ are in the same order relative to each other as are the two entries in $\mathfrak{q}_0$. We conclude that $\mathfrak{q}_0$ and $\mathfrak{q}'_0$ are corresponding queues.

We have shown that, when the respective subroutines are called on $R_0$ and $R'_0$, the tableaux are flips of each other and the queues referred to as $\mathfrak{q}$ by the subroutines are corresponding. Since the queues referred to as $\mathfrak{q}'$ by the subroutines are initially empty, they are also corresponding. Suppose that we move through both algorithms step by step, examining the tableaux and queues after every insertion. Suppose that, up to but not including the insertion of $x$ into row $s$ in the reverse insertion process, the queues referred to as $\mathfrak{q}$ in the subroutines are corresponding, the queues referred to as $\mathfrak{q}'$ in the subroutines are corresponding, and the two tableaux are flips of each other. Let $V$ be the tableau in the reverse insertion process right before the insertion of $x$ into row $s$ and $V'$ be the corresponding tableau in the forward insertion process (thus, $V' = \Flip(V)$). We will show that the two pairs of queues continue to correspond and that $V' = \Flip(V)$ after the insertion as well.

Suppose that $x$ takes the place of $y$ in row $s$ in $V$. Then $y$ is the rightmost entry smaller than $x$ in $s$. By our inductive hypothesis, in $V'$, $\overline{y}$ is in row $-s$, and, by definition, $\overline{y} > \overline{x}$. Suppose that there is an entry left of $\overline{y}$ in $V'$ that is greater than $\overline{x}$; call this entry $\overline{z}$. Then in $V$, $z < x$ and $z$ is to the right of $y$; this is a contradiction. Thus, $\overline{x}$ takes the place of $\overline{y}$ in row $-s$ in $V'$. The first entries in the queues called $\mathfrak{q}$ by the subroutines are deleted, so these queues continue to correspond. Corresponding entries ($(y,s - 1)$ and $(\overline{y},-s + 1)$) are added to the queues referred to as $\mathfrak{q}'$, so these queues also continue to correspond. Thus, our inductive step holds, as desired.

We have completed our induction. Since both algorithms end when their respective queues (which are corresponding and thus have the same length) are empty, the two algorithms terminate after the same total number of insertions, and, after each successive insertion, the two tableaux remain flips of one another.
\end{proof}

\begin{defin}
Let $\mu$ be the outer shape of a tableau $R$. Let $S$ be a set of boxes such that $R$ and $S$ satisfy the preconditions of Algorithm \ref{revfullmulti}. Let $R' = \ReverseFullMulti(R,S)$ and $\mu'$ be the outer shape of $R'$. Then the \textit{reverse new set} associated with $R$ and $S$, denoted $N^{-1}(R,S)$, is the set of boxes in $\mu / \mu'$.
\end{defin}

We will now create alternative inputs and outputs for the functions $\FullMulti$ and $\ReverseFullMulti$. Each algorithm can now take a pair $(R,S)$ (instead of two arguments $R$ and $S$) and returns a pair $(R',S')$, the first element being the tableau generated by the algorithms and the second being:
\begin{itemize}
\item $S' = N(R,S)$ for $\FullMulti$; and
\item $S' = N^{-1}(R,S)$ for $\ReverseFullMulti$.
\end{itemize}
We will also define $\Flip((R,S))$ to be $(\Flip(R),\Flip(S))$.

\begin{cor} \label{ffmf}
\begin{equation*}
\ReverseFullMulti((R,S)) = \Flip(\FullMulti(\Flip((R,S)))).
\end{equation*}
\end{cor}

\begin{proof}
Given $R$ and $S$, since, after all insertions have been accounted for, the tableaux returned by Algorithm \ref{revfullmulti} and Algorithm \ref{fullmulti} are flips of one another, the tableaux returned by two algorithms must be flips of one another. Thus, flipping the tableau returned by Algorithm \ref{fullmulti} results in the tableau returned by Algorithm \ref{revfullmulti}. Since $\FullMulti(\Flip((R,S)))$ is the flip of $\ReverseFullMulti((R,S))$, it follows that $N^{-1}(R,S)$ is the flip of $N(\Flip((R,S)))$, as desired.
\end{proof}

We can use Theorem \ref{flipfmflip} to prove about reverse multi-insertion much of what we proved about multi-insertion.

\begin{remark}
Given $R$ and $S$ satisfying the preconditions of Algorithm \ref{revfullmulti}, $R_0$ is a valid semistandard tableau. Furthermore, after every insertion that occurs while Algorithm \ref{revfullmulti} is being performed (including during the subroutine), $R$ is a valid semistandard tableau.
\end{remark}

\begin{proof}
We know from Theorem \ref{flipfmflip} that any value of $R$ that appears during the reverse multi-insertion process, starting from $R_0$, is the flip of a value of $R$ that appears during the forward multi-insertion process when it is applied to $\Flip(R)$ and $\Flip(S)$, starting from $\Flip(R)_0$. We know that all values of $R$ starting from $R_0$ during forward multi-insertion are indeed valid semistandard tableaux from Proposition \ref{validsst}. Since the $\Flip$ operator preserves valid semistandard tableaux (Remark \ref{flipvalidsst}), $R$ is a valid semistandard tableau after every insertion, as desired.
\end{proof}

\begin{remark}
The choice of $r$ in line \ref{revchooser} of Algorithm \ref{revfullmulti} is irrelevant.
\end{remark}

\begin{proof}
We have shown that $\ReverseFullMulti(R,S) = \Flip(\FullMulti(\Flip(R),\Flip(S)))$, where the $r$ chosen in the forward insertion process is the opposite of the $r$ chosen in the reverse insertion process. However, the $r$ chosen in the forward insertion process is immaterial (Remark \ref{forimmaterial}). Thus, the choice of $r$ in the reverse insertion process cannot matter either.
\end{proof}

\begin{remark} \label{revhorstrip}
If $(R',S') = \ReverseFullMulti((R,S))$, then $S'$ forms a horizontal strip.
\end{remark}

\begin{proof}
By Corollary \ref{ffmf}, we know that
\begin{equation*}
\ReverseFullMulti((R,S)) = \Flip(\FullMulti(\Flip((R,S)))).
\end{equation*}
Since the $\Flip$ function preserves horizontal strips (Remark \ref{fliphorstrip}) and $\FullMulti$ returns a set of boxes that forms a horizontal strip (Corollary \ref{horstrip}), we conclude that $\Flip(\FullMulti(\Flip((R,S))))$ returns a pair whose second element does indeed form a horizontal strip.
\end{proof}

\begin{remark}
Let $R$ be a tableau and $S$ a set of boxes on which reverse multi-insertion is performed. Then for all rows $r$ of $R$, the list of entries bumped out of $r$ in the order in which they were bumped out is weakly decreasing.
\end{remark}

\begin{proof}
Let $R' = \Flip(R)$ and $S' = \Flip(S)$. By Theorem \ref{flipfmflip}, the reverse insertion tableau is the flip of the forward insertion tableau after every time that an insertion is completed in both algorithms. Also, during the formation of $R_0$ and $R'_0$, entries are bumped out in the same order. This means that if $L = a_1, a_2, \dots, a_m$ is the list of entries bumped out of a row $r$ during the reverse multi-insertion process, the list of entries bumped out of $-r$ during the forward multi-insertion process at the same point during the process is $\overline{a_1}, \overline{a_2}, \dots, \overline{a_m}$.

By Lemma \ref{bumpincrease}, $\overline{a_1}, \overline{a_2}, \dots, \overline{a_m}$ is a weakly increasing list. Therefore, $L = a_1, a_2, \dots, a_m$ is a weakly decreasing list, as desired.
\end{proof}

Analogously, we can use Corollary \ref{insertincrease} to prove that, given a row $r$ of a tableau $R$ on which reverse multi-insertion is performed, the list of entries inserted into $r$ in the order in which they were inserted is weakly decreasing.

\begin{prop} \label{inverse}
$\ReverseFullMulti$ and $\FullMulti$ are inverse operations.
\end{prop}

\begin{proof}
In order to make notation less cumbersome, we will shorten the names of the operations as follows:
\begin{itemize}
\item $\ReverseFullMulti$ to $\RFM$;
\item $\FullMulti$ to $\FM$; and
\item $\Flip$ to $\F$.
\end{itemize}
We have shown that $\RFM(\FM((R,S))) = (R,S)$ (Remark \ref{inv1}) --- that is, $\RFM \circ \FM = \Id$. (The fact that the second element of the output pair is $S$ follows from the fact that the boxes added back to $R$ after being removed are the ones that are originally removed.) It remains to show only that $\FM \circ \RFM = \Id$.

Consider the following facts:
\begin{itemize}
\item From Remark \ref{inv1}:
\begin{equation} \label{eq1} \RFM \circ \FM = \Id\end{equation}
\item From Remark \ref{involflip}:
\begin{equation} \label{eq2} \F \circ \F = \Id\end{equation}
\end{itemize}
We proceed as follows:

Adding $\F$ to the left of both sides of (\ref{eq1}) and applying (\ref{eq2}), we have
\begin{equation*}
\F \circ \RFM = \F \circ \F \circ \FM \circ \F = \Id \circ \FM \circ \F = \FM \circ \F.
\end{equation*}
Adding $\FM$ to the right of both sides, we have
\begin{equation} \label{eq3}
\FM \circ \F \circ \FM = \F \circ \RFM \circ \FM = \F \circ \Id = \F.
\end{equation}
Now adding $\F$ to the right of both sides of (\ref{eq1}) and applying (\ref{eq2}), we have
\begin{equation*}
\RFM \circ \F = \F \circ \FM \circ \F \circ \F = \F \circ \FM \circ \Id = \F \circ \FM.
\end{equation*}
Adding $\FM$ to the left of both sides (here we need Remark \ref{revhorstrip} for well-definedness), we have
\begin{equation*}
\FM \circ \RFM \circ \F = \FM \circ \F \circ \FM.
\end{equation*}
Applying (\ref{eq3}), we have
\begin{equation*}
\FM \circ \RFM \circ \F = \FM \circ \F \circ \FM = \F.
\end{equation*}
Adding $\F$ to the right of both sides, we have
\begin{equation*}
\FM \circ \RFM \circ \F \circ \F = \F \circ \F.
\end{equation*}
Canceling out the $\F$'s, we have
\begin{equation*}
\FM \circ \RFM = \Id,
\end{equation*}
as desired.
\end{proof}

Having proven this, we have created a bijection between pairs consisting of a cylindric tableau and a set of boxes that forms a horizontal strip that can be added to the inner shape of the tableau while keeping the partition valid on the one hand, and pairs consisting of a cylindric tableau and a set of boxes that forms a horizontal strip that can be removed from the outer shape of the tableau while keeping the partition valid on the other. This is the basis of the cylindric Robinson-Schensted-Knuth (RSK) correspondence, which will be detailed in the following section.

\subsection{More Results about Reverse Insertion}
We can use the bijection proven above in order to prove more results about reverse multi-insertion, particularly about bumping routes. This is because, while we proved earlier that forward bumping routes are retraced with reverse multi-insertion (Theorem \ref{retrace}), before this Proposition \ref{inverse} we did not prove that $\FullMulti(\ReverseFullMulti((R,S))) = (R,S)$. This result allows us to use more extensively what we already know about forward multi-insertion in order to prove results about reverse multi-insertion.

\begin{remark}
Reverse bumping routes trend weakly right --- that is, given a reverse bumping route $G$ generated by performing reverse multi-insertion on $(R,S)$, $G(m)$ is to the right of $G_{m - 1}$ for all $1 < m \le L(G)$.
\end{remark}

\begin{proof}
Let $(R',S') = \ReverseFullMulti((R,S))$. We know that $\FullMulti((R',S')) = (R,S)$. Let $G'$ be the bumping route that consists of the same points as $G$ (we know such a $G'$ exists from Theorem \ref{retrace}). $G'$ trends weakly left, and $G$ has the same points as $G'$, but in reverse order. Therefore, $G$ trends weakly right, as desired.
\end{proof}

\begin{remark} \label{revcrblemma}
(Counterpart of Theorem \ref{crblemma}.) Let $R$ be a tableau and $S$ be a set of boxes such that $R$ and $S$ satisfy the input preconditions for Algorithm \ref{revfullmulti}. Let $G$ and $H$ be two reverse bumping routes created when Algorithm \ref{revfullmulti} is performed with $R$ and $S$ as input, such that $G_1 < H_1$. Then for all plane rows $r$ such that $G(r)$ and $H(r)$ are defined, $H(r)$ is strictly to the left of $G(r)$ and the bump that extended $G$ to include $G(r)$ came \emph{before} the bump that extended $H$ to include $H(r)$.
\end{remark}

\begin{proof}
Let $G'$ and $H'$ be the forward bumping routes that correspond to $G$ and $H$, respectively (obtained by performing Algorithm \ref{fullmulti} on the result of $\ReverseFullMulti((R,S))$). We know that $G'_{L(G')} = G_1 < H_1 = H'_{L(H')}$, so, by Corollary \ref{startend}, $G'_1 < H'_1$. By Theorem \ref{crblemma}, this means that, for all plane rows $r$ such that $G'(r)$ and $H'(r)$ are defined, $H'(r)$ is strictly to the left of $G'(r)$. Since $H'(r) = H(r)$ and $G'(r) = G(r)$ for all $r$, it follows that $H(r)$ is strictly to the left of $G(r)$.

Let $s$ be any number such that $G(s)$ and $H(s)$ are defined. If $G_1$ and $H_1$ are in different plane rows, then $G_1$ is above $H_1$, which means that the bump that extended $G$ to include $G(s)$ came in an earlier call of the subroutine than the bump that extended $H$ to include $H(s)$. If $G_1$ and $H_1$ are in the same plane row, let $t$ be the plane row of $G_1$ and $H_1$ and let $u = \pi(t)$. Let $a$ and $b$ be the entries in $\pi(G_1)$ and $\pi(H_1)$, respectively. Since $G_1 < H_1$, $\pi(G_1)$ is to the right of $\pi(H_1)$, which means that $(b,u - 1)$ is inserted into $\mathfrak{q}_0$ before $(a,u - 1)$ is. Because of how entries are taken out of and put into new queues in the subroutine, the pair corresponding to $G$ comes before the pair corresponding to $H$ in $\mathfrak{q}_m$ for all $m$. Thus, in this case as well, the bump that extended $G$ to include $G(s)$ came before the bump that extended $H$ to include $H(s)$. Therefore, for all plane rows $r$ such that $G(r)$ and $H(r)$ are defined, the bump that extended $G$ to include $G(r)$ came before the bump that extended $H$ to include $H(r)$.
\end{proof}

\begin{remark}
(Counterpart of Corollary \ref{oner}.) Let $G$ and $H$ be two reverse bumping routes such that, for some plane row $r$, $H(r)$ is strictly to the left of $G(r)$. Then $G_1 < H_1$, and for all plane rows $s$ such that $G(s)$ and $H(s)$ are defined, $H(s)$ is strictly to the left of $G(s)$.
\end{remark}

\begin{proof}
Let $G'$ and $H'$ be the forward bumping routes corresponding to $G$ and $H$. Since $H(r)$ is strictly to the left of $G(r)$, $H'(r)$ is strictly to the left of $G(r)$. By Corollary \ref{oner}, this implies that, for all plane rows $s$ such that $G'(s)$ and $H'(s)$ are defined, $H'(s)$ is strictly to the left of $G'(s)$, which implies that, for all plane rows $s$ such that $G(s)$ and $H(s)$ are defined, $H(s)$ is strictly to the left of $G(s)$. It also implies that $G'_1 < H'_1$, which, by Corollary \ref{startend}, implies that $G'_{L(G')} < H'_{L(H')}$, which implies that $G_1 < H_1$, as desired.
\end{proof}

\begin{remark}
(Counterpart of Corollary \ref{onerlater}; the proof is entirely analogous, so it is not given here.) Let $G$ and $H$ be two reverse bumping routes such that, for some plane row $r$, the bump that extends $G$ to include $G(r)$ came after the bump that extended $H$ to include $H(r)$. Then $G_1 < H_1$, and for all plane rows $s$ such that $G(s)$ and $H(s)$ are defined, the bump that extends $G$ to include $G(s)$ came before the bump that extended $H$ to include $H(s)$.
\end{remark}

\begin{remark} \label{revpointonce}
(Counterpart of Corollary \ref{pointonce}.) No point can be part of two different reverse bumping routes that are created by the same application of multi-insertion.
\end{remark}

\begin{proof}
Suppose that a point $P$ is part of two reverse bumping routes $G$ and $H$. Let $G'$ and $H'$ be the corresponding forward bumping routes to $G$ and $H$. Then $P$ is on both $G$ and $H$. By Corollary \ref{pointonce}, this is a contradiction.
\end{proof}

\begin{defin}
Given any reverse bumping route $H$, the \textit{cylindric reverse bumping route} $\pi(H)$ is the list of boxes $\pi(H_1)$, $\pi(H_2)$, $ \dots$, $\pi(H_{L(H)})$. Given any cylindric reverse bumping route $K$, the set of bumping routes $H$ such that $\pi(H) = k$ is denoted $\pi^{-1}(K)$.
\end{defin}

\begin{remark}
(Counterpart of Corollary \ref{boxonce}.) No box is part of two different cylindric reverse bumping routes. No box appears twice of more among the elements of a cylindric bumping route.
\end{remark}

\begin{proof}
Suppose there is a box $B$ that is part of two different cylindric reverse bumping routes $J$ and $K$. Let $P$ be any element of $\pi^{-1}(B)$. Then there exists a reverse bumping route $G \in \pi^{-1}(J)$ and a reverse bumping route $H \in \pi^{-1}(K)$ such that $P \in G$ and $P \in H$. By Remark \ref{revpointonce}, this is a contradiction.

Suppose that a box $B$ appears twice of more among the elements of a cylindric reverse bumping route $K$. Choose any $H \in \pi^{-1}(K)$. Let $H'$ be the corresponding forward bumping route to $H$. Then $\pi(H')$ contains $B$ twice or more. By Corollary \ref{boxonce}, this is a contradiction.
\end{proof}

\begin{remark}
(Counterpart of Corollary \ref{startend}.) Let $G$ and $H$ be two reverse bumping routes. Then $G_1 < H_1$ if and only if $G_{L(G)} < H_{L(H)}$.
\end{remark}

\begin{proof}
Let $G'$ and $H'$ be the corresponding forward bumping routes to $G$ and $H$, respectively. We know that $G_1 = G'_{L(G')}$, $H_1 = H'_{L(H')}$, $G_{L(G)} = G'_1$, and $H_{L(H)} = H'_1$. From Corollary \ref{startend}, we know that $G'_1 < H'_1$ if and only if $G'_{L(G')} < H'_{L(H')}$. Therefore, $G_1 < H_1$ if and only if $G_{L(G)} < H_{L(H)}$.
\end{proof}

\section{The Cylindric RSK Correspondence}
\subsection{The Correspondence}
We now present an analog of the Robinson-Schensted-Knuth correspondence for cylindric tableaux. The idea for this correspondence is based on that of Sagan and Stanley \cite[\S6, Theorem 6.1]{SagStan}.
\begin{defin}
Let $S$ be a set and $T_s$ be a set for any $s \in S$. We define the \textit{disjoint union} operator $\bigsqcup$ as follows:
\begin{equation*}
\bigsqcup \limits_{s \in S} T_s = \{(t,s) | s \in S, t \in T_s\}
\end{equation*}
\end{defin}
\begin{theorem} \label{crsk}
Fix two partitions $\alpha$ and $\beta$. Then there exists a bijective mapping
\begin{equation*}
\CRSK: \bigsqcup \limits_{\substack{\mu \in \Cylpar; \\ \mu \subseteq \alpha; \ \mu \subseteq \beta}} \SSCT(\alpha / \mu) \times \SSCT(\beta / \mu) \to \bigsqcup \limits_{\substack{\lambda \in \Cylpar; \\ \alpha \subseteq \lambda; \ \beta \subseteq \lambda}} \SSCT(\lambda / \beta) \times \SSCT(\lambda / \alpha)
\end{equation*}
such that for any $\mu \in \Cylpar$ satisfying $\mu \subseteq \alpha$ and $\mu \subseteq \beta$, $T \in \SSCT(\alpha / \mu)$, and $U \in \SSCT(\beta / \mu)$, if $\CRSK(((T,U),\mu))=((P,Q),\lambda)$, then $\wt(T)=\wt(P)$ and $\wt(U)=\wt(Q)$.
\end{theorem}

\begin{proof}
For $\mu$, $T$, and $U$ as in Theorem \ref{crsk}, we will construct $\CRSK(((T,U),\mu))$ as follows:

\vspace{.5cm}
\begin{alg}[Cylindric RSK] \label{crskalg}
\begin{center}\end{center}
\begin{algorithmic}[1]
\vspace{.5cm}
Function CylindricRSK(((tableau $T$, tableau $U$), partition $\mu$)): \Comment{$\mu$ must be the inner shape of both $T$ and $U$}.
\State $P := T$.
\State $\alpha :=$ outer shape of $T$.
\State $Q :=$ empty tableau with shape $\alpha / \alpha$.
\State $i :=$ smallest entry that is present in any of the boxes in $U$.
\While{true}:
\State $S :=$ set of all boxes that, in $U$, contain $i$.
\State $(P,S') := \FullMulti((P,S))$. \label{callsub} \Comment{That is, $P$ gets assigned to the tableau returned by Algorithm \ref{fullmulti} and $S'$ gets assigned to the returned set of boxes.}
\State Put $i$ into $Q$ in all boxes in $S'$; expand the outer shape of $Q$ to include all boxes in $S'$.
\If{$i$ is less than the largest entry that is present in any of the boxes in $U$}:
\State $i :=$ smallest entry greater than $i$ that is present in any of the boxes in $U$.
\Else:
\State Exit the loop.
\EndIf.
\EndWhile.
\State $\lambda :=$ outer shape of $P$.
\State \Return $((P,Q),\lambda)$.
\end{algorithmic}
\end{alg}

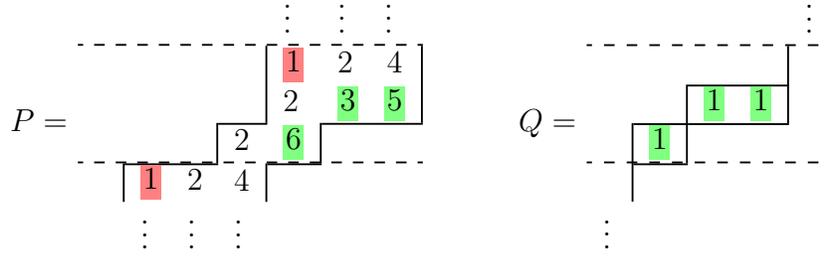
\vspace{.5cm}
\begin{example} \label{bigex}
Suppose that we start with the following $T$ and $U$:
\begin{center}
$T=$
\begin{tabular}{m{.2cm}m{.2cm}m{.2cm}m{.2cm}m{.2cm}m{.2cm}m{.2cm}}
&&&&$\vdots$&$\vdots$&$\vdots$\\
\hdashline&&&\multicolumn{1}{c|}{}&2&3&\multicolumn{1}{c|}{5}\\
\cline{4-4}\cline{6-7}&&\multicolumn{1}{c|}{}&2&\multicolumn{1}{c|}{6}\\
\cline{2-3}\cline{5-5}\multicolumn{1}{c|}{}&1&2&\multicolumn{1}{c|}{4}\\
\hdashline\multicolumn{1}{c|}{}&2&3&\multicolumn{1}{c|}{5}\\
&$\vdots$&$\vdots$&$\vdots$
\end{tabular}
\hspace{1cm}
$U=$
\begin{tabular}{m{.2cm}m{.2cm}m{.2cm}m{.2cm}m{.2cm}m{.2cm}}
&&&&$\vdots$&$\vdots$\\
\hdashline&&&\multicolumn{1}{c|}{}&2&\multicolumn{1}{c|}{4}\\
\cline{4-4}&&\multicolumn{1}{c|}{}&1&3&\multicolumn{1}{c|}{5}\\
\cline{2-3}\cline{6-6}\multicolumn{1}{c|}{}&1&1&3&\multicolumn{1}{c|}{4}\\
\hdashline\cline{4-5}\multicolumn{1}{c|}{}&2&\multicolumn{1}{c|}{4}\\
&$\vdots$&$\vdots$
\end{tabular}
\end{center}
Setting $P$ equal to $T$ and $Q$ equal to the (empty) tableau with shape $\alpha / \alpha$, we have:
\begin{center}
$P=$
\begin{tabular}{m{.2cm}m{.2cm}m{.2cm}m{.2cm}m{.2cm}m{.2cm}m{.2cm}}
&&&&$\vdots$&$\vdots$&$\vdots$\\
\hdashline&&&\multicolumn{1}{c|}{}&2&3&\multicolumn{1}{c|}{5}\\
\cline{4-4}\cline{6-7}&&\multicolumn{1}{c|}{}&\begin{subtab}\cellcolor{red!50}2\end{subtab}&\multicolumn{1}{c|}{6}\\
\cline{2-3}\cline{5-5}\multicolumn{1}{c|}{}&\begin{subtab}\cellcolor{red!50}1\end{subtab}&\begin{subtab}\cellcolor{red!50}2\end{subtab}&\multicolumn{1}{c|}{4}\\
\hdashline\multicolumn{1}{c|}{}&2&3&\multicolumn{1}{c|}{5}\\
&$\vdots$&$\vdots$&$\vdots$
\end{tabular}
\hspace{1cm}
$Q=$
\begin{tabular}{m{.2cm}m{.2cm}m{.2cm}m{.2cm}m{.2cm}}
&&&&$\vdots$\\
\hdashline&&&\multicolumn{1}{c|}{}\\
\cline{3-4}\cline{3-4}&\multicolumn{1}{c|}{}\\
\cline{2-2}\cline{2-2}\multicolumn{1}{c|}{}\\
\hdashline\multicolumn{1}{c|}{}\\
$\vdots$
\end{tabular}
\end{center}
First, $i=1$. We let $S$ be the set of boxes that, in $U$, contain $1$ (shown in red in the diagram of $P$ above). We then perform Algorithm \ref{fullmulti} on $(P,S)$, letting $P$ be equal to the return tableau and letting $S'$ be the returned set of boxes. We put $1$ in $Q$ in all boxes in $S'$. The result is shown below. Boxes that are new to $P$ are shown in green; boxes that are outside (to the right of) the inner shape of $P$, but will be inside (to the left of) the inner shape of $P$ after the next call of Algorithm \ref{fullmulti}, are shown in red.
\begin{center}
$P=$
\begin{tabular}{m{.2cm}m{.2cm}m{.2cm}m{.2cm}m{.2cm}m{.2cm}m{.2cm}}
&&&&$\vdots$&$\vdots$&$\vdots$\\
\hdashline&&&\multicolumn{1}{c|}{}&\begin{subtab}\cellcolor{red!50}1\end{subtab}&2&\multicolumn{1}{c|}{4}\\
&&&\multicolumn{1}{c|}{}&2&\begin{subtab}\cellcolor{green!50}3\end{subtab}&\multicolumn{1}{c|}{\begin{subtab}\cellcolor{green!50}5\end{subtab}}\\
\cline{4-4}\cline{6-7}&&\multicolumn{1}{c|}{}&2&\multicolumn{1}{c|}{\begin{subtab}\cellcolor{green!50}6\end{subtab}}\\
\hdashline\cline{2-3}\cline{5-5}\multicolumn{1}{c|}{}&\begin{subtab}\cellcolor{red!50}1\end{subtab}&2&\multicolumn{1}{c|}{4}\\
&$\vdots$&$\vdots$&$\vdots$
\end{tabular}
\hspace{1cm}
$Q=$
\begin{tabular}{m{.2cm}m{.2cm}m{.2cm}m{.2cm}m{.2cm}}
&&&&$\vdots$\\
\hdashline&&&\multicolumn{1}{c|}{}\\
\cline{3-4}&\multicolumn{1}{c|}{}&\begin{subtab}\cellcolor{green!50}1\end{subtab}&\multicolumn{1}{c|}{\begin{subtab}\cellcolor{green!50}1\end{subtab}}\\
\cline{2-2}\cline{3-4}\multicolumn{1}{c|}{}&\multicolumn{1}{c|}{\begin{subtab}\cellcolor{green!50}1\end{subtab}}\\
\hdashline\cline{2-2}\multicolumn{1}{c|}{}\\
$\vdots$
\end{tabular}
\end{center}
This process continues as follows.
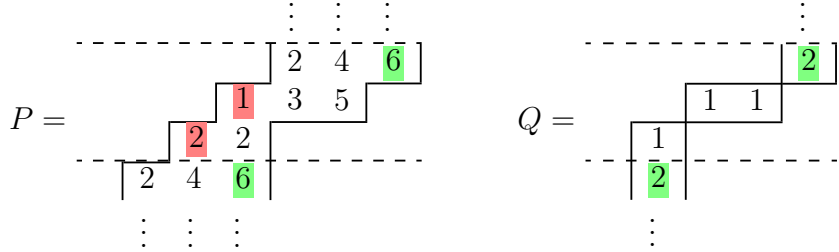
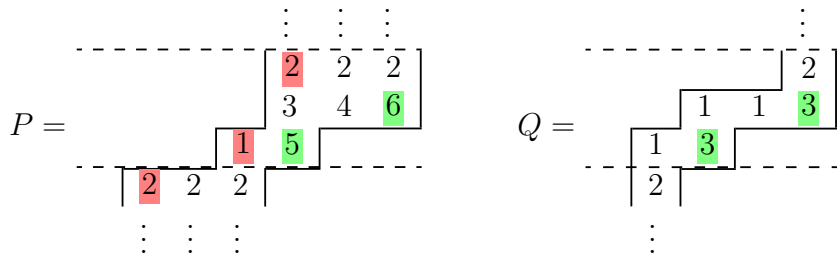
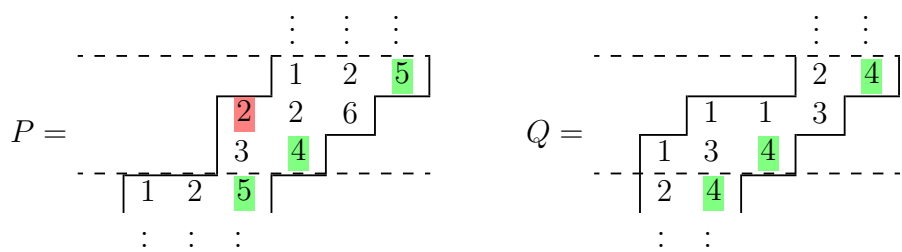
\begin{center}
$P=$
\begin{tabular}{m{.2cm}m{.2cm}m{.2cm}m{.2cm}m{.2cm}m{.2cm}m{.2cm}}
&&&&$\vdots$&$\vdots$&$\vdots$\\
\hdashline&&&\multicolumn{1}{c|}{}&2&4&\multicolumn{1}{c|}{\begin{subtab}\cellcolor{green!50}6\end{subtab}}\\
\cline{4-4}\cline{7-7}&&\multicolumn{1}{c|}{}&\begin{subtab}\cellcolor{red!50}1\end{subtab}&3&\multicolumn{1}{c|}{5}\\
\cline{3-3}\cline{5-6}&\multicolumn{1}{c|}{}&\begin{subtab}\cellcolor{red!50}2\end{subtab}&\multicolumn{1}{c|}{2}\\
\hdashline\cline{2-2}\multicolumn{1}{c|}{}&2&4&\multicolumn{1}{c|}{\begin{subtab}\cellcolor{green!50}6\end{subtab}}\\
&$\vdots$&$\vdots$&$\vdots$
\end{tabular}
\hspace{1cm}
$Q=$
\begin{tabular}{m{.2cm}m{.2cm}m{.2cm}m{.2cm}m{.2cm}}
&&&&$\vdots$\\
\hdashline&&&\multicolumn{1}{c|}{}&\multicolumn{1}{c|}{\begin{subtab}\cellcolor{green!50}2\end{subtab}}\\
\cline{3-4}\cline{5-5}&\multicolumn{1}{c|}{}&1&\multicolumn{1}{c|}{1}\\
\cline{2-2}\cline{3-4}\multicolumn{1}{c|}{}&\multicolumn{1}{c|}{1}\\
\hdashline\multicolumn{1}{c|}{}&\multicolumn{1}{c|}{\begin{subtab}\cellcolor{green!50}2\end{subtab}}\\
&$\vdots$
\end{tabular}
\end{center}
\begin{center}
$P=$
\begin{tabular}{m{.2cm}m{.2cm}m{.2cm}m{.2cm}m{.2cm}m{.2cm}m{.2cm}}
&&&&$\vdots$&$\vdots$&$\vdots$\\
\hdashline&&&\multicolumn{1}{c|}{}&\begin{subtab}\cellcolor{red!50}2\end{subtab}&2&\multicolumn{1}{c|}{2}\\
&&&\multicolumn{1}{c|}{}&3&4&\multicolumn{1}{c|}{\begin{subtab}\cellcolor{green!50}6\end{subtab}}\\
\cline{4-4}\cline{6-7}&&\multicolumn{1}{c|}{}&\begin{subtab}\cellcolor{red!50}1\end{subtab}&\multicolumn{1}{c|}{\begin{subtab}\cellcolor{green!50}5\end{subtab}}\\
\hdashline\cline{2-3}\cline{5-5}\multicolumn{1}{c|}{}&\begin{subtab}\cellcolor{red!50}2\end{subtab}&2&\multicolumn{1}{c|}{2}\\
&$\vdots$&$\vdots$&$\vdots$
\end{tabular}
\hspace{1cm}
$Q=$
\begin{tabular}{m{.2cm}m{.2cm}m{.2cm}m{.2cm}m{.2cm}}
&&&&$\vdots$\\
\hdashline&&&\multicolumn{1}{c|}{}&\multicolumn{1}{c|}{2}\\
\cline{3-4}&\multicolumn{1}{c|}{}&1&1&\multicolumn{1}{c|}{\begin{subtab}\cellcolor{green!50}3\end{subtab}}\\
\cline{2-2}\cline{4-5}\multicolumn{1}{c|}{}&1&\multicolumn{1}{c|}{\begin{subtab}\cellcolor{green!50}3\end{subtab}}\\
\hdashline\cline{3-3}\multicolumn{1}{c|}{}&\multicolumn{1}{c|}{2}\\
&$\vdots$
\end{tabular}
\end{center}
\begin{center}
$P=$
\begin{tabular}{m{.2cm}m{.2cm}m{.2cm}m{.2cm}m{.2cm}m{.2cm}m{.2cm}}
&&&&$\vdots$&$\vdots$&$\vdots$\\
\hdashline&&&\multicolumn{1}{c|}{}&1&2&\multicolumn{1}{c|}{\begin{subtab}\cellcolor{green!50}5\end{subtab}}\\
\cline{4-4}\cline{7-7}&&\multicolumn{1}{c|}{}&\begin{subtab}\cellcolor{red!50}2\end{subtab}&2&\multicolumn{1}{c|}{6}\\
\cline{6-6}&&\multicolumn{1}{c|}{}&3&\multicolumn{1}{c|}{\begin{subtab}\cellcolor{green!50}4\end{subtab}}\\
\hdashline\cline{2-3}\cline{5-5}\multicolumn{1}{c|}{}&1&2&\multicolumn{1}{c|}{\begin{subtab}\cellcolor{green!50}5\end{subtab}}\\
&$\vdots$&$\vdots$&$\vdots$
\end{tabular}
\hspace{1cm}
$Q=$
\begin{tabular}{m{.2cm}m{.2cm}m{.2cm}m{.2cm}m{.2cm}m{.2cm}}
&&&&$\vdots$&$\vdots$\\
\hdashline&&&\multicolumn{1}{c|}{}&2&\multicolumn{1}{c|}{\begin{subtab}\cellcolor{green!50}4\end{subtab}}\\
\cline{3-4}\cline{6-6}&\multicolumn{1}{c|}{}&1&1&\multicolumn{1}{c|}{3}\\
\cline{2-2}\cline{5-5}\multicolumn{1}{c|}{}&1&3&\multicolumn{1}{c|}{\begin{subtab}\cellcolor{green!50}4\end{subtab}}\\
\hdashline\cline{4-4}\multicolumn{1}{c|}{}&2&\multicolumn{1}{c|}{\begin{subtab}\cellcolor{green!50}4\end{subtab}}\\
&$\vdots$&$\vdots$
\end{tabular}
\end{center}
Finally, we have $i = 5$. Since $5$ is the largest entry in $U$, these are our final values of $P$ and $Q$, and they are returned (in addition to their common outer shape) by Algorithm \ref{crskalg}.
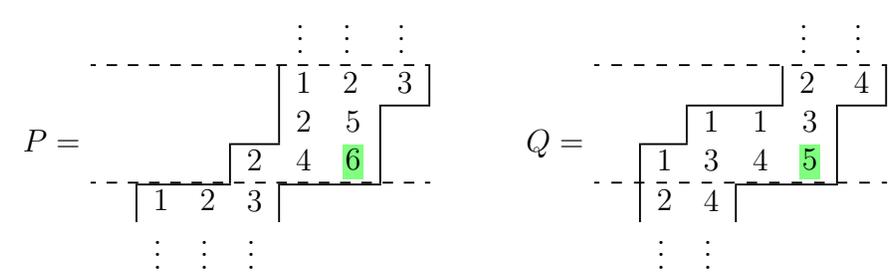
\begin{center}
$P=$
\begin{tabular}{m{.2cm}m{.2cm}m{.2cm}m{.2cm}m{.2cm}m{.2cm}m{.2cm}}
&&&&$\vdots$&$\vdots$&$\vdots$\\
\hdashline&&&\multicolumn{1}{c|}{}&1&2&\multicolumn{1}{c|}{3}\\
\cline{7-7}&&&\multicolumn{1}{c|}{}&2&\multicolumn{1}{c|}{5}\\
\cline{4-4}&&\multicolumn{1}{c|}{}&2&4&\multicolumn{1}{c|}{\begin{subtab}\cellcolor{green!50}6\end{subtab}}\\
\hdashline\cline{2-3}\cline{5-6}\multicolumn{1}{c|}{}&1&2&\multicolumn{1}{c|}{3}\\
&$\vdots$&$\vdots$&$\vdots$
\end{tabular}
\hspace{1cm}
$Q=$
\begin{tabular}{m{.2cm}m{.2cm}m{.2cm}m{.2cm}m{.2cm}m{.2cm}}
&&&&$\vdots$&$\vdots$\\
\hdashline&&&\multicolumn{1}{c|}{}&2&\multicolumn{1}{c|}{4}\\
\cline{3-4}\cline{6-6}&\multicolumn{1}{c|}{}&1&1&\multicolumn{1}{c|}{3}\\
\cline{2-2}\multicolumn{1}{c|}{}&1&3&4&\multicolumn{1}{c|}{\begin{subtab}\cellcolor{green!50}5\end{subtab}}\\
\hdashline\cline{4-5}\multicolumn{1}{c|}{}&2&\multicolumn{1}{c|}{4}\\
&$\vdots$&$\vdots$
\end{tabular}
\end{center}
\end{example}

\begin{example}
The following example is shown to demonstrate what happens when, at some time during the execution of Algorithm \ref{crskalg} (specifically its subroutine, Algorithm \ref{fullmulti}), the inner shape of $P$ is expanded to include boxes that are not present in $P$ (this is the ``Else" in Algorithm \ref{fullmulti} on line \ref{fordegen}).
\begin{center}
$T=$
\begin{tabular}{m{.2cm}m{.2cm}m{.2cm}m{.2cm}}
&&&$\vdots$\\
\hdashline&&\multicolumn{1}{c|}{}&\multicolumn{1}{c|}{2}\\
\cline{2-3}\cline{3-4}\multicolumn{1}{c|}{}&\multicolumn{1}{c|}{1}\\
\hdashline\multicolumn{1}{c|}{}&\multicolumn{1}{c|}{2}\\
&$\vdots$
\end{tabular}
\hspace{1cm}
$U=$
\begin{tabular}{m{.2cm}m{.2cm}m{.2cm}m{.2cm}}
&&&$\vdots$\\
\hdashline&&\multicolumn{1}{c|}{}&\multicolumn{1}{c|}{2}\\
\cline{2-3}\multicolumn{1}{c|}{}&1&1&\multicolumn{1}{c|}{3}\\
\hdashline\cline{3-4}\multicolumn{1}{c|}{}&\multicolumn{1}{c|}{2}\\
&$\vdots$
\end{tabular}
\end{center}
We have
\begin{center}
$P=$
\begin{tabular}{m{.2cm}m{.2cm}m{.2cm}m{.2cm}}
&&&$\vdots$\\
\hdashline&&\multicolumn{1}{c|}{}&\multicolumn{1}{c|}{2}\\
\cline{2-3}\cline{3-4}\multicolumn{1}{c|}{}&\multicolumn{1}{c|}{\begin{subtab}\cellcolor{red!50}1\end{subtab}}&\begin{subtab}\cellcolor{red!50}\phantom{1}\end{subtab}\\
\hdashline\multicolumn{1}{c|}{}&\multicolumn{1}{c|}{2}\\
&$\vdots$
\end{tabular}
\hspace{1cm}
$Q=$
\begin{tabular}{m{.2cm}m{.2cm}m{.2cm}m{.2cm}}
&&&$\vdots$\\
\hdashline&&\multicolumn{1}{c|}{}\\
\cline{2-3}\cline{2-3}\multicolumn{1}{c|}{}\\
\hdashline\multicolumn{1}{c|}{}\\
$\vdots$
\end{tabular}
\end{center}
Here, the inner shape of $P$ expands to include a box that is not in $P$. This means that, despite the fact that two boxes are added to $P$'s inner shape, only one entry (which is $(1,1)$) is added to $\mathfrak{q}$ in the subroutine. The $1$ takes the place of the $2$ in row $0$, putting $(2,0)$ into $\mathfrak{q}$. The $2$ is placed in row $1$ (which is empty) and gets placed in the first box to the right of $P$'s new inner shape (the box just to the right of the right red box above). However, $1$ is added in $Q$ into both boxes that were added to the inner shape of $P$ during the subroutine. Thus, we have:
\begin{center}
$P=$
\begin{tabular}{m{.2cm}m{.2cm}m{.2cm}m{.2cm}}
&&&$\vdots$\\
\hdashline&&\multicolumn{1}{c|}{}&\multicolumn{1}{c|}{1}\\
&&\multicolumn{1}{c|}{}&\multicolumn{1}{c|}{2}\\
\hdashline\cline{2-3}\cline{3-4}\multicolumn{1}{c|}{}&\multicolumn{1}{c|}{1}\\
&$\vdots$
\end{tabular}
\hspace{1cm}
$Q=$
\begin{tabular}{m{.2cm}m{.2cm}m{.2cm}m{.2cm}}
&&&$\vdots$\\
\hdashline&&\multicolumn{1}{c|}{}\\
\cline{2-3}\multicolumn{1}{c|}{}&1&\multicolumn{1}{c|}{1}\\
\hdashline\cline{2-3}\multicolumn{1}{c|}{}\\
$\vdots$
\end{tabular}
\end{center}
The reader may check that the following is the end result of the algorithm:
\begin{center}
$P=$
\begin{tabular}{m{.2cm}m{.2cm}m{.2cm}m{.2cm}}
&&&$\vdots$\\
\hdashline&&\multicolumn{1}{c|}{}&\multicolumn{1}{c|}{1}\\
&&\multicolumn{1}{c|}{}&\multicolumn{1}{c|}{2}\\
\hdashline\cline{2-3}\cline{3-4}\multicolumn{1}{c|}{}&\multicolumn{1}{c|}{1}\\
&$\vdots$
\end{tabular}
\hspace{1cm}
$Q=$
\begin{tabular}{m{.2cm}m{.2cm}m{.2cm}m{.2cm}}
&&&$\vdots$\\
\hdashline&&\multicolumn{1}{c|}{}&\multicolumn{1}{c|}{2}\\
\cline{2-3}\multicolumn{1}{c|}{}&1&1&\multicolumn{1}{c|}{3}\\
\hdashline\cline{3-4}\multicolumn{1}{c|}{}&\multicolumn{1}{c|}{2}\\
&$\vdots$
\end{tabular}
\end{center}
The reader should note that $P \neq T$ and that $Q \neq U$. The boxes in $P$ are shifted down and to the right by one unit from the boxes in $T$; the boxes in $Q$ are shifter to the right by one unit from the boxes in $U$.
\end{example}

For $\lambda$, $P$, and $Q$ as in Theorem \ref{crsk}, we will construct $\CRSK^{-1}(((P,Q),\lambda))$ as follows:

\vspace{.5cm}
\begin{alg}[Inverse Cylindric RSK] \label{invcrskalg}
\begin{center}\end{center}
\begin{algorithmic}[1]
\vspace{.5cm}
Function InverseCylindricRSK(((tableau $P$, tableau $Q$), partition $\lambda$)): \Comment{$\lambda$ must be the outer shape of both $P$ and $Q$}.
\State $T := P$.
\State $\beta :=$ inner shape of $P$.
\State $U :=$ empty tableau with shape $\beta / \beta$.
\State $i :=$ largest entry that is present in any of the boxes in $Q$.
\While{true}:
\State $S :=$ set of all boxes that, in $Q$, contain $i$.
\State $(T,S') := \ReverseFullMulti((T,S))$. \Comment{That is, $T$ gets assigned to the tableau returned by Algorithm \ref{revfullmulti} and $S'$ gets assigned to the returned set of boxes.}
\State Put $i$ into $U$ in all boxes in $S'$; shrink the inner shape of $U$ to exclude all boxes in $S'$.
\If{$i$ is greater than the smallest entry that is present in any of the boxes in $Q$}:
\State $i :=$ largest entry less than $i$ that is present in any of the boxes in $Q$.
\Else:
\State Exit the loop.
\EndIf.
\EndWhile.
\State $\mu :=$ inner shape of $T$.
\State \Return $((T,U),\mu)$.
\end{algorithmic}
\end{alg}

For an example of our inverse cylindric RSK algorithm, the reader may look at Example \ref{bigex}, starting with the last two tableaux drawn as $P$ and $Q$, respectively. The pairs of tableaux obtained after each call of the subroutine (Algorithm \ref{revfullmulti}) and update of $U$ are the same as the pairs of tableaux in Example \ref{bigex}, except in reverse order. The algorithm ends with $T$ and $U$ equalling the $T$ and $U$ given at the beginning of Example \ref{bigex}, respectively.

\vspace{.5cm}
Now we will prove that Algorithm \ref{crskalg} and Algorithm \ref{invcrskalg} are valid algorithms for $\CRSK((T,U),\mu)$ and $\CRSK^{-1}((P,Q),\lambda)$, respectively.

Given any $\alpha, \beta, \mu \in \Cylpar$ such that $\mu \subseteq \alpha$ and $\mu \subseteq \beta$, as well as any $T \in \SSCT(\alpha / \mu)$ and $U \in \SSCT(\beta / \mu)$, let $\CylindricRSK(((T,U),\mu)) = ((P,Q),\lambda)$. Based on our theorem statement, in order to show that Algorithm \ref{crskalg} is a valid algorithm for the $\CRSK$ function, we need to show that:
\begin{itemize}
\item Every time that Algorithm \ref{fullmulti} is called from Algorithm \ref{crskalg}, all preconditions for Algorithm \ref{fullmulti} are met.
\item Algorithm \ref{crskalg} terminates;
\item $\lambda \in \Cylpar$;
\item $\alpha \subseteq \lambda$ and $\beta \subseteq \lambda$;
\item $P \in \SSCT(\lambda / \beta)$ and $Q \in \SSCT(\lambda / \alpha)$; and
\item $\wt(T) = \wt(P)$ and $\wt(U) = \wt(Q)$.
\end{itemize}
Given any $\alpha, \beta, \lambda \in \Cylpar$ such that $\alpha \subseteq \lambda$ and $\beta \subseteq \lambda$, as well as any $P \in \SSCT(\lambda / \beta)$ and $Q \in \SSCT(\lambda / \alpha)$, let $\InverseCylindricRSK(((P,Q),\lambda)) = ((T,U),\mu)$. Based on our theorem statement, in order to show that Algorithm \ref{invcrskalg} is a valid algorithm for the inverse of the $\CRSK$ function, we need to show that:
\begin{itemize}
\item Every time that Algorithm \ref{revfullmulti} is called from Algorithm \ref{invcrskalg}, all preconditions for Algorithm \ref{revfullmulti} are met.
\item Algorithm \ref{invcrskalg} terminates;
\item $\mu \in \Cylpar$;
\item $\mu \subseteq \alpha$ and $\mu \subseteq \beta$;
\item $T \in \SSCT(\alpha / \mu)$ and $U \in \SSCT(\beta / \mu)$; and
\item $\wt(T) = \wt(P)$ and $\wt(U) = \wt(Q)$.
\end{itemize}
Finally, we need to show that Algorithm \ref{invcrskalg} does indeed invert Algorithm \ref{crskalg}. That is, we need to show that, for a valid input $((T,U),\mu)$ into Algorithm \ref{crskalg},
\begin{equation*}
\InverseCylindricRSK(\CylindricRSK(((T,U),\mu))) = ((T,U),\mu),
\end{equation*}
and that, for a valid input $((P,Q),\lambda)$ into Algorithm \ref{invcrskalg},
\begin{equation*}
\CylindricRSK(\InverseCylindricRSK(((P,Q),\lambda))) = ((P,Q),\lambda).
\end{equation*}

We will prove these facts below.

\begin{remark}
Every time that Algorithm \ref{fullmulti} is called from Algorithm \ref{crskalg}, all preconditions for Algorithm \ref{fullmulti} are met.
\end{remark}

\begin{proof}
We will proceed by induction. We know that the original $P$ (which equals $T$) is a valid semistandard tableau. We also know that $\mu$ plus the boxes in $U$ that contain the smallest entry of $U$ (call it $m$) is a valid partition (because otherwise there would be a box in $U$ containing $m$ that would be below or to the right of a box containing an entry less than $m$). Furthermore, any set of boxes that, in a tableau, contain the same entry, forms a horizontal strip (since otherwise there would be an entry below an entry of the same value in the tableau). Thus, the preconditions for Algorithm \ref{fullmulti} are satisfied the first time that the subroutine is called.

Suppose that, up to but not including the call of Algorithm \ref{fullmulti} when $i = j$ (for some $j > m$), the preconditions of the algorithm are satisfied. We know, then, that $P$ is a valid semistandard tableau at the time of the call when $i = j$ (because Algorithm \ref{fullmulti} returns valid semistandard tableaux). We know that $S$ forms a horizontal strip of boxes (by the same reasoning as above).

Let $\mu_j$ be the inner shape of $P$ when $i = j$ but before Algorithm \ref{fullmulti} is called in that iteration of the loop. We know that $\mu_j$ is a valid partition. Suppose, for contradiction, that adding to $\mu_j$ the boxes in $S$ will result in an invalid partition. Let $\mu'_j$ be the result when these boxes are added (by our assumption, $\mu'_j$ is not a valid partition). Then one of the boxes in $S$ (call it $B$) has a box above it or to its left (call it $C$) that is neither in $\mu_j$, nor in $S$.

The boxes that are not in $U$ fall into two categories: boxes in $\mu$ and boxes not in $\beta$. We know that any box above or to the left of $B$ that is not in $U$ falls into the first category. Thus, if $C$ is not in $U$, $C \in \mu$, and $\mu \subseteq \mu_j$. It follows that $C \in \mu_j$, a contradiction. If $C$ is in $U$ but is not in $\mu_j$, then, in $U$, $C$ contains an entry that is greater than or equal to $j$. We know that $C$ cannot contain $j$ because then it would have been in $S$. This means that $C$ contains an entry greater than $j$ in $U$. This is a contradiction, as $B$ contains $j$ and $C$ is above or to the left of $B$. Therefore, adding $S$ to $\mu_j$ results in a valid partition. This means that all preconditions for Algorithm \ref{fullmulti} are satisfied.

Having completed our induction, we have shown that all preconditions for Algorithm \ref{fullmulti} are met every time that the algorithm is called from Algorithm \ref{crskalg}.
\end{proof}

\begin{remark}
Algorithm \ref{crskalg} terminates.
\end{remark}

\begin{proof}
By Remark \ref{fmtermin}, we know that the subroutine called from Algorithm \ref{crskalg} always terminates. The loop in Algorithm \ref{crskalg} terminates as well, as $Q$ has a finite number of boxes (since $\lambda / \mu$ has a finite number of boxes for all $\mu \subseteq \lambda$) and therefore has a finite number of distinct entries (values that $i$ can take on). Therefore, Algorithm \ref{crskalg} always terminates.
\end{proof}

\begin{remark}
$\lambda \in \Cylpar$.
\end{remark}

\begin{proof}
Since the subroutine (Algorithm \ref{fullmulti}) always returns a valid tableau, the final tableau returned by the subroutine (which is the final value of $P$) is bounded by a valid outer shape. This is the partition $\lambda$ that is returned by the algorithm.
\end{proof}

\begin{remark}
$\alpha \subseteq \lambda$ and $\beta \subseteq \lambda$.
\end{remark}

\begin{proof}
Suppose, for contradiction, that $\alpha \not \subseteq \lambda$. Then there is a box that is in $\alpha$, but is not in $\lambda$. Since at the beginning of the algorithm, $\alpha$ is the outer shape of $P$, and at the end, $\lambda$ is the outer shape of $P$, it follows that the outer shape of $P$ lost at some point during the algorithm. $P$ is only modified in the subroutine, but in the subroutine the outer shape of the tableau can only expand. This is a contradiction. Therefore, $\alpha \subseteq \lambda$.

Suppose, for contradiction, that $\beta \not \subseteq \lambda$. Then $P$ would not be a valid tableau. This is a contradiction, as $P$ is returned by Algorithm \ref{fullmulti}, which always returns valid tableaux. Therefore, $\beta \subseteq \lambda$.
\end{proof}

\begin{remark} \label{pqvalid}
$P \in \SSCT(\lambda / \beta)$ and $Q \in \SSCT(\lambda / \alpha)$.
\end{remark}

\begin{proof}
The outer shape of $P$ is $\lambda$ by definition. The inner shape of $P$ is $\mu$, plus the boxes of $U$ (i.e. the boxes in $\beta / \mu$). Therefore, the inner shape of $P$ is $\beta$. $P$ is a semistandard tableau, as it is returned by Algorithm \ref{fullmulti}, which always returns semistandard tableaux. Therefore, $P \in \SSCT(\lambda / \beta)$.

The boxes of $Q$ are the boxes that are in $P$ but not in $T$ (since they are the boxes that were added to $P$ after every call of Algorithm \ref{fullmulti}). Since $T \in \SSCT(\alpha / \mu)$ and $P \in \SSCT(\lambda / \beta)$, $Q$ consists of all boxes that are in $\lambda$ but not in $\alpha$ --- that is, $Q$'s inner shape is $\alpha$ and its outer shape is $\lambda$.

Suppose, for contradiction, that $Q$ is not a semistandard tableau. Then there is a box $B$ in $Q$ that contains an entry (call it $b$) that is below a box that contains $b$, or is below or to the right of a box that contains an entry greater than $b$.

The first case would imply that $S'$ does not form a horizontal strip. This is a contradiction, by Corollary \ref{horstrip}. In the second case, let $C$ be a box that is above or to the left of $B$ that contains an entry greater than $b$. Consider the tableau that is returned by Algorithm \ref{fullmulti} during the iteration of the loop of Algorithm \ref{crskalg} when $i = b$; call this tableau $R$. $B$ is in $S'$ (and is therefore in the outer shape of $R$), but $C$ gets added to the outer shape of $P$ on a later call of the loop, and is therefore not in the outer shape of $R$. This means that the outer shape of $R$ is not a valid partition. This is a contradiction. Therefore, $Q$ is semistandard. It follows that $Q \in \SSCT(\lambda / \alpha)$.
\end{proof}

\begin{remark}
$\wt(T) = \wt(P)$ and $\wt(U) = \wt(Q)$.
\end{remark}

\begin{proof}
Algorithm \ref{fullmulti} clearly preserves the weight of a tableau. It follows that $\wt(T) = \wt(P)$. In any given iteration of the loop in Algorithm \ref{crskalg}, $\left|S\right|$ is the number of times that $i$ appears in $U$, and $\left|S'\right|$ is the number of times that $i$ appears in $Q$. Clearly, $\left|S\right| = \left|S'\right|$. Therefore, $\wt(U) = \wt(Q)$.
\end{proof}

We have proven all necessary facts about Algorithm \ref{crskalg}. The proofs of the facts outlined above about Algorithm \ref{invcrskalg} have entirely analogous proofs, and, for this reason, are not provided in this paper.

Next, we show that our algorithms are indeed inverses of each other.

\begin{remark}
$\InverseCylindricRSK(\CylindricRSK(((T,U),\mu))) = ((T,U),\mu)$
\end{remark}

\begin{proof}
Let $((P,Q),\lambda) = \CylindricRSK(((T,U),\mu))$. We first note that $((P,Q),\lambda)$ satisfies the input preconditions for Algorithm \ref{invcrskalg}: $P$ and $Q$ are both valid semistandard tableaux with outer shape $\lambda$ (Remark \ref{pqvalid}).

Define $P_j$ and $Q_j$ to be the values of $P$ and $Q$, respectively, during the execution of Algorithm \ref{crskalg} on $((T,U),\mu)$ while $i = j$, but \emph{before} Algorithm \ref{fullmulti} is called inside the loop. Define $T_j$ and $U_j$ to be the values of $T$ and $U$, respectively, during the execution of Algorithm \ref{invcrskalg} on $((P,Q),\lambda)$ while $i = j$, \emph{after} Algorithm \ref{revfullmulti} is called inside the loop. Define $S_j$ to be the set of boxes that, in $U$, contain $j$, and $N_j$ to be the set of boxes that, in $Q$, contain $j$.

Suppose that, for all $j > m$ for some $m$, $P_j = T_j$ and $U_j$ is the tableau of boxes that, in $U$, contain entries greater than or equal to $j$, with those boxes holding the same entries as they do in $U$. We will show that $P_m = T_m$ and that $U_m$ is the tableau of boxes that, in $U$, contain entries greater than or equal to $m$, with those boxes holding the same entries as they do in $U$.

If $m$ is the largest value that $i$ takes on (that is, it is the largest entry in any box in $U$), let $V = P$ and $W = Q$. Otherwise, let $M$ be the smallest entry in $U$ that is greater than $m$ (that is, it is the value that $i$ takes on after $m$ in Algorithm \ref{crskalg}, and it is the value that $i$ takes on before $m$ in Algorithm \ref{invcrskalg}); let $V = T_M$ and $W = U_M$. We know that $\FullMulti(P_m,S_m) = V$, because we know that the tableau returned by Algorithm \ref{fullmulti} is $P_M = T_M = V$ (or is $V$ if $m$ is the largest value that $i$ takes on). Since $N_m$ is the set of boxes that, in $Q$, contain $m$, it follows that $\FullMulti((P_m,S_m)) = (V,N_m)$. $V$ is the tableau from which, during Algorithm \ref{invcrskalg}, the boxes that, in $Q$, contain $m$ are removed (in the iteration of the loop following the one in which $V$ is defined). Thus, $\ReverseFullMulti((V,N_m)) = (T_m,X)$ for some set of boxes $X$. By Theorem \ref{inverse}, it follows that $T_m = P_m$ and that $X = S_m$. From the latter fact it follows that $m$ is placed in $U_m$ in the boxes of $S_m$. By definition, $S_m$ is the set of boxes that, in $U$, contain $m$. Since $W$ is the tableau of boxes that, in $U$, contains entries \emph{greater} than $m$, with those boxes holding the same entries as they do in $U$ (by inductive hypothesis or trivially in the case that $m$ is the largest value that $i$ takes on), it follows that $U_m$ is the tableau of boxes that, in $U$, contains entries greater than or equal to $m$, with those boxes holding the same entries as they do in $U$. Our induction being complete, we have proven that $P_t = T_t$ and $U_t = U$, where $t$ is the smallest value that $i$ takes on. Since $T_t$ and $U_t$ are the final values that $T$ and $U$ take on in Algorithm \ref{invcrskalg}, $P_t = T$ in Algorithm \ref{crskalg}, and $\mu$ is the inner shape of $T$, it follows that $((T,U),\mu)$ is returned by Algorithm \ref{invcrskalg}, as desired.
\end{proof}

The proof that $\CylindricRSK(\InverseCylindricRSK(((P,Q),\lambda))) = ((P,Q),\lambda)$ is entirely analogous and, for this reason, is not provided in this paper.

Having proven the validity of the outputs of Algorithms \ref{crskalg} and \ref{invcrskalg}, the outputs' compliance with the theorem statement, and the fact that the algorithms are inverses of each other, we have shown the existence of a cylindric RSK correspondence.\footnote{Having proven the validity of $\CylindricRSK$ and $\InverseCylindricRSK$, we can now refer to them as they are referred to in the statement of Theorem \ref{crsk}: $\CRSK$ and $\CRSK^{-1}$, respectively.}
\end{proof}

\subsection{Consequences of the Cylindric RSK Correspondence}
Just as the RSK correspondence can be used in order to prove many results about regular tableaux, the cylindric RSK correspondence can be used in order to prove results about cylindric tableaux. Some such results are shown below.

\begin{theorem}[Cylindric Cauchy Identity] \label{schurpols}
Given two cylindric partitions $\alpha$ and $\beta$, we have
\begin{equation} \label{schureq}
\sum \limits_{\substack{\mu \in \Cylpar; \\ \mu \subseteq \alpha; \mu \subseteq \beta}} s_{\alpha / \mu}(\mathbf{x}) s_{\beta / \mu}(\mathbf{y}) = \sum \limits_{\substack{\lambda \in \Cylpar; \\ \alpha \subseteq \lambda; \beta \subseteq \lambda}} s_{\lambda / \beta}(\mathbf{x}) s_{\lambda / \alpha}(\mathbf{y}).
\end{equation}
where $\mathbf{x}$ and $\mathbf{y}$ are variable sets.
\end{theorem}

\begin{proof}
It is easy to see that
\begin{equation} \label{schur1}
\sum \limits_{\substack{\mu \in \Cylpar; \\ \mu \subseteq \alpha; \mu \subseteq \beta}} s_{\alpha / \mu}(\mathbf{x}) s_{\beta / \mu}(\mathbf{y}) = \sum \limits_{\substack{\mu \in \Cylpar; \\ \mu \subseteq \alpha; \mu \subseteq \beta; \\ T \in \SSCT(\alpha / \mu); \\ U \in \SSCT(\beta / \mu)}} \mathbf{x}^{\wt(T)} \mathbf{y}^{\wt(U)}
\end{equation}
and that
\begin{equation} \label{schur2}
\sum \limits_{\substack{\lambda \in \Cylpar; \\ \alpha \subseteq \lambda; \beta \subseteq \lambda}} s_{\lambda / \beta}(\mathbf{x}) s_{\lambda / \alpha}(\mathbf{y}) = \sum \limits_{\substack{\lambda \in \Cylpar; \\ \alpha \subseteq \lambda; \beta \subseteq \lambda; \\ P \in \SSCT(\lambda / \beta); \\ Q \in \SSCT(\lambda / \alpha)}} \mathbf{x}^{\wt(P)} \mathbf{y}^{\wt(Q)}.
\end{equation}
The cylindric RSK correspondence establishes a bijection between such pairs $(T,U)$ and $(P,Q)$ as in equations \ref{schur1} and \ref{schur2}. Furthermore, the bijection takes a pair $(T,U)$ to a pair $(P,Q)$ such that $\wt(T) = \wt(P)$ and $\wt(U) = \wt(Q)$, meaning that $\mathbf{x}^{\wt(T)} = \mathbf{x}^{\wt(P)}$ and $\mathbf{y}^{\wt(U)} = \mathbf{y}^{\wt(Q)}$. Therefore,
\begin{equation*}
\sum \limits_{\substack{\mu \in \Cylpar; \\ \mu \subseteq \alpha; \mu \subseteq \beta; \\ T \in \SSCT(\alpha / \mu); \\ U \in \SSCT(\beta / \mu)}} \mathbf{x}^{\wt(T)} \mathbf{y}^{\wt(U)} = \sum \limits_{\substack{\lambda \in \Cylpar; \\ \alpha \subseteq \lambda; \beta \subseteq \lambda; \\ P \in \SSCT(\lambda / \beta); \\ Q \in \SSCT(\lambda / \alpha)}} \mathbf{x}^{\wt(P)} \mathbf{y}^{\wt(Q)}.
\end{equation*}
Our theorem statement is therefore proven, by the transitive property of equality.
\end{proof}

\begin{defin}
A \textit{standard cylindric tableau} is a semistandard cylindric tableau whose alphabet is ${1,2, \dots,m}$ for some nonnegative integer $m$, such that each entry of the tableau's alphabet appears in exactly one box of the tableau. An example of a standard cylindric tableau is shown below.
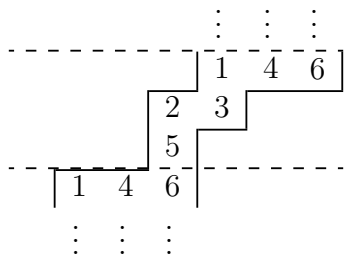
\begin{center}
\begin{tabular}{m{.2cm}m{.2cm}m{.2cm}m{.2cm}m{.2cm}m{.2cm}m{.2cm}}
&&&&$\vdots$&$\vdots$&$\vdots$\\
\hdashline&&&\multicolumn{1}{c|}{}&1&4&\multicolumn{1}{c|}{6}\\
\cline{4-4}\cline{6-7}&&\multicolumn{1}{c|}{}&2&\multicolumn{1}{c|}{3}\\
\cline{5-5}&&\multicolumn{1}{c|}{}&\multicolumn{1}{c|}{5}\\
\hdashline\cline{2-3}\multicolumn{1}{c|}{}&1&4&\multicolumn{1}{c|}{6}\\
&$\vdots$&$\vdots$&$\vdots$
\end{tabular}
\end{center}
\end{defin}

\begin{defin}
Given two cylindric partitions $\lambda$ and $\mu$, $f_{\lambda / \mu}$ will denote the number of standard tableaux with shape $\lambda / \mu$.
\end{defin}

\begin{cor} \label{fcount}
Let $\alpha$ and $\beta$ be two cylindric partitions and $m$ be a nonnegative integer. Let $M$ be the set of all partitions $\mu$ such that $\mu \subseteq \alpha$, $\mu \subseteq \beta$, and $\alpha / \mu$ contains $m$ distinct boxes. Let $\Lambda$ be the set of all partitions $\lambda$ such that $\alpha \subseteq \lambda$, $\beta \subseteq \lambda$, and $\lambda / \beta$ contains $m$ distinct boxes. Then
\begin{equation} \label{standardcount}
\sum \limits_{\mu \in M} f_{\alpha / \mu} f_{\beta / \mu} = \sum \limits_{\lambda \in \Lambda} f_{\lambda / \alpha} f_{\lambda / \beta}.
\end{equation}
\end{cor}

\begin{proof}
Let $a$ be the number of boxes that are in $\alpha$, but not in $\beta$; let $b$ be the number of boxes that are in $\beta$, but not in $\alpha$. Then the number of boxes in $\beta / \mu$ is $m - a + b$. Let $p = m - a + b$. Similarly, the number of boxes in $\lambda / \alpha$ is $p$.

For any standard tableau $T$ with shape $\alpha / \mu$, we have $\mathbf{x}^{\wt(T)} = x_1 x_2 \dots x_m$. Similarly, for any standard tableau $U$ with shape $\beta / \mu$, we have $\mathbf{y}^{\wt(U)} = y_1 y_2 \dots y_p$.

The left hand side of (\ref{standardcount}) is the number of pairs of standard tableaux $(T,U)$ with a common inner shape and outer shapes of $\alpha$ and $\beta$, respectively. It follows that this is the coefficient of $x_1 x_2 \dots x_m \cdot y_1 y_2 \dots y_p$ on the left hand side of (\ref{schureq}).

Similarly, we have that the right hand side of (\ref{standardcount}) is the coefficient of $x_1 x_2 \dots x_m \cdot y_1 y_2 \dots y_p$ on the right hand side of (\ref{schureq}). Therefore, by Theorem \ref{schurpols}, the left hand side and the right hand side of (\ref{standardcount}) have the same value, as desired.
\end{proof}

\subsection{The Symmetry Property of CRSK}
We now prove an important property of $\CRSK$, from which we can then derive more identities.

\begin{theorem} \label{symmetry}
Given $T$, $U$, and $\mu$ that satisfy the preconditions of $\CRSK(((T,U),\mu))$, let $\CRSK(((T,U),\mu)) = ((P,Q),\lambda)$. Then $\CRSK(((U,T),\mu)) = ((Q,P),\lambda)$.
\end{theorem}

\begin{proof}
We first note that it is clear that if the preconditions of $\CRSK(((T,U),\mu))$ are met, then the preconditions of $\CRSK(((U,T),\mu))$ are also met.

We will prove this theorem $T$ and $U$ with alphabets $(1,2,\dots,m_T)$ and $(1,2, \dots,m_U)$, respectively; the result for tableaux with other alphabets follows directly by analogy. We will also assume, without loss of generality, that all elements of $\{1,2,\dots,m_T\}$ appear in $T$ and that all elements of $\{1,2,\dots,m_U\}$ appear in $U$ (this can be assumed because otherwise the entries in $T$ or $U$ can be ``compressed" and $m_T$ or $m_U$ reduced).

Define $T_0$ to be $T$ and $T_j$, for $1 \le j \le m_U$, to be the value of $P$ after line \ref{callsub} of Algorithm \ref{crskalg} is performed (with input $((T,U),\mu)$), when $i = j$. Define $U_0$ to be $U$ and $U_j$, for $1 \le j \le m_T$, to be the value of $P$ after line \ref{callsub} of Algorithm \ref{crskalg} is performed (with input $((U,T),\mu)$), when $i = j$.

Define $\lambda_{i,j}$, for $0 \le i \le m_U$ and $0 \le j \le m_T$, to be the outer shape of the tableau with the same inner shape as $T_i$ and containing exactly the boxes that, in $T_i$, contain entries that are less than or equal to $j$. (That is, $\lambda_{i,j}$ is the right bound on the entries $1$ through $j$ in $T_i$.) Define $\nu_{i,j}$, for $0 \le i \le m_T$ and $0 \le j \le m_U$, to be the outer shape of the tableau with the same inner shape as $U_i$ and containing exactly the boxes that, in $U_i$, contain entries that are less than or equal to $j$.

We will perform induction on two statements, for all $1 \le c \le m_U$ and $1 \le d \le m_T$:
\begin{enumerate}[label=(\arabic*)]
\item Let $B$ be a box. Then $d$ is bumped from $B$ during the formation of $T_c$ if and only if $c$ is bumped from $B$ during the formation of $U_d$; and ~\label{corbump}
\item $\lambda_{c,d} = \nu_{d,c}$. ~\label{lanu}
\end{enumerate}

Suppose that $d = 0$. Then $\lambda_{c,d}$ is the inner shape of $T_c$. Based on how $\CRSK$ works, we know that $\lambda_{c,0}$ is the outer shape of the tableau with inner shape that of $U$, containing exactly the boxes that, in $U$, contain entries that are less than or equal to $c$. It is easy to see that this equals $\nu_{0,c} = \nu_{d,c}$. Similarly, if $c = 0$, then $\lambda_{c,d} = \nu_{d,c}$. These two facts serve as our base cases (as it will turn out, base cases for statement \ref{corbump} are not necessary).

Suppose that statements \ref{corbump} and \ref{lanu} are true for all $d < j$, for some $1 \le j \le m_T$. We will show that they are true for $d = j$ as well.

Suppose that statements \ref{corbump} and \ref{lanu} are true for $d = j$ and all $c < i$, for some $1 \le i \le m_U$. We will show that they are true for $d = j$ and $c = i$ as well.

We will first prove statement \ref{corbump} for $i$ and $j$ --- that is, that, for any box $B$, $j$ is bumped from $B$ during the formation of $T_i$ if and only if $i$ is bumped from $B$ during the formation of $U_j$.

Choose a particular box $B$. We will first prove that if $j$ is bumped from $B$ during the formation of $T_i$, then $i$ is bumped from $B$ during the formation of $U_j$.

\textbf{Case 1: In $T$, $B$ contains $j$.} Since $j$ is bumped from $B$ during the formation of $T_i$ (and is replaced by an entry that is less than $j$), we know that $B \in \lambda_{i,j-1} / \lambda_{i-1,j-1}$. It follows from our inductive hypothesis for statement \ref{lanu} that $B \in \nu_{j-1,i} / \nu_{j-1,i-1}$, which means that, in $U_{j-1}$, $B$ contains $i$. By how $\CRSK$ works, we know that an entry is bumped out of $B$ during the formation of $U_j$ (since $T$ contains $j$ in $B$). Thus, $i$ is bumped from $B$ during the formation of $U_j$, as desired.

\textbf{Case 2: In $T$, $B$ does not contain $j$.} $B$ contains $j$ in $T_{i-1}$, so $B \in \lambda_{i-1,j} / \lambda_{i-1,j-1}$. By our inductive hypothesis on statement \ref{lanu}, it follows that $B \in \nu_{j,i-1} / \nu_{j-1,i-1}$. Thus, an entry that is greater than or equal to $i$ is bumped out of $B$ during the formation of $U_j$. (From Corollary \ref{boxonce}, we know that only one entry is bumped out of $B$ during the formation of $U_j$.) It remains to show that this entry is indeed $i$.

Suppose that in $U$, $B$ contains $i$. Then in $U_{j-1}$, $B$ contains an entry that is less than or equal to $i$. However, an entry that is greater than or equal to $i$ is bumped out of $B$ during the formation of $U_j$. Thus, $i$ must be the entry that is bumped out of $B$ during the formation of $U_j$, as desired.

Now suppose that in $U$, $B$ does not contain $i$. This means that, during the formation of $T_i$, $B$ is not one of the boxes that is originally removed in the subroutine (that is, $B \not \in S$). That is, $j$ is bumped out of $B$ by an entry $j''$ during the formation of $T_i$. This means that $B \in \lambda_{i,j''} / \lambda_{i-1,j''}$. By our inductive hypothesis on statement \ref{lanu}, this means that $B \in \nu_{j'',i} / \nu_{j'',i-1}$ --- that is, that $B$ contains $i$ in $U_{j''}$. For all $j'' < x < j$, $x$ is not bumped from $B$ during the formation of $T_i$; it follows from our inductive hypothesis on statement \ref{corbump} that for all $j'' < x < j$, $i$ is not bumped from $B$ during the formation of $U_x$. Therefore, $i$ is in $B$ in $U_{j-1}$. We conclude that $i$ is the entry that is bumped out of $B$ during the formation of $U_j$ in this case as well.

Note that we used our inductive hypothesis on statement \ref{corbump} once, and we used it for $c = i$ and $d = x$ for $j'' < x < j$ (that is, for defined values of $c$ and $d$); for this reason, we did not need a base case for this proof.

The converse of the statement that we have proven follows by symmetry. This is not immediately trivial, because our induction is not symmetric (we induct on $d$ and, for each $d$, induct on $c$); however, the reader may notice that, in all cases, the inductive hypothesis was used only for $c \le i$ and $d \le j$. Thus, the converse of the statement does indeed follow by symmetry.

We have shown statement \ref{corbump} to be true for $c = i$ and $d = j$. Now we show that statement \ref{lanu} is true for $c = i$ and $d = j$.

Consider the process of formation of $T_i$ and consider a particular row $r$. By Corollary \ref{insertincrease}, we know that all $j$'s that are bumped into $r$ are bumped into $r$ after all entries less than $j$ are bumped into $r$. It follows that all $j$'s that are bumped into $r$ are inserted just to the right of the rightmost of the following:
\begin{itemize}
\item The rightmost entry less-than-or-equal-to $j$ that is in $T_{i-1}$; or
\item The rightmost entry less than $j$ that is in $T_i$ (in the case that, during the formation of $T_i$, entries less than $j$ bump out all $j$'s that were in $T_{i-1}$).
\end{itemize}
That is, the $j$'s are placed just to the right of the rightmost of $\lambda_{i-1,j}$ and $\lambda_{i,j-1}$ in row $r$.

Similarly, during the formation of $U_j$, the $i$'s are placed in row $r$ just to the right of the rightmost of $\nu_{j-1,i}$ and $\nu_{j,i-1}$. By our inductive hypothesis, $\lambda_{i-1,j} = \nu_{j,i-1}$ and $\lambda_{i,j-1} = \nu_{j-1,i}$.

Since, for any box $B$ in row $r - 1$, $j$ is bumped from $B$ during the formation of $T_i$ if and only if $i$ is bumped from $B$ during the formation of $U_j$, it follows that the same number of $j$'s are inserted into $r$ during the formation of $T_i$ as the number of $i$'s that are inserted into $r$ during the formation of $U_j$. Since the leftmost of these insertion locations is the same (see the previous paragraph), we conclude that the $j$'s inserted during the formation of $T_i$ are inserted into the same boxes as the $i$'s inserted during the formation of $U_j$. It follows that $\lambda_{i,j} = \nu_{j,i}$, as desired.

We have completed our induction on $i$ and have shown that statements \ref{corbump} and \ref{lanu} are true for $d = j$ and all values of $c$.

We have completed our induction on $j$ and have shown that statements \ref{corbump} and \ref{lanu} are true for all values $d$ and all values of $c$.

Let $\CRSK(((T,U),\mu)) = ((P,Q),\lambda)$ and $\CRSK(((U,T),\mu)) = ((Q',P'),\lambda')$. We will show that $P = P'$ (it will follow by analogy that $Q = Q'$, and it will follow that $\lambda = \lambda'$).

Define $\gamma_i$, for $1 \le i \le m_T$, to be the outer shape of the tableau with the same inner shape as $P$ and containing exactly the boxes that, in $P$, contain entries that are less than or equal to $i$; define $\gamma_0$ to be the inner shape of $P$. Define $\gamma'_i$ analogously for $P'$. It suffices to show that, for all $i$, $\gamma_i = \gamma'_i$.

$\gamma_0$ and $\gamma'_0$ are both the outer shape of $U$, so clearly $\gamma_0 = \gamma'_0$. For any $1 \le i \le m_T$, $\gamma_i = \lambda_{m_U,i}$ (since $P = T_{m_U}$). By how $\CRSK$ works, we know that, for any $1 \le i \le m_T$, $\gamma'_i$ is the outer shape of $U_i$. Thus, $\gamma'_i = \nu_{i,m_U}$. By statement \ref{lanu}, we know that $\lambda_{m_U,i} = \nu_{i,m_U}$. Therefore, $\gamma_i = \gamma'_i$, as desired.

Therefore, for any two tableaux $T$ and $U$ and partition $\mu$ that satisfy the preconditions of $\CRSK(((T,U),\mu))$ and $\CRSK((U,T),\mu)$, if $\CRSK(((T,U),\mu)) = ((P,Q),\lambda)$, then $\CRSK(((U,T),\mu)) = ((Q,P),\lambda)$, as desired.
\end{proof}

\begin{cor}
Given two tableaux $P$ and $Q$ and a partition $\lambda$ that satisfy the preconditions of $\CRSK^{-1}(((P,Q),\lambda))$ and $\CRSK^{-1}((Q,P),\lambda)$, let $\CRSK^{-1}(((P,Q),\lambda)) = ((T,U),\mu)$. Then $\CRSK^{-1}(((Q,P),\lambda)) = ((U,T),\mu)$.
\end{cor}

\begin{proof}
We have that $\CRSK(((T,U),\mu)) = ((P,Q),\lambda)$. By Theorem \ref{symmetry}, we have that $\CRSK(((U,T),\mu)) = ((Q,P),\lambda)$. Therefore, $\CRSK^{-1}(((Q,P),\lambda)) = ((U,T),\mu)$, as desired.
\end{proof}

\begin{cor} \label{tabtotab}
Given a tableau $T$ with inner shape $\mu$, $\CRSK(((T,T),\mu)) = ((P,P),\lambda)$, for some $P$ and $\lambda$.
\end{cor}

\begin{proof}
Let $\CRSK(((T,T),\mu)) = ((P,Q),\lambda)$. We are to show that $P = Q$. By Theorem \ref{symmetry}, we have $\CRSK(((T,T),\mu)) = ((Q,P),\lambda)$. Thus, $P = Q$, as desired.
\end{proof}

Corollary \ref{tabtotab} shows that $\CRSK$ can be viewed as a bijection that takes a tableau with outer shape $\alpha$ to a tableau with inner shape $\alpha$ with the weight. This allows us to obtain some more important results:

\begin{cor} \label{oneschurcor}
Given any cylindric partition $\alpha$, we have
\begin{equation} \label{oneschur}
\sum \limits_{\substack{\mu \in \Cylpar; \\ \mu \subseteq \alpha}} s_{\alpha / \mu}(\mathbf{x}) = \sum \limits_{\substack{\lambda \in \Cylpar; \\ \alpha \subseteq \lambda}} s_{\lambda / \alpha}(\mathbf{x}).
\end{equation}
\end{cor}

It is worth mentioning that setting $\beta := \alpha$ and $\mathbf{x} = \mathbf{y}$ in Theorem \ref{schurpols} gives us the following related (but different) result:

\begin{equation*}
\sum \limits_{\substack{\mu \in \Cylpar; \\ \mu \subseteq \alpha}} s^2_{\alpha / \mu}(\mathbf{x}) = \sum \limits_{\substack{\lambda \in \Cylpar; \\ \alpha \subseteq \lambda}} s^2_{\lambda / \alpha}(\mathbf{x}).
\end{equation*}

We can use (\ref{oneschur}) in order to obtain the following result, analogous to Corollary \ref{fcount}:

\begin{cor} \label{onefcor}
Let $\alpha$ be a cylindric partition and $m$ be a nonnegative integer. Let $M$ be the set of all partitions $\mu$ such that $\mu \subseteq \alpha$ and $\alpha / \mu$ contains $m$ distinct boxes. Let $\Lambda$ be the set of all partitions $\lambda$ such that $\alpha \subseteq \lambda$ and $\lambda / \alpha$ contains $m$ distinct boxes. Then
\begin{equation*}
\sum \limits_{\mu \in M} f_{\alpha / \mu} = \sum \limits_{\lambda \in \Lambda} f_{\lambda / \alpha}.
\end{equation*}
\end{cor}

Relatedly, we can set $\beta := \alpha$ in (\ref{standardcount}) to obtain:

\begin{equation*}
\sum \limits_{\mu \in M} f^2_{\alpha / \mu} = \sum \limits_{\lambda \in \Lambda} f^2_{\lambda / \alpha}.
\end{equation*}

\section{A Marble-Game Interpretation of Cylindric Tableaux}
\begin{defin}
Consider any partition $\alpha$. The \textit{marble arrangement} of $\alpha$ (denoted $\Arr(\alpha)$) is constructed as follows. We begin with $k$ people --- call them $p_0$, $p_1$, $\dots$, $p_{k - 1}$ --- in a circle ($k$ being the vertical period of the cylinder) such that $p_i$ is clockwise from $p_{i - 1}$ (note that, just as with rows, $p_i$ refers to $p_{i \pmod{k}}$, and thus $p_0$ is clockwise from $p_{k - 1}$). For all $i$, $p_i$ has $\alpha_{i - 1} - \alpha_i$ marbles ($\alpha_i$ is the $i$'th term in the $\alpha$'s sequence).\footnote{Clearly, there is a total of $n - k$ marbles ($n - k$ is the horizontal period of the cylinder).}
\end{defin}

Let $R$ be a cylindric tableau with alphabet $\{1,2,\dots,t\}$, inner shape $\mu$, and outer shape $\lambda$. Define $\Arr_0(R) := \Arr(\mu)$. Next, for $1 \le j \le t$, obtain $\Arr_j(R)$ from $\Arr_{j-1}(R)$ as follows: for every row $r$, let $x$ be the number of times $j$ appears in row $r$ in $R$. Then $p_r$ in $\Arr_{j-1}(R)$ passes $x$ marbles to $p_{r + 1}$.

For $0 \le j \le t$, let $\lambda_j(R)$ be the outer shape of the tableau with inner shape $\mu$ and containing exactly the boxes that, in $R$, contains entries less than or equal to $j$.

\begin{remark}
For $0 \le j \le t$, $\Arr_j(R) = \Arr(\lambda_j(R))$.
\end{remark}

\begin{proof}
We will proceed by induction. We know that $\lambda_0(R) = \mu$, so we have $\Arr_0(R) = \Arr(\lambda_0(R))$. Suppose that $\Arr_j(R) = \Arr(\lambda_j(R))$ for all $0 \le j < m$. We will show that $\Arr_m(R) = \Arr(\lambda_m(R))$.

Consider any integer $r$. Say that $m$ appears $x_1$ times in row $r$ (the projection of plane row $r$ onto the cylinder) and $x_2$ times in row $r - 1$. Then $(\lambda_m(R))_r - (\lambda_{m - 1}(R))_r = x_1$ and $(\lambda_m(R))_{r - 1} - (\lambda_{m - 1}(R))_{r - 1} = x_2$. Subtracting these two equations, we have

\begin{equation*}
(\lambda_m(R))_{r - 1} - (\lambda_m(R))_r = (\lambda_{m - 1}(R))_{r - 1} - (\lambda_{m - 1}(R))_r + x_2 - x_1.
\end{equation*}

Therefore, we have that $p_r$ has $x_2 - x_1$ more marbles in $\Arr(\lambda_m(R))$ than he did in $\Arr(\lambda_{m - 1}(R))$.

When $\Arr_m(R)$ is obtained from $\Arr_{m - 1}(R)$, $p_r$ passes $x_1$ marbles to $p_{r+1}$ and receives $x_2$ marbles from $p_{r - 1}$. Thus, we have that $p_r$ has $x_2 - x_1$ more marbles in $\Arr_m(R)$ than he did in $\Arr_{m - 1}(R)$. Since this is true for all integers $r$, and by our inductive hypothesis, we have that $\Arr_m(R) = \Arr(\lambda_m(R))$, as desired. Having completed our induction, we have proven that for $0 \le j \le t$, $\Arr_j(R) = \Arr(\lambda_j(R))$.
\end{proof}

\begin{remark} \label{validgame}
For $1 \le j \le t$, the formation of $\Arr_j(R)$ never entails a person passing more marbles than he has to this right.
\end{remark}

\begin{proof}
The number of marbles that $p_r$ has, being the difference between the length of the row above $r$ and the length of row $r$ in the sub-tableau of $R$ consisting of the entries $1$ through $j - 1$, describes how many entries $j$ row $r$ can possibly have without the tableau being non-semistandard.
\end{proof}

\begin{defin}
A \textit{turn} is a combination of simultaneous marble passes among the $k$ people, such that each person passes at most as many marbles as he has. A turn is denoted $(a_0,a_1,\dots,a_{k - 1})$, where $a_i$ is the number of marbles that $p_i$ passes to $p_{i + 1}$.
\end{defin}

\begin{defin}
A \textit{marble game} of length $t$ is an initial arrangement of $n - k$ marbles among the $k$ people, along with a sequence of $t$ successive turns.
\end{defin}

\begin{prop} \label{marblebijection}
Given a partition $\mu$ and a nonnegative integer $t$, there is a bijection between marble games of length $t$ with initial arrangement $\Arr(\mu)$ and cylindric tableaux with alphabet $\{1,2,\dots,t\}$ and inner shape $\mu$.
\end{prop}

\begin{proof}
We have shown a construction of a marble game given a tableau, and, as stated in Remark \ref{validgame}, the game thus constructed is a valid one. Given a marble game $G$, one can obtain a unique tableau that produces $G$ with the process described above. This is because turn $j$ in $G$ describes the horizontal strip of $j$'s that is in its corresponding tableau: from $\lambda_{j - 1}(R)$ we obtain $\lambda_j(R)$ by adding to row $r$ the number of $j$'s equal to the number of marbles that $p_r$ gives to $p_{r + 1}$. Because $p_r$ cannot give more marbles than he has to $p_{r + 1}$, the tableau thus produced is indeed a valid semistandard tableau.
\end{proof}

\begin{example}
Let $R$ be the following cylindric tableau, with $t = 6$ (recall that the top row of the tableau is row $0$).

\begin{center}
$R =$
\begin{tabular}{m{.2cm}m{.2cm}m{.2cm}m{.2cm}m{.2cm}m{.2cm}m{.2cm}m{.2cm}m{.2cm}m{.2cm}}
&&&&&$\vdots$&$\vdots$&$\vdots$&$\vdots$&$\vdots$\\
\hdashline&&&&\multicolumn{1}{c|}{}&1&2&2&5&\multicolumn{1}{c|}{6}\\
\cline{5-5}\cline{10-10}&&&\multicolumn{1}{c|}{}&1&2&6&6&\multicolumn{1}{c|}{6}\\
\cline{3-4}\cline{7-9}&\multicolumn{1}{c|}{}&1&1&4&\multicolumn{1}{c|}{5}\\
\hdashline\cline{2-2}\multicolumn{1}{c|}{}&1&2&2&5&\multicolumn{1}{c|}{6}\\
&$\vdots$&$\vdots$&$\vdots$&$\vdots$&$\vdots$
\end{tabular}
\end{center}

Originally, we have $p_0 = 1$, $p_1 = 1$, and $p_2 = 2$. On the first turn, one $1$ is added to rows $0$ and $1$ and two $1$'s are added to row $2$; thus, after the first turn (represented as $(1,1,2)$), we have $p_0 = 2$, $p_1 = 1$, and $p_2 = 1$. The next turn is $(2,1,0)$, so after turn $2$, we have $p_0 = 0$, $p_1 = 2$, and $p_2 = 2$. There are no $3$'s in $R$, so the next turn --- $(0,0,0)$ --- results in $p_0 = 0$, $p_1 = 2$, and $p_2 = 2$. Turn $4$ is $(0,0,1)$; turn $5$ is $(1,0,1)$; and turn $6$ is $(1,3,0)$. Thus, the marble game corresponding to $R$ has the following sequence of turns:

\begin{equation*}
(1,1,2), (2,1,0), (0,0,0), (0,0,1), (1,0,1), (1,3,0).
\end{equation*}

If we know $\mu$, we can use this information to retrace $R$: the first turn encodes that $R$ has one $1$ in rows $0$ and $1$ and two $1$'s in row $2$; the second turn encodes that $R$ has two $2$'s in row $0$, one $2$ in row $1$, and no $2$'s in row $2$; and so on.
\end{example}

Note that knowing $\mu$ and $\Arr_j(R)$ for $1 \le j \le t$ does not uniquely determine $R$. For example, if $\Arr_1(R) = \Arr(\mu)$, one cannot necessarily determine whether there are no $1$'s in $R$, each row has one $1$ in $R$, each row has two $1$'s in $R$, etc. Thus, a game, and not just the series of arrangements attained by a game, is necessary in order to describe a tableau.

We can use Proposition \ref{marblebijection} in order to obtain some results relating the marble game and cylindric tableaux.

\begin{cor} \label{forward}
Let $\mu$ be a cylindric partition and $t$ be a nonnegative integer. The number of cylindric tableaux with inner shape $\mu$ and alphabet $\{1,2,\dots,t\}$ equals the number of possible marble games with $t$ turns that begins with the arrangement $\Arr(\mu)$.
\end{cor}

\begin{cor} \label{forwardstandard}
Let $\mu$ be a cylindric partition and $t$ be a nonnegative integer. The number of standard cylindric tableaux with inner shape $\mu$ and alphabet $\{1,2,\dots,t\}$ equals the number of possible marble games with $t$ turns that begins with the arrangement $\Arr(\mu)$ and in which exactly one marble changes hands on every turn.
\end{cor}

\begin{cor} \label{backward}
Let $\lambda$ be a cylindric partition and $t$ be a nonnegative integer. The number of cylindric tableaux with outer shape $\lambda$ and alphabet $\{1,2,\dots,t\}$ equals the number of possible marble games with $t$ turns that ends with the arrangement $\Arr(\lambda)$.
\end{cor}

Note that the number of possible marble games with $t$ turns that ends with the arrangement $\Arr(\lambda)$ is the same as the number of possible marble games with $t$ turns that begins with the arrangement $\Arr(\lambda)$, except that marbles are passed counterclockwise instead of clockwise. This is because this modified marble game retraces the steps of any marble game that ends with $\Arr(\lambda)$.

\begin{cor} \label{backwardstandard}
Let $\lambda$ be a cylindric partition and $t$ be a nonnegative integer. The number of standard cylindric tableaux with outer shape $\lambda$ and alphabet $\{1,2,\dots,t\}$ equals the number of possible marble games with $t$ turns that ends with the arrangement $\Arr(\lambda)$ and in which exactly one marble changes hands on every turn.
\end{cor}

\begin{cor}
Given a cylindric partition $\alpha$ and nonnegative integer $t$, the number of possible marble games with $t$ turns that starts with the arrangement $\Arr(\alpha)$ is equal to the number of possible marble games with $t$ turns that ends with the arrangement $\Arr(\alpha)$ (which is in turn equal to the number of possible marble games with $t$ turns that starts with the arrangement $\Arr(\alpha)$, but in which marbles are passed counterclockwise). This remains true if we restrict all turns so that exactly one marble changes hands on every turn.
\end{cor}

\begin{proof}
If in Corollary \ref{oneschurcor}, we let $\mathbf{x} = (x_1,x_2,\dots,x_t)$ and we let $x_1 = x_2 = \dots = x_t = 1$, we find that the number of cylindric tableaux with alphabet $\{1,2,\dots,t\}$ and inner shape $\alpha$ equals the number of cylindric tableaux with alphabet $\{1,2,\dots,t\}$ and outer shape $\alpha$. The first part of our proposition now follows directly from corollaries \ref{forward} and \ref{backward}. The second part of our proposition follows analogously from corollaries \ref{onefcor}, \ref{forwardstandard}, and \ref{backwardstandard}.
\end{proof}

It is noteworthy that this marble-game construction only works for cylindric tableaux. A similar concept can be defined for regular tableaux, but it would involve an infinite line of people in which the first person can obtain an arbitrary number of marbles --- a construction that is not nearly as interesting and much less likely to be of use. The purely combinatorial construction which we have described in this section shows the promise of cylindric tableaux in potential applications outside of tableau theory.

One such application may be in information theory, where marbles are bits and players are computers that are linked in a ring in which communication is only allowed in one direction. In such an interpretation, the numbers in a column represent the times when a particular bit is transferred from one computer to the next. The ease of tracking particular bits with tableau representations of communication in unidirectional rings makes it possible for such representations to be useful in optimizing such communication.

\section{Applying Results Concerning Cylindric Tableaux to Skew Tableaux} \label{cylskew}
At the beginning of the paper, we fixed $n$ and $k$. For this section, we will unfix $n$ and $k$.

It turns out that results concerning skew tableaux can be proven from analogous results concerning cylindric tableaux. The typical construction for such proofs is, given a skew tableau or shape, to create a cylindric tableau that looks like the skew tableau or shape, but with a very large $k$ (as large as necessary) and an even larger $n$ (so that $n - k$ is as large as necessary). As an example, given the skew shape below on the left, the cylindric shape that would be produced looks like the cylindric shape below on the right.

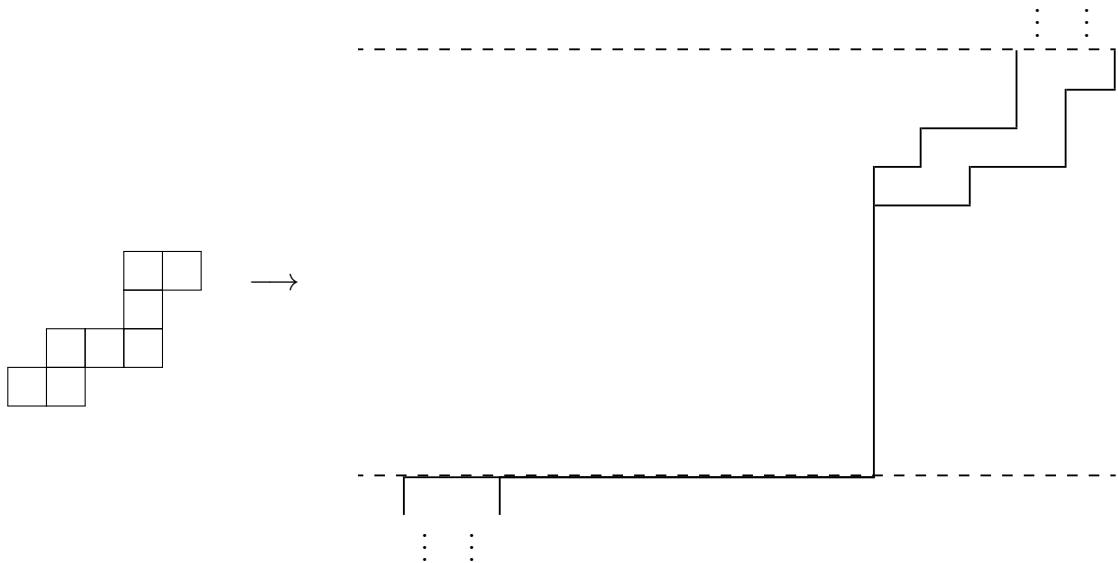
\begin{center}
\ytableausetup{boxsize=.5cm}
\ydiagram{3+2,3+1,1+3,2}
\hspace{.5cm}
$\longrightarrow$
\hspace{.5cm}
\begin{tabular}{m{.2cm}m{.2cm}m{.2cm}m{.2cm}m{.2cm}m{.2cm}m{.2cm}m{.2cm}m{.2cm}m{.2cm}m{.2cm}m{.2cm}m{.2cm}m{.2cm}m{.2cm}m{.2cm}}
&&&&&&&&&&&&&&$\vdots$&$\vdots$\\
\hdashline&&&&&&&&&&&&&\multicolumn{1}{c|}{}&\textcolor{white}{1}&\multicolumn{1}{c|}{\textcolor{white}{1}}\\
\cline{16-16}&&&&&&&&&&&&&\multicolumn{1}{c|}{}&\multicolumn{1}{c|}{\textcolor{white}{1}}\\
\cline{13-14}&&&&&&&&&&&\multicolumn{1}{c|}{}&\textcolor{white}{1}&\textcolor{white}{1}&\multicolumn{1}{c|}{\textcolor{white}{1}}\\
\cline{12-12}\cline{14-15}&&&&&&&&&&\multicolumn{1}{c|}{}&\textcolor{white}{1}&\multicolumn{1}{c|}{\textcolor{white}{1}}\\
\cline{12-13}&&&&&&&&&&\multicolumn{1}{c|}{}\\
&&&&&&&&&&\multicolumn{1}{c|}{}\\
&&&&&&&&&&\multicolumn{1}{c|}{}\\
&&&&&&&&&&\multicolumn{1}{c|}{}\\
&&&&&&&&&&\multicolumn{1}{c|}{}\\
&&&&&&&&&&\multicolumn{1}{c|}{}\\
&&&&&&&&&&\multicolumn{1}{c|}{}\\
\hdashline\cline{2-11}\cline{4-11}\multicolumn{1}{c|}{}&\textcolor{white}{1}&\multicolumn{1}{c|}{\textcolor{white}{1}}\\
&$\vdots$&$\vdots$
\end{tabular}
\end{center}

\begin{defin}
If a partition $\lambda$ is turned into a cylindric partition as shown above, where the vertical period of the cylindric partition is $k$ and the horizontal period of the cylindric partition is $n - k$, we denote the cylindric partition $\Cyl_{k,n}(\lambda)$. (Thus, if above our skew shape is $\lambda / \mu$, then the cylindric shape shown above is $\Cyl_{k,n}(\lambda) / \Cyl_{k,n}(\mu)$.)
\end{defin}

We will illustrate an example of such a proof. The proposition below is not a new result and follows directly from \cite[\S6, Corollary 6.12]{SagStan}, given \cite[\S4.3, Equation (3)]{Ful}. The purpose of the proposition below is, instead, to demonstrate a proof technique that uses results concerning cylindric tableaux in order to prove results about skew tableaux.

\begin{prop}
Given two (non-cylindric) partitions $\alpha$ and $\beta$, we have
\begin{equation*}
\sum \limits_{\mu} s_{\alpha / \mu}(\mathbf{x}) s_{\beta / \mu}(\mathbf{y}) \cdot \sum \limits_{\gamma} s_\gamma(\mathbf{x}) s_\gamma(\mathbf{y}) = \sum \limits_{\lambda} s_{\lambda / \beta}(\mathbf{x}) s_{\lambda / \alpha}(\mathbf{y}).
\end{equation*}
\end{prop}

\begin{proof}
Given a power series $P$, let $\hg_d P$ be the degree-$d$ homogeneous component of $P$ (the part of the power series consisting of only degree-$d$ terms). We show that for all $d$,
\begin{equation*}
\hg_d \sum \limits_{\mu} s_{\alpha / \mu}(\mathbf{x}) s_{\beta / \mu}(\mathbf{y}) \cdot \sum \limits_{\gamma} s_\gamma(\mathbf{x}) s_\gamma(\mathbf{y}) = \hg_d \sum \limits_{\lambda} s_{\lambda / \beta}(\mathbf{x}) s_{\lambda / \alpha}(\mathbf{y}).
\end{equation*}
Consider any particular $d$. Choose $k$ and $n$ such that $k$ and $n - k$ are both very large compared to $d$ and the dimensions of $\alpha$ and $\beta$. Specifically, if $\alpha = (a_1,a_2,\dots,a_m)$ and $\beta = (b_1,b_2,\dots,b_p)$, then the following values of $k$ and $n$ should be large enough:
\begin{equation*}
k = \max(m,p) + 2d + 1; \qquad n = k + \max(a_1,b_1) + 2d + 1.
\end{equation*}

\begin{lemma}
Given two partitions $\alpha$ and $\beta$, let $\alpha' = \Cyl_{k,n}(\alpha)$ and $\beta' = \Cyl_{k,n}(\beta)$. Then
\begin{equation*}
\hg_d \sum \limits_{\mu} s_{\alpha / \mu}(\mathbf{x}) s_{\beta / \mu}(\mathbf{y}) \sum \limits_{\gamma} s_\gamma(\mathbf{x}) s_\gamma(\mathbf{y}) = \hg_d \sum \limits_{\substack{\mu' \in \Cylpar; \\ \mu' \subseteq \alpha'; \mu' \subseteq \beta'}} s_{\alpha' / \mu'}(\mathbf{x}) s_{\beta' / \mu'}(\mathbf{y}).
\end{equation*}
(Partitions on the left side above are regular; those on the right side are cylindric.)
\end{lemma}

\begin{proof}
Take any particular $\mu'$ that produce terms of degree $d$ on the right hand side of the above equation. Note that this puts a limit on the number of boxes in $\alpha' / \mu'$ (and $\beta' / \mu'$). This means that $\alpha' / \mu'$ looks something like this:

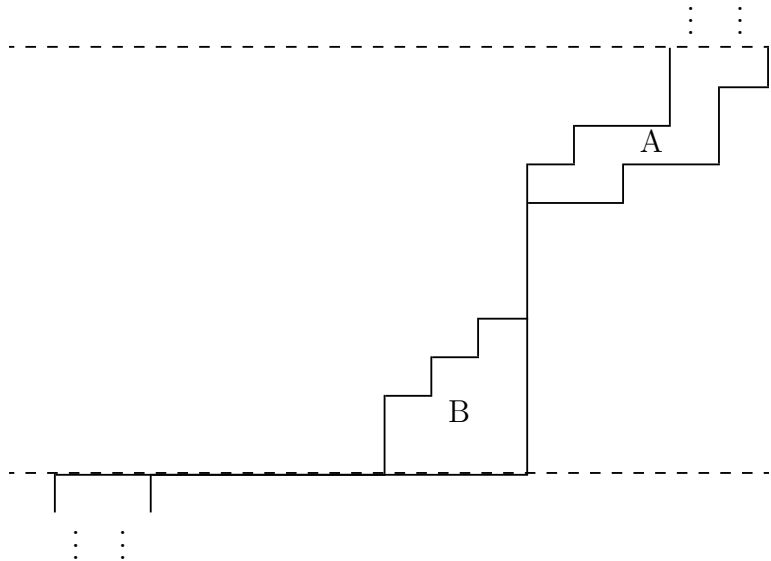
\begin{center}
\begin{tabular}{m{.2cm}m{.2cm}m{.2cm}m{.2cm}m{.2cm}m{.2cm}m{.2cm}m{.2cm}m{.2cm}m{.2cm}m{.2cm}m{.2cm}m{.2cm}m{.2cm}m{.2cm}m{.2cm}}
&&&&&&&&&&&&&&$\vdots$&$\vdots$\\
\hdashline&&&&&&&&&&&&&\multicolumn{1}{c|}{}&\textcolor{white}{1}&\multicolumn{1}{c|}{\textcolor{white}{1}}\\
\cline{16-16}&&&&&&&&&&&&&\multicolumn{1}{c|}{}&\multicolumn{1}{c|}{\textcolor{white}{1}}\\
\cline{13-14}&&&&&&&&&&&\multicolumn{1}{c|}{}&\textcolor{white}{1}&A&\multicolumn{1}{c|}{\textcolor{white}{1}}\\
\cline{12-12}\cline{14-15}&&&&&&&&&&\multicolumn{1}{c|}{}&\textcolor{white}{1}&\multicolumn{1}{c|}{\textcolor{white}{1}}\\
\cline{12-13}&&&&&&&&&&\multicolumn{1}{c|}{}\\
&&&&&&&&&&\multicolumn{1}{c|}{}\\
&&&&&&&&&&\multicolumn{1}{c|}{}\\
\cline{11-11}&&&&&&&&&\multicolumn{1}{c|}{}&\multicolumn{1}{c|}{\textcolor{white}{1}}\\
\cline{10-10}&&&&&&&&\multicolumn{1}{c|}{}&\textcolor{white}{1}&\multicolumn{1}{c|}{\textcolor{white}{1}}\\
\cline{9-9}&&&&&&&\multicolumn{1}{c|}{}&\textcolor{white}{1}&B&\multicolumn{1}{c|}{\textcolor{white}{1}}\\
&&&&&&&\multicolumn{1}{c|}{}&\textcolor{white}{1}&\textcolor{white}{1}&\multicolumn{1}{c|}{\textcolor{white}{1}}\\
\hdashline\cline{2-8}\cline{4-11}\multicolumn{1}{c|}{}&\textcolor{white}{1}&\multicolumn{1}{c|}{\textcolor{white}{1}}\\
&$\vdots$&$\vdots$
\end{tabular}
\end{center}

Because $k$ and $n - k$ are large compared to $d$ and the dimensions of $\alpha'$, the regions designated above as $A$ and $B$ are not connected.

Both regions $A$ and $B$ are in $\alpha'$; neither is in $\mu'$. Let $\nu'$ be the partition such that $\mu' \subseteq \nu' \subseteq \alpha'$, $B$ is in $\nu'$, and $A$ is entirely outside of $\nu'$. Since $A$ and $B$ are disconnected regions, we have
\begin{equation*}
s_{\alpha' / \mu'}(\mathbf{x}) = s_{\alpha' / \nu'}(\mathbf{x}) s_{\nu' / \mu'}(\mathbf{x}).
\end{equation*}
Similarly, we have
\begin{equation*}
s_{\beta' / \mu'}(\mathbf{x}) = s_{\beta' / \nu'}(\mathbf{x}) s_{\nu' / \mu'}(\mathbf{x}).
\end{equation*}

Since $\mu'$ can be thought of in terms of as $\alpha'$ sans region $A$, and then sans region $B$ (for small enough sizes of $A$ and $B$), we have
\begin{equation} \label{nuterms}
\hg_d \sum \limits_{\substack{\mu' \in \Cylpar; \\ \mu' \subseteq \alpha'; \mu' \subseteq \beta'}} s_{\alpha' / \mu'}(\mathbf{x}) s_{\beta' / \mu'}(\mathbf{y}) = \hg_d \sum \limits_{\nu'} s_{\alpha' / \nu'}(\mathbf{x}) s_{\beta' / \nu'}(\mathbf{y}) \sum \limits_{\mu'} s_{\nu' / \mu'}(\mathbf{x}) s_{\nu' / \mu'}(\mathbf{y}),
\end{equation}
where $\nu'$ ranges over all cylindric partitions such that $\Cyl_{k,n}(\nu) = \nu'$ for some partition $\nu$ such that $\alpha / \nu$ (a skew shape) contains $d$ or fewer boxes, and where $\mu'$ ranges over all cylindric partitions that are $\nu'$, except with a ``corner" (such as region $B$ above) --- which has at most $d$ boxes --- removed from $\nu'$.

For any particular $\nu'$ and $\mu'$ (satisfying the limitations described in the above paragraph), let $\ScrB$ be the skew shape that looks like the region $B$ between $\mu'$ and $\nu'$ (so in the above example, $\ScrB = (3,3,3,3)/(2,1)$). Clearly,
\begin{equation*}
s_{\nu' / \mu'}(\mathbf{x}) = s_\ScrB(\mathbf{x}).
\end{equation*}

Let $\ScrB'$ be the shape that is $\ScrB$ rotated by $180^\circ$ (because of the shape of the outer shape of $\ScrB$, we have that $\ScrB'$ is a straight (non-skew) partition; in our example, we would have $\ScrB' = (3,3,2,1)$). As it turns out, $s_\ScrB(\mathbf{x}) = s_{\ScrB'}(\mathbf{x})$ \cite[\S2, Exercise 2.22 (b)]{GriRei}. This is because if a skew tableau $R$ with shape $\ScrB$ has content $(a_1,a_2,\dots,a_m)$, we can create the tableau $R'$ with shape $\ScrB'$ that has content $(a_m,a_{m - 1},\dots,a_1)$ by switching all $1$'s with $m$'s, $2$'s with $m - 1$'s, etc., and by symmetry of Schur polynomials, we have that the number of tableaux with shape $\ScrB'$ and content $(a_m,a_{m - 1},\dots,a_1)$ equals the number of tableaux with shape $\ScrB'$ and content $(a_1,a_2,\dots,a_m)$. Therefore, we have that, given $\nu'$, for small enough degrees (including all up to $d$),
\begin{equation} \label{numugamma}
\sum \limits_{\mu'} s_{\nu' / \mu'}(\mathbf{x}) s_{\nu' / \mu'}(\mathbf{y}) = \sum \limits_{\ScrB'} s_{\ScrB'}(\mathbf{x}) s_{\ScrB'}(\mathbf{y})= \sum \limits_{\gamma} s_\gamma(\mathbf{x}) s_\gamma(\mathbf{y}).
\end{equation}
Combining equations \ref{nuterms} and \ref{numugamma}, we have
\begin{align*}
\hg_d \sum \limits_{\substack{\mu' \in \Cylpar; \\ \mu' \subseteq \alpha'; \mu' \subseteq \beta'}} s_{\alpha' / \mu'}(\mathbf{x}) s_{\beta' / \mu'}(\mathbf{y}) &= \hg_d \sum \limits_{\nu'} s_{\alpha' / \nu'}(\mathbf{x}) s_{\beta' / \nu'}(\mathbf{y}) \sum \limits_{\gamma} s_\gamma(\mathbf{x}) s_\gamma(\mathbf{y}) \\
&= \hg_d \sum \limits_{\nu'} s_{\alpha' / \nu'}(\mathbf{x}) s_{\beta' / \nu'}(\mathbf{y}) \sum \limits_{\gamma} s_\gamma(\mathbf{x}) s_\gamma(\mathbf{y}).
\end{align*}

Finally, since $s_{\alpha' / \nu'}(\mathbf{x})$, for every $\nu'$, describes the ways to fill the corresponding region $A$ (or its corresponding skew shape), and analogously for $\beta'$, we have that, for small enough degrees (including all up to $d$),
\begin{equation*}
\sum \limits_{\nu'} s_{\alpha' / \nu'}(\mathbf{x}) s_{\beta' / \nu'}(\mathbf{y}) = \sum \limits_{\mu} s_{\alpha / \mu}(\mathbf{x}) s_{\beta / \mu}(\mathbf{y}).
\end{equation*}
Therefore, 
\begin{equation*}
\hg_d \sum \limits_{\mu} s_{\alpha / \mu}(\mathbf{x}) s_{\beta / \mu}(\mathbf{y}) \sum \limits_{\gamma} s_\gamma(\mathbf{x}) s_\gamma(\mathbf{y}) = \hg_d \sum \limits_{\substack{\mu' \in \Cylpar; \\ \mu' \subseteq \alpha'; \mu' \subseteq \beta'}} s_{\alpha' / \mu'}(\mathbf{x}) s_{\beta' / \mu'}(\mathbf{y}),
\end{equation*}
as desired.
\end{proof}

By Theorem \ref{schurpols}, we have that
\begin{equation*}
\hg_d \sum \limits_{\substack{\mu' \in \Cylpar; \\ \mu' \subseteq \alpha'; \mu' \subseteq \beta'}} s_{\alpha' / \mu'}(\mathbf{x}) s_{\beta' / \mu'}(\mathbf{y}) = \hg_d \sum \limits_{\substack{\lambda' \in \Cylpar; \\ \alpha' \subseteq \lambda'; \beta' \subseteq \lambda'}} s_{\lambda' / \beta'}(\mathbf{x}) s_{\lambda' / \alpha'}(\mathbf{y}).
\end{equation*}
Thus, it suffices to show that
\begin{equation*}
\hg_d \sum \limits_{\substack{\lambda' \in \Cylpar; \\ \alpha' \subseteq \lambda'; \beta' \subseteq \lambda'}} s_{\lambda' / \beta'}(\mathbf{x}) s_{\lambda' / \alpha'}(\mathbf{y}) = \hg_d \sum \limits_{\lambda} s_{\lambda / \beta}(\mathbf{x}) s_{\lambda / \alpha}(\mathbf{y}).
\end{equation*}
his is true because the Schur polynomials on either side of the above equation (for a skew partition $\lambda$ and $\lambda' = \Cyl_{k,n}(\lambda)$) describe the ways of filling the same region with letters (since $k$ and $n - k$ are large enough compared to $d$ and the dimensions of $\alpha$ and $\beta$). Thus, we have

\begin{equation*}
\hg_d \sum \limits_{\mu} s_{\alpha / \mu}(\mathbf{x}) s_{\beta / \mu}(\mathbf{y}) \cdot \sum \limits_{\gamma} s_\gamma(\mathbf{x}) s_\gamma(\mathbf{y}) = \hg_d \sum \limits_{\lambda} s_{\lambda / \beta}(\mathbf{x}) s_{\lambda / \alpha}(\mathbf{y}).
\end{equation*}

This being true for all nonnegative integers $d$ (since we can always find large enough values of $k$ and $n$), we have proven our proposition.
\end{proof}

The proof style used above can be used to prove other facts about skew tableaux. This means that cylindric tableaux have the potential to be very useful to regular tableau theory.

\section{A Note on Knuth Equivalence for Cylindric Tableaux}
\subsection{Words and Knuth Equivalence}
This subsection introduces the concepts of tableau words, Knuth transformations, and Knuth equivalence. These concepts are well-established in tableau theory, but are described here in order to establish conventions. These conventions are the ones used in William Fulton's \textit{Young Tableaux} \cite[\S2.1]{Ful}. \emph{In this subsection, tableaux will refer to regular (non-cylindric) tableaux.}
\begin{defin}
The \textit{word} of a tableau is the sequence of letters (entries) in the tableau, reading left to right across the tableau's rows and then from bottom to top.
\end{defin}
For example, the following (skew) tableau's word is $3346354$.

\ytableausetup{boxsize=normal}
\begin{center}
\ytableaushort{\none\none\none\none4,\none\none\none35,3346}
\end{center}

\begin{defin} \label{kprime}
Let $yzx$ be a sequence of three consecutive letters of a word $w$ such that $x < y \le z$. Then a transformation of type $K'$ takes $w$ into the word $w'$ obtained by replacing the three letters $yzx$ with $yxz$.

For example, we have $3346354 \overset{K'}{ \rightarrow} 3343654$. Note that $K'$ may not be able to be applied to a word, or might be applicable to a word in more than one way.
\end{defin}

\begin{defin} \label{k2prime}
Let $xzy$ be a sequence of three consecutive letters of a word $w$ such that $x \le y < z$. Then a transformation of type $K''$ takes $w$ into the word $w'$ obtained by replacing the three letters $xzy$ with $zxy$.

For example, we have $3346354 \overset{K''}{ \rightarrow} 3346534$. Just like $K'$, the transformation $K''$ may not be able to be applied to a word, or might be applicable to a word in more than one way.
\end{defin}

\begin{defin}
An \textit{elementary Knuth transformation} is a transformation that is $K'$, $K''$, or the inverse transformation of $K'$ or $K''$. Two words are \textit{Knuth equivalent} if one can be obtained from another through a series of elementary Knuth transformations. Two tableaux are Knuth equivalent if their words are Knuth equivalent.
\end{defin}

Knuth equivalence is an equivalence relation; thus, for any two words $w$ and $v$, if $w$ is equivalent to $v$, then $v$ is equivalent to $w$. The significance of Knuth equivalence is that there is a unique semistandard tableau that is Knuth equivalent to a given skew semistandard tableau, and that tableau is that skew tableau's rectification \cite[\S2.1, Corollary 2]{Ful}. (See \cite[\S1.2]{Ful} for an explanation of tableau rectification.) In general, Knuth equivalence is an extraordinarily useful notion in tableau theory.

\subsection{Cyclic Knuth Equivalence}
The concept of tableau words makes sense for cylindric tableaux as well as regular tableaux.

\begin{defin}
The \textit{word} of a cylindric tableau is the sequence of letters (entries) in the tableau, reading left to right across the tableau's rows and then from row $k$ to row $1$.
\end{defin}

Due to the cyclic nature of cylindric tableaux, it is desirable that two tableaux be Knuth equivalent if one can be obtained from the other by shifting all cells by the same amount in the same direction. For example, the following two tableaux should be Knuth equivalent:

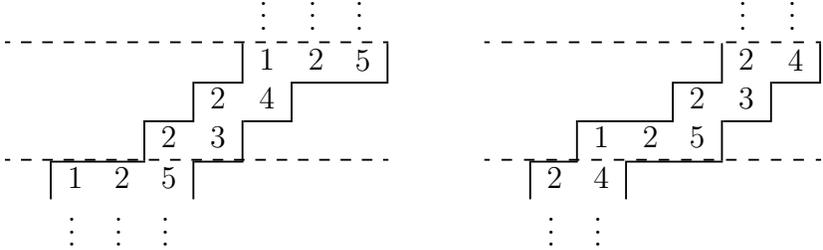
\begin{center}
\begin{tabular}{m{.2cm}m{.2cm}m{.2cm}m{.2cm}m{.2cm}m{.2cm}m{.2cm}m{.2cm}}
&&&&&$\vdots$&$\vdots$&$\vdots$\\
\hdashline&&&&\multicolumn{1}{c|}{}&1&2&\multicolumn{1}{c|}{5}\\
\cline{5-5}\cline{7-8}&&&\multicolumn{1}{c|}{}&2&\multicolumn{1}{c|}{4}\\
\cline{4-4}\cline{6-6}&&\multicolumn{1}{c|}{}&2&\multicolumn{1}{c|}{3}\\
\hdashline\cline{2-3}\cline{5-5}\multicolumn{1}{c|}{}&1&2&\multicolumn{1}{c|}{5}\\
&$\vdots$&$\vdots$&$\vdots$
\end{tabular}
\hspace{1cm}
\begin{tabular}{m{.2cm}m{.2cm}m{.2cm}m{.2cm}m{.2cm}m{.2cm}m{.2cm}}
&&&&&$\vdots$&$\vdots$\\
\hdashline&&&&\multicolumn{1}{c|}{}&2&\multicolumn{1}{c|}{4}\\
\cline{5-5}\cline{7-7}&&&\multicolumn{1}{c|}{}&2&\multicolumn{1}{c|}{3}\\
\cline{3-4}\cline{6-6}&\multicolumn{1}{c|}{}&1&2&\multicolumn{1}{c|}{5}\\
\hdashline\cline{2-2}\cline{4-5}\multicolumn{1}{c|}{}&2&\multicolumn{1}{c|}{4}\\
&$\vdots$&$\vdots$
\end{tabular}
\end{center}

Given our current definition of Knuth equivalence, this is not always the case. For example, the following two tableaux are not Knuth equivalent, as their words ($123$ and $312$, respectively) are not Knuth equivalent.

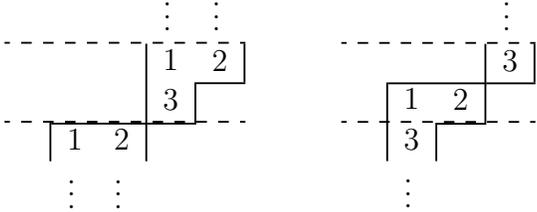
\begin{center}
\begin{tabular}{m{.2cm}m{.2cm}m{.2cm}m{.2cm}m{.2cm}}
&&&$\vdots$&$\vdots$\\
\hdashline&&\multicolumn{1}{c|}{}&1&\multicolumn{1}{c|}{2}\\
\cline{5-5}&&\multicolumn{1}{c|}{}&\multicolumn{1}{c|}{3}\\
\hdashline\cline{2-3}\cline{4-4}\multicolumn{1}{c|}{}&1&\multicolumn{1}{c|}{2}\\
&$\vdots$&$\vdots$
\end{tabular}
\hspace{1cm}
\begin{tabular}{m{.2cm}m{.2cm}m{.2cm}m{.2cm}}
&&&$\vdots$\\
\hdashline&&\multicolumn{1}{c|}{}&\multicolumn{1}{c|}{3}\\
\cline{2-3}\cline{4-4}\multicolumn{1}{c|}{}&1&\multicolumn{1}{c|}{2}\\
\hdashline\cline{3-3}\multicolumn{1}{c|}{}&\multicolumn{1}{c|}{3}\\
&$\vdots$
\end{tabular}
\end{center}

It is natural, then, to define a rotation operation $R$ on a word, which places the rightmost letter of a word to the left of the word. For example, we have $3346354 \overset{R}{ \rightarrow} 4334635$. By applying $R$ multiple times to a word, we can rotate the letters of the word by any amount.

\begin{defin}
An \textit{elementary cyclic Knuth transformation} is any transformation that is an elementary Knuth transformation or is $R$. Two words are \textit{cyclic Knuth equivalent} if one can be obtained from another through a series of elementary cyclic Knuth transformation. Two cylindric tableaux are cyclic Knuth equivalent if their words are cyclic Knuth equivalent.
\end{defin}

It turns out, however, that cyclic Knuth equivalence is not a useful equivalence relation for cylindric tableau theory, as all tableaux of the same content are cyclic Knuth equivalent. This fact is a consequence of Theorem 5.6.7 in Lothaire's \textit{Algebraic Combinatorics on Words} \cite[\S5.6, Theorem 5.6.7]{Loth}. Here, we provide a proof of this fact that relies exclusively on basic techniques.

\begin{theorem} \label{permequiv}
For any positive integer $m$, all words of length $m$ that are permutations (i.e. words of length $m$ that consist of the letters $1$ through $m$) are cyclic Knuth equivalent.
\end{theorem}

\begin{proof}
In this proof, ``word" will refer exclusively to permutations. We will also define a word as an equivalence class of words modulo $R$, since $R$ allows us to rotate any word to any desired position. We will identify each word with its representative that has $1$ as its leftmost letter. For any word $w$, we will say that $w_1 = 1$ and $w_i$ is the letter directly to the right of $w_{i-1}$ for $1 < i \le m$. In addition, $w_1$ is directly to the right of $w_m$.

We will create an algorithm that transforms any word $w$ into $x=1234 \dots m$. This is sufficient because, for any $w$ and $v$, if $w$ is equivalent to $x$ and $v$ is equivalent to $x$, then $w$ is equivalent to $v$.

Note that the elementary Knuth transformations state that we can switch two consecutive letters if the letter on either side of them is in between the two letters in value.
\begin{defin}
A \textit{catalyst} is a letter that facilitates the swapping of two other letters (i.e. the letter $y$ in definitions \ref{kprime} and \ref{k2prime}). We say that a catalyst \textit{catalyzes} the swapping of two letters.
\end{defin}

Our theorem may be easily verified for $m=1$, $m=2$, and $m=3$. For $m \ge 4$, our algorithm is as follows:

\begin{alg}[Word Transformation Algorithm] \label{wordtrans}
\begin{center}\end{center}
\vspace{.5cm}
\begin{algorithmic}[1]
Function WordTrans(word $w$) \Comment{$w$ must be a permutation of length $m$.}
\While{$w \neq x$}: \label{notdone}
\State $i := 1$.
\While{$i < m$}: \label{incri}
\If{$w_i$ can switch with $w_{i + 1}$ under an elementary Knuth transformation}:
\State Switch $w_i$ with $w_{i + 1}$.
\State Exit out of the loop beginning on line \ref{incri}.
\EndIf.
\State $i := i + 1$.
\EndWhile.
\EndWhile.
\State \Return $w$. \Comment{Clearly, $x$ is returned (if this line is reached).}
\end{algorithmic}
\end{alg}

\vspace{.4cm}
An example of this algorithm is shown below, with $w=159362847$. The two letters that are switched in the following step are bolded. Steps following those in which $w_1$ switches with $w_2$ are in red. The reader will gain intuition about the algorithm by following these steps; however, following all of these steps is not necessary for the proof.

\vspace{.4cm}
$15\mathbf{93}62847 \rightarrow \mathbf{15}3962847 \rightarrow \textcolor{red}{1\mathbf{39}628475} \rightarrow \mathbf{19}3628475 \rightarrow \textcolor{red}{13\mathbf{62}84759} \rightarrow \mathbf{13}2684759 \rightarrow \textcolor{red}{126\mathbf{84}7593} \rightarrow 1\mathbf{26}487593 \rightarrow \mathbf{16}2487593 \rightarrow \textcolor{red}{12\mathbf{48}75936} \rightarrow 1\mathbf{28}475936 \rightarrow \mathbf{18}2475936 \rightarrow \textcolor{red}{12\mathbf{47}59368} \rightarrow 1\mathbf{27}459368 \rightarrow \mathbf{17}2459368 \rightarrow \textcolor{red}{1245\mathbf{93}687} \rightarrow 124\mathbf{53}9687 \rightarrow 1\mathbf{24}359687 \rightarrow \mathbf{14}2359687 \rightarrow \textcolor{red}{123\mathbf{59}6874} \rightarrow 12\mathbf{39}56874 \rightarrow 1\mathbf{29}356874 \rightarrow \mathbf{19}2356874 \rightarrow \textcolor{red}{1235\mathbf{68}749} \rightarrow 123\mathbf{58}6749 \rightarrow 12\mathbf{38}56749 \rightarrow 1\mathbf{28}356749 \rightarrow \mathbf{18}2356749 \rightarrow \textcolor{red}{12356\mathbf{74}98} \rightarrow 1235\mathbf{64}798 \rightarrow 12\mathbf{35}46798 \rightarrow 1\mathbf{25}346798 \rightarrow \mathbf{15}2346798 \rightarrow \textcolor{red}{12346\mathbf{79}85} \rightarrow 1234\mathbf{69}785 \rightarrow 123\mathbf{49}6785 \rightarrow 12\mathbf{39}46785 \rightarrow 1\mathbf{29}346785 \rightarrow \mathbf{19}2346785 \rightarrow \textcolor{red}{123467\mathbf{85}9} \rightarrow 12346\mathbf{75}89 \rightarrow 123\mathbf{46}5789 \rightarrow 12\mathbf{36}45789 \rightarrow 1\mathbf{26}345789 \rightarrow \mathbf{16}2345789 \rightarrow \textcolor{red}{1234578\mathbf{96}} \rightarrow 123457\mathbf{86}9 \rightarrow 1234\mathbf{57}689 \rightarrow 123\mathbf{47}5689 \rightarrow 12\mathbf{37}45689 \rightarrow 1\mathbf{27}345689 \rightarrow \mathbf{17}2345689 \rightarrow \textcolor{red}{1234568\mathbf{97}} \rightarrow 12345\mathbf{68}79 \rightarrow 1234\mathbf{58}679 \rightarrow 123\mathbf{48}5679 \rightarrow 12\mathbf{38}45679 \rightarrow 1\mathbf{28}345679 \rightarrow \mathbf{18}2345679 \rightarrow \textcolor{red}{123456\mathbf{79}8} \rightarrow 12345\mathbf{69}78 \rightarrow 1234\mathbf{59}678 \rightarrow 123\mathbf{49}5678 \rightarrow 12\mathbf{39}45678 \rightarrow 1\mathbf{29}345678 \rightarrow \mathbf{19}2345678 \rightarrow \textcolor{red}{123456789}$

\vspace{.4cm}
\begin{defin}
For any word $w$, define $M(w)$ to be word $w_1 w_2 \dots w_j$ such that $w_1 < w_2 < \dots < w_j$, but $w_j > w_{j+1}$ (recall that $w_1 = 1 < w_2$). We will call $j$ the \textit{increase length} of $w$. (If $w$ has length $m$ and $w_1 < w_2 <  \dots < w_m$ (in other words, $w=x$), then $M(w) = w_1 w_2 \dots w_m$ and the increase length is $m$.)
\end{defin}

With this concept, we can illustrate exactly when the ``If" clause of the above algorithm is entered. It is possible that $w_1$ can be switched with $w_2$, with $w_n$ as a catalyst. If not, let $j$ be the increase length of $w$ (we assume that $j \neq m$, because in that case we are done). We know that $w_j > w_{j-1}$ and that $w_j > w_{j+1}$, so, out of $w_{j - 1}$, $w_j$, and $w_{j + 1}$, either $w_{j - 1}$ or $w_{j + 1}$ is the middle term in value. Thus, $w_j$ can be switched with the smaller of its two neighbors. The switch cannot occur earlier, because two consecutive terms in the middle of $w$'s increasing sequence do not have a term in between them in value to either side.

We conclude that the ``If" clause is entered in every pass through the loop beginning on line \ref{notdone}, because $M(w)$ ends with a letter that is not $w_m$ (since we are taking $w \neq x$).

One might observe from the example that the location of the switch always shifts to the left, unless $w_1$ switched with $w_2$ in the previous step. This is indeed the case.

\begin{lemma}
The location of the switch made by the algorithm always shifts to the left, unless $w_1$ switched with $w_2$ in the previous step.
\end{lemma}

\begin{proof}
Let $w$ be a word that is obtained during the execution of Algorithm \ref{wordtrans} (but not necessarily the word that is taken as input). We will consider two cases.

\textbf{Case 1: $w = 1cb \dots$, where $b$ and $c$ just switched (that is, $w_2$ and $w_3$ just switched values).} In this case, $c$ cannot be between $1$ and $b$ in value, because then $1$ and $b$ would have switched in the previous step. $c$ cannot be less than $1$. Thus, $1 < b < c$, which means that $b$ will catalyze the switch between $1$ and $c$.

\textbf{Case 2: $w = \dots abdc \dots$, where $c$ and $d$ just switched (that is, $w_l$ and $w_{l + 1}$ just switched for some $l > 2$).} In this case, $d$ cannot be between $b$ and $c$, because then $b$ and $c$ would have switched in the previous step. We also know that $a < b < c$, because at least one element of the increasing sequence beginning with $w_1$ is present in every switch (which implies that the increasing sequence continues up to, if not beyond, $c$). Thus, we have that either $a < b < c < d$, $a < d < b < c$, or $d < a < b < c$. In the first and third cases, $b$ switches with $d$ (which means that the location of the switch shifts one place to the left). In the second case, $a$ switches with $b$ (which means that the location of the switch shifts two places to the left).

Thus, the location of a switch is always to the left of the location of the previous switch, unless the previous switch was between $w_1$ and $w_2$.
\end{proof}

\begin{remark} \label{1or2}
Note that the switch in question always shifts either $1$ or $2$ places to the left (if the previous switch was not between $w_1$ and $w_2$) as the algorithm proceeds; this fact will be used later.
\end{remark}

From the fact that the switch always shifts to the left, it follows that, from any word, a switch between $w_1$ and $w_2$ will eventually happen (unless we arrive at $x$ beforehand). This means that, in performing the algorithm, an infinite loop in which $w_1$ never switches with $w_2$ will never be reached. Thus, if there exists an infinite loop in any possible execution of the algorithm, such a loop must contain a word that was obtained by switching $w_1$ with $w_2$.

\begin{defin}
A \textit{critical word} is a word obtained by switching $w_1$ with $w_2$. Note that no particular word is critical or non-critical; this adjective describes how the word was obtained within a particular execution of Algorithm \ref{wordtrans}.
\end{defin}

In the example above, the critical words were written in red. We will examine these critical words more closely.

\vspace{.4cm}
$ \dots \rightarrow \textcolor{red}{139628475} \rightarrow \dots \rightarrow \textcolor{red}{136284759} \rightarrow \dots \rightarrow \textcolor{red}{126847593} \rightarrow \dots \rightarrow \textcolor{red}{124875936} \rightarrow \dots \rightarrow \textcolor{red}{124759368} \rightarrow \dots \rightarrow \textcolor{red}{124593687} \rightarrow \dots \rightarrow \textcolor{red}{123596874} \rightarrow \dots \rightarrow \textcolor{red}{123568749} \rightarrow \dots \rightarrow \textcolor{red}{123567498} \rightarrow \dots \rightarrow \textcolor{red}{123467985} \rightarrow \dots \rightarrow \textcolor{red}{123467859} \rightarrow \dots \rightarrow \textcolor{red}{123457896} \rightarrow \dots \rightarrow \textcolor{red}{123456897} \rightarrow \dots \rightarrow \textcolor{red}{123456798} \rightarrow \dots \rightarrow \textcolor{red}{123456789}$

\vspace{.4cm}
For any $1 \le l \le m$, let $P_w(l)$ be the integer $t$ such that $w_t = l$. To each word $w$ of length $m$ we can assign a base-$m + 1$ number $N(w)$ whose base-$m + 1$ representation is $P_w(1)P_w(2) \dots P_w(m)$. For example, $N(159362847) = 164825973$. We will write down the numbers associated with the critical words in our example.

\vspace{.4cm}
$152794863 \rightarrow 142683759 \rightarrow 129573648 \rightarrow 128369547 \rightarrow 127358496 \rightarrow 126347985 \rightarrow 123946875 \rightarrow 123845769 \rightarrow 123745698 \rightarrow 123495687 \rightarrow 123485679 \rightarrow 123459678 \rightarrow 123456978 \rightarrow 123456798 \rightarrow 123456789$

\vspace{.4cm}
In our example, these numbers strictly decrease. It suffices to prove that they strictly decrease, because then a loop cannot be entered, and the smallest possible value of $N(w)$ for any word $w$, $123 \dots n$, will eventually be reached. It is indeed the case that these numbers are strictly decreasing.

\begin{lemma}
Let $w$ be a word and let $v$ be any critical word reached by following the algorithm starting with $w$, except for $123 \dots n$. Let $v'$ be the next critical word obtained by the algorithm, assuming that such a word exists.\footnote{If $v'$ does not exist, then $123 \dots n$ is reached, and so we are done.} Then $N(v') < N(v)$.
\end{lemma}

\begin{proof}
Let $l$ be the number such that $v_l \neq l$, but for all $i$ such that $1 \le i < l$, $v_i = i$. Then $N(v) = 123 \dots (l-1) \dots$. We have identified the position of the switch that happens to $v$ as either between $v_1$ and $v_2$, catalyzed by $v_m$, or between the local maximum\footnote{By this we mean a letter that is greater than the two letters to either side of it.} with the smallest position number and its smallest adjacent entry.

\textbf{Case 1: $v_1$ switches with $v_2$ with catalyst $v_m$.} Since $v_1 = 1$, $v_2$ cannot be $2$, as nothing can catalyze a switch between $1$ and $2$. Thus, $N(v)$'s second digit, $P_v(2)$, is some number that does not equal $1$ or $2$. When $v_1$ switches with $v_2$, the position of $2$ in the new word (call it $v'$) is smaller by $1$ (since $P_{v'}(v_2) = m$, and all other position numbers have decreased by $1$). Hence, in this case, $N(v') < N(v)$, as desired.

\textbf{Case 2: The switch that happens to $v$ is between the local maximum with the smallest position number in $v$ and the smallest entry adjacent to it.} The first local maximum cannot be $l$, as $v_l > l$ and $v_l$ is part of $M(v)$. Furthermore, since $l$ is smaller than some entry to its left, the first local maximum must be to the left of $l$. This means that the switch that the algorithm performs on $v$ decreases $P_v(l)$ by $1$ or keeps it the same. By our discussion above, the switch that happens to the word obtained from $v$ will be to the left of the switch that happens to $v$, so $P_v(l)$ will either decrease by $1$ or stay the same. This will continue until some letter switches with $1$, when $P_v(l)$ will decrease by $1$ (since the letter occupying the second position now occupies the $m$'th position). Note that the letter that switches with $1$ cannot be $l$, since, for this to happen, $l$ must have switched with $l - 1$ at some point, which is impossible (as there are no letters that could catalyze this switch).

Let $v'$ be the next critical word obtained from $v$. We have shown that $l$ will have a smaller position number in the next critical word; that is, $P_{v'}(l) < P_v(l)$. Now we will prove that for all $1 \le i < l$, $P_v(i) = P_{v'}(i)$.

By Remark \ref{1or2}, in the process in which $v$ is transforms into $v'$, the location of the switch shifts to the left by $1$ or $2$ places every time. Combined with the fact that the switch that happens to $v$ is to the right of $v_{l - 1}$ (or includes $v_{l - 1}$), we conclude that there is a switch that switches $v_{l - 1}$ with another letter during the process that transforms $v$ into $v'$. Consider the first such switch and let $u$ be the word to which this switch occurs. Note that for all $1 \le i < l$, $u_i = v_i = i$. Since the switch that happens to $u$ involves $u_{l - 1} = l - 1$, the switch cannot be between $u_{l - 1}$ and $u_{l - 2}$ (since, again, two letters whose value differs by $1$ cannot be switched). Hence, the switch occurs between $u_{l - 1}$ and $u_l$. Let $u'$ be the word that is obtained after this switch. Then $u' = 123 \dots (l - 2)u_l(l - 1)$. Since $l - 2 < l - 1 < u_l$, $u_l$ switches with $l - 2$. If the new word is called $u''$ (and $u''$ is not a critical word), we have $u'' = 123 \dots (l - 3)u_l(l - 2)$, so $u_l$ switches with $l - 3$. This continues until $u_l$ switches with $1$. When this happens, $v'$ is obtained, and by this analysis we know that $P_{v'}(i) = i$ for $1 \le i < l$.

Thus, we know that $N(v) = 123 \dots (l-1)P_v(l) \dots$ and $N(v')=123 \dots (l-1)P_{v'}(l) \dots$. As proved earlier, $P_{v'}(l) < P_v(l)$. Thus, $N(v') < N(v)$ in this case as well, as desired.
\end{proof}

We have constructed an algorithm that transforms words into other words. We have shown that this algorithm does not get stuck in a loop where no change is made. We have shown that this algorithm does not get stuck in a loop devoid of critical words. We have shown, through a decreasing monovariant, that no critical word can be repeated. Since the word with the smallest monovariant is $123 \dots m$, we have shown that this word is obtained through our algorithm. Since the Knuth transformations are reversible, if one can get from any word $p$ to another word $q$, one can get from $q$ to $p$ through these transformations. Thus, for any two words $w$ and $v$, we have shown how to get from $w$ to $123 \dots m$ and from $123 \dots m$ to $v$, thus constructing a method to get from $w$ to $v$. As this is applicable to all words, we have successfully shown that all permutations are cyclic Knuth equivalent.
\end{proof}

We now extend our theorem, which applies only to permutations, to all words, and thus show that any two cylindric tableaux with the same content are cyclic Knuth equivalent.

\begin{cor}
For any two words $w$ and $v$ such that $v$ is a permutation of $w$, $w$ and $v$ are cyclic Knuth equivalent.
\end{cor}

\begin{proof}
For any word $w$ of length $m$ that consists of more than one distinct letter (if all letters are the same, then clearly $w$ is cyclic Knuth equivalent to all of its permutations), let $s(w)$ be the smallest letter of $w$. Construct the word $w'$ as follows: let $t$ be the integer such that:

\begin{itemize}
\item If the last letter of $w$ is not $s(w)$, then $t$ is such that the $t$'th letter of $w$ (from the left) is the leftmost instance of $s(w)$ in $w$.
\item If the last letter of $w$ is $s(w)$, then $t$ is the smallest positive integer greater than $1$ such that the $t$'th letter of $w$ is $s(w)$, but the $(t-1)$'th letter of $w$ is not $s(w)$.
\end{itemize}

Then, for all $1 \le i \le m$, let the $i$'th letter of $w'$ be the $i$'th letter of $w$, plus $\frac{(i - t) \pmod m}{m}$. Next, construct the permutation $p(w)$ of length $m$ such that, for all $1 \le i,j \le m$, if the $i$'th letter of $w'$ is smaller than the $j$'th letter of $w'$, then the $i$'th letter of $p(w)$ is smaller than the $j$'th letter of $p(w)$. For instance, if $w = 43242$, then $s(w) = 2$, $t = 3$, $w' = (4.6)(3.8)(2)(4.2)(2.4)$, and $p(w) = 53142$.

We know that for any word $w$ of length $m$, $p(w)$ is transformed by Algorithm \ref{wordtrans} into $x = 12 \dots m$. Clearly, for any two distinct words $w$ and $v$ that are permutations of each other, $p(w) \neq p(v)$. We show $p(w)$ and $p(v)$ to be equivalent by applying Algorithm \ref{wordtrans} to each one, transforming them into $x$. Thus, if we can show that for each switch that occurs during the execution of Algorithm \ref{wordtrans} on $p(w)$, the corresponding switch of letters in $w$ is permissible under the cyclic Knuth transformations, then it will follow that $w$ and $v$ are cyclic Knuth equivalent. We will now show this.

For $1 \le i \le m$, define $w_i$ to be the letter corresponding to $p(w)_i$. (In our example above, the first $2$ is $w_1$.) For any $1 \le i,j \le m$, if $p(w)_i < p(w)_j$, then $w_i \le w_j$. Thus, any restrictions on switches in $w$ that do not exist in $p(w)$ arise only when two letters that are equal in $w$ are part of the switch (either being switched or acting as a catalyst).

The restrictions on switches involving equal letters set by definitions \ref{kprime} and \ref{k2prime} can be summarized as follows:

\begin{enumerate}[label=(\arabic*)]
\item Equal letters cannot be switched; \label{equalno}
\item A letter acting as a catalyst from the left cannot be equal to the smaller of the two letters being switched; and \label{leftcat}
\item A letter acting as a catalyst from the right cannot be equal to the larger of the two letters being switched. \label{rightcat}
\end{enumerate}

Suppose that all switches performed on $w$ corresponding to switches done by Algorithm \ref{wordtrans} on $p(w)$ have been legitimate up through the formation of a word $u$ ($u$ may equal $w$). (We think of $w$ and $p(w)$ as words that change value over time, as opposed to $u$, which is a particular word.) We will show (assuming that Algorithm \ref{wordtrans} has not yet terminated) that the switch that takes $u$ to a new word $u'$, analogous to the corresponding switch done by Algorithm \ref{wordtrans} with input $p(w)$, is also legitimate.

First, we note that, since all switches so far have been legitimate, it is the case that, for all $j$, the $j$'s in $w$ have stayed in the same order relative to each other. This means that the letters of $p(w)$ corresponding to the $j$'s in $w$ have also stayed in the same order relative to each other. In the original $p(w)$, this order is increasing (in our example above, the letters of $p(w)$ corresponding to the $4$'s in $w$ are $p(w)_2 = 4$ and $p(w)_4 = 5$). This means that no two of these letters can switch (because, if they are adjacent, they differ by $1$, which means that no letter can catalyze their switch). Consequently, condition \ref{equalno} is satisfied for the switch that takes $u$ to $u'$.

Suppose that condition \ref{leftcat} is violated by the switch that would take $u$ to $u'$. Due to the increasing nature of the letters of $p(w)$ corresponding to the $j$'s in $w$ (for any $j$) and the fact that $p(w)_1$ never switches with $p(w)_m$ during the execution of Algorithm \ref{wordtrans}, the only possibility of this is if $p(w)_m$ catalyzes the switch between $p(w)_1$ and $p(w)_2$. In order for this scenario to be a violation of condition \ref{leftcat}, $u_m$ must equal $s(w)$ (note that although $w$ changes, $s(w)$ is constant). We now show that this cannot be the case.

Suppose that $u_m = s(w)$. For the original value of $w$, $w_m \neq s(w)$ (because $t$ is originally chosen such that the $(t-1)$'th letter of $w$ --- that is, $w_m$ --- is not $s(w)$). Thus, $w_m$ becomes $s(w)$ after some switch that happens to a value of $w$ that precedes $u$. In order for this to happen, one of two things must happen: $p(w)_1$ must switch with $p(w)_2$, with $p(w)_2$ corresponding to $s(w)$ in $w$, --- this is impossible, as previously discussed --- or $p(w)_m$ must switch with $p(w)_{m - 1}$, with $p(w)_{m - 1}$ corresponding to $s(w)$ in $w$. This cannot be the case, because neither $p(w)_{m - 1}$ nor $p(w)_m$ can be the first local maximum in $p(w)$ (as we discussed in our proof of Theorem \ref{permequiv}, every switch not between $p(w)_1$ and $p(w)_2$ involves the first local maximum). Thus, $u_m \neq s(w)$ and condition \ref{leftcat} is not violated.

Suppose that condition \ref{rightcat} is violated by the switch that would take $u$ to $u'$. Due to the increasing nature of the letters of $p(w)$ corresponding to the $j$'s in $w$ (for any $j$) and the fact that $p(w)_1$ never switches with $p(w)_m$ during the execution of Algorithm \ref{wordtrans}, the only possibility of this is if $p(w)_1$ catalyzes a switch between $p(w)_{m -1}$ and $p(w)_m$; this is clearly impossible, as $p(w)_1$ cannot act as a catalyst.

Having proven that no condition is violated by the switch that takes $u$ to $u'$, we have completed our induction and have shown that $w$ and $v$ are cyclic Knuth equivalent for all $w$ and $v$ that are permutations of each other.
\end{proof}

It would be very helpful to have an analog to Knuth equivalence for cylindric tableaux; however, cyclic Knuth equivalence is clearly not the desired analog.

\newpage
\begin{center}
\textbf{\large{Acknowledgments}}
\end{center}

\vspace{.2cm}
First and foremost, I would like to thank my mentor, MIT graduate student Darij Grinberg, for introducing me to this topic, answering all of my questions, and helping to proofread my paper. I would also like to thank MIT professors Slava Gerovitch, Pavel Etingof, and Tanya Khovanova for organizing and managing PRIMES, the program that facilitated my research. Finally, I would like to thank MIT professor Alexander Postnikov for suggesting this subject of research.
\end{document}